\documentclass[11pt]{amsart}
\usepackage{amsmath,amssymb,graphicx,mathrsfs}
\usepackage{fullpage}
\usepackage{hyperref}
\usepackage{tikz}
\usepackage{color}
\def\al{\alpha}
\def\nn{\nonumber}
\def\nr{\nn\\}
\def\what#1#2#3{\widetilde{W}_{#1,#2}(#3)}

\def\nn{\nonumber}

\def\Z{{\mathbb Z}}
\def\C{{\mathbb C}}

\def\R{\mathbb R}
\def\eps{\varepsilon}

\def\pTR#1{p_{T,R}^{(#1)}}

\def\halpha{\widehat{\alpha}}
\def\hbeta{\widehat{\beta}}
\def\hn{\widehat{n}}
\DeclareMathOperator{\res}{Res}
\DeclareMathOperator{\len}{length}
\DeclareMathOperator{\stab}{Stab}
\def\rddots{\displaystyle\cdot^{\displaystyle\cdot^{\displaystyle\cdot}}}
\def\wlong{w_{\mathrm{long}}}
\DeclareMathOperator{\diag}{diag}
\DeclareMathOperator{\Ad}{Ad}
\newcommand\smallf[2]{{\textstyle{\frac{#1}{#2}}}}

\newtheorem{theorem}[equation]{Theorem}
\newtheorem{lemma}[equation]{Lemma}
\newtheorem{proposition}[equation]{Proposition}
\newtheorem{corollary}[equation]{Corollary}
\newtheorem{conjecture}[equation]{Conjecture}
\theoremstyle{remark}

\newtheorem{remark}[equation]{\bf Remark}
\newtheorem{definition}[equation]{\bf Definition}

\numberwithin{equation}{subsection}
\numberwithin{table}{section}

\usepackage{apptools}
\AtAppendix{
\numberwithin{equation}{section}
}

\newcommand{\FR}[1]{\mathcal{F}_R^{(#1)}}
\DeclareMathOperator{\GL}{GL}
\DeclareMathOperator{\SL}{SL}
\DeclareMathOperator{\re}{Re}
\DeclareMathOperator{\im}{Im}

\title{\boldmath An Asymptotic Orthogonality Relation for $\GL(n, \mathbb R)$}
\author{Dorian Goldfeld \and Eric Stade \and Michael Woodbury}
\address{Department of Mathematics\\Columbia University\\ 2990 Broadway \\ New York, NY 10027, USA}
\email{goldfeld@columbia.edu} 
\address{Department of Mathematics\\
University of Colorado, Boulder\\
Colorado 80309, USA}
\email{eric.stade@colorado.edu}
\address{Department of Mathematics \\Rutgers, The State University of New Jersey\\
110 Frelinghuysen Rd\\
Piscataway, NJ 08854-8019, USA}
\email{michael.woodbury@rutgers.edu}

\thanks{Dorian Goldfeld is partially supported by Simons Collaboration Grant Number 567168.}

\begin{document}

\begin{abstract}Orthogonality is a fundamental theme in representation theory and Fourier analysis. An orthogonality relation for characters of finite abelian groups (now recognized as an orthogonality relation on GL(1)) was used by Dirichlet to prove infinitely many primes in arithmetic progressions. Asymptotic orthogonality relations for GL$(n)$, with $n\le 3$, and applications to number theory, have been considered by various researchers over the last 45 years. Recently, the authors of the present work have derived an explicit asymptotic orthogonality relation, with a power savings error term, for GL$(4,\R)$.  Here we we extend those results to GL$(n,\mathbb  R)$ $(n\ge2)$.  

For $n\le 5$, our results are contingent on the Ramanujan conjecture at the infinite place, but otherwise are unconditional. In particular, the case $n=5$ represents a new result. The key new ingredient for the proof of the case $n=5$ is  the theorem of Kim-Shahidi that functorial products of cusp forms on GL(2)$\times$GL(3) are automorphic on GL(6). For $n>5$ (assuming again the Ramanujan conjecture holds at the infinite place), our results are conditional on two conjectures, both of which have been verified in various special cases.  The first of these conjectures regards lower bounds for Rankin-Selberg L-functions, and the second concerns recurrence relations for Mellin transforms of GL$(n,\mathbb R)$ Whittaker functions.

Central to our proof is an application of the Kuznetsov Trace formula, and a detailed analysis, utilizing a number of novel techniques, of the various entities---Hecke-Maass cusp forms, Langlands Eisenstein series, spherical principal series Whittaker functions and their Mellin transforms, and so on---that arise in this application.
\end{abstract}

\maketitle
\tableofcontents

\section{\large\bf Introduction} 
\subsection{Brief description of the main result of this paper}
Let $n\ge 1$ be a rational integer,  $s\in\mathbb C$, and  $\mathbb A_{\mathbb  Q} =\mathbb R\times \mathbb A_f$ denote the ring of adeles over $\mathbb Q$ where $\mathbb A_f$ denotes the finite adeles.
The family of unitary cuspidal automorphic representations  $\pi$ of $\text{\rm GL}(n,\mathbb A_{\mathbb  Q})$ and their standard L-functions 
$$L(s,\pi)= L_\infty(s,\pi)\cdot \prod_{p} L_p(s,\pi)$$ were first introduced by Godement and Jacquet \cite{MR0342495} and have played a major role in modern number theory. In the special case of $n=1$ the Euler products $\prod_{p} L_p(s,\pi)$  are just Dirichlet L-functions.

 In this paper we focus on the unitary cuspidal automorphic representations of $\text{\rm GL}(n,\mathbb A_{\mathbb  Q})$ with trivial central character which are globally unramified. For $n\ge 2,$ these can be studied classically in terms of Hecke-Maass cusp forms on
$$\text{\rm SL}(n, \mathbb Z)\backslash \text{\rm GL}(n, \mathbb R)/\left(\text{\rm O}(n, \mathbb R)\cdot \mathbb R^\times\right)$$
where $$\mathfrak h^n := \text{\rm GL}(n, \mathbb R)/\left(\text{\rm O}(n, \mathbb R)\cdot \mathbb R^\times\right)$$ is a generalization of the classical upper half-plane. In fact $\mathfrak h^2 := 
\left\{\left(\begin{smallmatrix} y&x\\0&1\end{smallmatrix}\right)\big\vert \; y>0, \, x\in\mathbb R \right\}$
is isomorphic to the classical upper half-plane.
\vskip 2pt
For $n\ge 2$, Hecke-Maass cusp forms are smooth functions $\phi: \mathfrak h^n\to \mathbb C$
which are automorphic for $\text{\rm SL}(n, \mathbb Z)$ with moderate growth and which are joint eigenfunctions of the full ring of invariant differential operators on $\text{\rm GL}(n, \mathbb R)$ and are  also joint eigenfunctions of the Hecke operators. Such globally unramified Hecke-Maass forms can be classified in terms of Langlands parameters
which (assuming the cusp form is tempered) are  $n$ pure imaginary numbers $(\alpha_1, \alpha_2, \ldots,\alpha_n) \in\mathbb (i\cdot \mathbb R)^n$ that sum to zero.  Further, the Hecke-Maass cusp forms $\phi$ with Langlands parameters $(\alpha_1, \ldots, \alpha_{n})$ can be ordered in terms of their Laplace eigenvalues $\lambda_\Delta(\phi)$ given by  
$$\lambda_\Delta(\phi) = \frac{n^3-n}{24} - \frac{\alpha_1^2+\alpha_2^2+\cdots +\alpha_n^2}{2}$$ 
as proved by Stephen Miller \cite{Miller_2002}.

\vskip 5pt
Let $\phi$ be a Hecke-Maass cusp form for $\text{\rm SL}(n, \mathbb Z)$ for $n \ge 2$ and set $$\langle \phi,\phi\rangle\, := \hskip-5pt\int\limits_{\text{\rm SL}(n, \mathbb Z)\backslash\mathfrak h^n}\hskip-5pt \phi(g) \overline{\phi(g)} \, dg$$ to denote the Petersson norm of $\phi.$  The Hecke-Maass cusp forms form a Hilbert space over $\mathbb C$   with respect to the Petersson inner product.  

\vskip 10pt
\begin{definition}[\bf L-function of a Hecke-Maass cusp form] Let $\phi$ be a Hecke-Maass cusp form for $\text{\rm SL}(n,\mathbb Z)$. Then for $s\in\mathbb C$ with $\text{\rm Re}(s)$ sufficiently large we define the L-function
$L(s,\phi) :=  \sum\limits_{k=1}^\infty \lambda(k) k^{-s}$
where $\lambda(k)$ is the $k^{th}$ Hecke eigenvalue of $\phi.$
\end{definition}

\begin{definition}[\bf Asymptotic orthogonality relation for  $\text{\rm \bf GL}(n,\mathbb R)$]
\label{GLnOrthoDef} 
Let $\{\phi_j\}_{j=1,2,\ldots}$ (with associated Langlands parameters $\alpha^{(j)} = (\alpha_1^{(j)}, \alpha_2^{(j)}, \ldots,\alpha_n^{(j)})$) denote an orthogonal basis of Hecke-Maass cusp forms for $\SL(n, \Z)$   with L-function given by
$L(s, \phi_j) := \sum\limits_{k=1}^\infty \lambda_j(k) k^{-s}.$
 Fix positive integers $\ell,m$.  Then,  for $T\to\infty$, we have
$$
\lim_{T\to\infty} \, \frac{\sum\limits_{j=1}^\infty  \lambda_j(\ell) \, \overline{\lambda_j(m)} \,\frac{h_T\left(\alpha^{(j)}\right)   }{\mathcal L_j}}{ \sum\limits_{j=1}^\infty  \frac{h_T\left(\alpha^{(j)}\right)   }{\mathcal L_j}} = \begin{cases} 1 + o\left(1\right) & \text{if}\;\; \ell = m,\\
 \;\;o\left(1\right) &\text{if}\;\; \ell\ne m.
\end{cases}  
$$
where $\mathcal L_j = L(1, \Ad{ \phi_j})$ and   $h_T\left(\alpha^{(j)}\right)$ is a smooth function of the variables  $\alpha^{(j)}, T$ (for $T>0$)  with support on the Laplace eigenvalues  $\lambda_\Delta(\phi_j)$ where  $0<\lambda_\Delta(\phi_j)\ll T.$
\end{definition}

\begin{remark}[\bf Power savings error term]  The asymptotic orthogonality relation has a power savings error term if $o(1)$ can be replaced with $\mathcal O\left(T^{-\theta}   \right)$ for some fixed $\theta > 0.$
 The error terms $o(1), \mathcal O\left(T^{-\theta}\right)$ will generally depend on $L, M$. This type of asymptotic orthogonality relation was first conjectured by Fan Zhou \cite{Zhou2014}.
\end{remark}

\begin{remark}[\bf Normalization of Hecke-Maass cusp forms] The approach we take in proving asymptotic orthogonality relations for $\text{\rm GL}(n,\mathbb R)$ is the Kuznetsov trace formula presented in \S \ref{KuznetsovTraceFormula} where $\frac{\lambda_j(\ell)\overline{\lambda_j(m)}}{\langle\phi_j, \phi_j\rangle}$
(which are independent of the way the $\phi_j$ are normalized)   appears naturally on the spectral side of the trace formula leading to an asymptotic orthogonality relation of the form
\begin{equation} \label{PhiInnerProduct}
\lim_{T\to\infty} \, \frac{\sum\limits_{j=1}^\infty  \lambda_j(\ell) \, \overline{\lambda_j(m)} \,\frac{h_T\left(\alpha^{(j)}\right)   }{\langle \phi_j,\phi_j\rangle}}{ \sum\limits_{j=1}^\infty  \frac{h_T\left(\alpha^{(j)}\right)   }{\langle \phi_j,\phi_j\rangle}} = \begin{cases} 1 + o\left(1\right) & \text{if}\;\; \ell = m,\\
 \;\;o\left(1\right) &\text{if}\;\; \ell\ne m.
\end{cases}  
\end{equation}
 If we normalize  $\phi_j$ so that its first Fourier coefficient is equal to one then  it is shown in
 Proposition \ref{PropFirstCoeff} that
 $$\langle\phi_j,\phi_j\rangle = c_n L(1, \Ad{ \phi_j}) \prod\limits_{1\le i\ne k\le n}\Gamma\left(\frac{1+\alpha^{(j)}_i-\alpha^{(j)}_k}{2}  \right), \qquad\quad (c_n\ne0).$$ This allows us (with a modification of the test function $h_T$) to replace the inner product $\langle\phi_j,\phi_j\rangle$ appearing in (\ref{PhiInnerProduct}) with the adjoint L-function $\mathcal L_j$ as in
  Definition \ref{GLnOrthoDef}. The main reason for doing this is that there are much better techniques developed for bounding special values of L-functions as opposed to bounding inner products of cusp forms. So having  $\mathcal L_j^{-1}$ in the asymptotic orthogonality relation instead of  $\langle\phi_j,\phi_j\rangle^{-1}$ will allow us to obtain better error terms in applications.
\end{remark}

Orthogonality relations  as in Definition \ref{GLnOrthoDef} have a long history going back to Dirichlet (for the case of GL(1)) who introduced the orthogonality relation for Dirichlet characters to prove infinitely many primes in arithmetic progressions. Bruggeman \cite{Bruggeman1978} was the first to obtain an asymptotic orthogonality relation for GL(2) which he presented in  the form
$$\lim_{T\to\infty} \sum_{j=1}^\infty \frac{\lambda_{j}(\ell) \overline{\lambda_{j}(m)}\cdot 4\pi^2 e^{-\frac{\lambda_\Delta(\phi_j)}{T} }}{T \cosh\left(\pi \sqrt{\lambda_{\Delta}(\phi_j) -\tfrac14}\right) } = \begin{cases} 1 & \text{if}\; \ell=m,\\
0 & \text{if} \; \ell\ne m.\end{cases}$$
where  $\big\{\phi_j\big\}_{j=1,2,\ldots}$ goes over an orthogonal basis of Hecke-Maass cusp forms for $\text{\rm SL}(2,\mathbb Z)$. This is not quite in the form of Definition \ref{GLnOrthoDef} but it can be put into that form with some work. Other versions of GL(2) type orthogonality relations with important applications were obtained by Sarnak \cite{Sarnak1984}, and, for holomorphic Hecke modular forms, by Conrey-Duke-Farmer \cite{CDF1997} and J.P. Serre \cite{Serre1997}.

The first asymptotic orthogonality relations for GL(3) with power savings error term were proved independently by Blomer \cite{Blomer_2013} and Goldfeld-Kontorovich \cite{GK2013}
in 2013. In 2021 Goldfeld-Stade-Woodbury \cite{GSW21} were the first to obtain a power-saving asymptotic orthogonality relation as in Definition \ref{GLnOrthoDef} for GL(4).

\vskip 10pt

A major breakthrough was obtained by
Matz-Templier \cite{MR4324829} who unconditionally proved an asymptotic orthogonality relation  for ${\rm SL}(n,\mathbb Z)$, as in (\ref{PhiInnerProduct}), for a wide class of test functions for all $n\ge 2$ (with power savings)  but without the harmonic weights given by the inverse of the adjoint L-function at 1. Their results were further strengthened in 
Finis-Matz \cite{MR4297181}. 
The principal tool used to prove the asymptotic orthogonality relation in  
\cite{MR4324829} was the Arthur-Selberg trace formula, whereas our approach is the natural generalization of the earlier results \cite{Blomer_2013}, \cite{GK2013}, \cite{GSW21}, which were based on the Kuznetsov trace formula. 
\vskip 8pt
Blomer \cite{MR4611943} presented a very nice exposition comparing the Arthur-Selberg and Kuznetsov trace formulae which we now briefly summarize for the application to asymptotic orthogonality relations. 
\vskip 8pt
\hspace*{0.6cm}
\begin{minipage}{.85\textwidth}
 $\bullet$ The first key difference  between these trace formulae  is that the spectral side of the Kuznetsov trace formula has harmonic weights $\mathcal L_j^{-1}$ while the Arthur-Selberg trace formula does not have these harmonic weights. For $\GL(n)$ with $n>3$ it is not currently known how to remove these weights (see \cite{MR4073945} for how to remove the weights on $\GL(3)$). In  \cite{MR4611943} Blomer remarks that {\it ``for applications to L-functions involving period formulae it is often desirable to have an additional factor $1/L(1, {\rm Ad}\; \phi)$ in the cuspidal spectrum, but in other situations one may prefer a summation formula without an extra L-value.''}
\end{minipage}
\vskip 8pt
\hspace*{0.6cm}
\begin{minipage}{.85\textwidth}
 $\bullet$ The second major difference  between these trace formulae  is that the spectral side of the Kuznetsov trace formula does not contain residual spectrum while the Arthur-Selberg trace formula does. As pointed out by a referee the bulk of the work in Matz-Templier \cite{MR4324829} consists  in bounding the unipotent contribution on the geometric side of the Arthur trace formula
 so that it stays in line with the error term coming from the residual Eisenstein contribution  on the spectral side given by Lapid-Mueller \cite{MR2541128}. These residual Eisenstein series do not appear in the Kuznetsov trace formula which leads to a very strong conjectural error term in Theorem \ref{Main-Theorem}. In fact, the largest error term on the spectral side of the Kuznetsov trace formula arises from the tempered Eisenstein series coming from the maximal parabolic having $(n-1,1)$ Levi block decomposition. For explicit comparisons between our main theorem and the results of  \cite{MR4324829} see Remark \ref{RemarkWeights}.
\end{minipage}
\vskip 8pt
\hspace*{0.6cm}
\begin{minipage}{.85\textwidth}
 $\bullet$ There are certain applications of our results using the Kuznetsov trace formula approach that go beyond the results in \cite{MR4324829}, \cite{MR4297181}. Recall that $\lambda_j(p)$  denotes the $p^{th}$ Hecke eigenvalue of the Maass form $\phi_j$.  Fan's thesis concerns the so-called vertical Sato-Tate problem which is a conjecture about the distribution of $\lambda_j(p)$ where $p$ is fixed and $j$ varies. This problem was studied by Bruggeman \cite{Bruggeman1978} and Sarnak \cite{Sarnak1984} (for Maass forms), and Serre \cite{Serre1997} and Conrey-Duke-Farmer \cite{CDF1997} (for holomorphic forms), who showed by fixing $p$ and varying $j$, that $\lambda_j(p)$ is an equidistributed sequence with respect to the Plancherel measure which depends on $p$.  Strikingly, as observed by Fan Zhou (\cite{Zhou2014}), if we give each Hecke eigenvalue $\lambda_j(p)$ the weight
$\mathcal{L}_j^{-1}$, then the distribution involves the Sato-Tate measure which is independent of $p$.  Jana, in \cite{MR4297178}, generalized the results of Zhou, but he only obtained an asymptotic formula without a power savings error term. A problem for the future would be to combine Jana's approach with the methods of this paper.  Jana also obtains bounds toward Sarnak's density hypothesis using this strategy that are stronger than anything known using the Arthur-Selberg trace formula.
\end{minipage}
\vskip 10pt
The main aim of this paper is to explicitly work out an asymptotic orthogonality relation for $\SL(n,\mathbb Z)$ via the Kuznetsov trace formula for a special choice of test function $h_{T,R}^{(n)}$ whose form is that of a Gaussian times a fixed polynomial. We do not address applications in this paper and leave that to future research. See \cite{MR4611943} for various applications of the Arthur-Selberg and Kuznetsov trace formulae and how they compare. We also point out that the Kuznetsov trace formula was generalized by Jacquet and Lai \cite{MR0783512} who developed the relative trace formula which has had a wide following with new types of applications.
\vskip 10pt
See \S \ref{Main-Theorem} for the statement of our main theorem. The proof we give assumes the Ramanujan conjecture at $\infty$ but it is possible to prove a weaker result by dropping this assumption. Otherwise the proof is unconditional for $n\le 5$.  In particular, the case $n=5$ represents a complete, new result. For $n>5$, our result is conditional on  two conjectures.

\vskip 10pt\noindent
\subsection{Ishii-Stade Conjecture} The Ishii-Stade Conjecture (see \S\ref{sec:IS-Conj}) concerns the  normalized Mellin transform $\what{n}{\alpha}{s}$ of the  $\text{\rm GL}(n,\mathbb R)$ Whittaker function $W_{n,\alpha}(y)$ defined in  Definition \ref{def:JacWhittFunction}.  Here, $s=(s_1, s_2,\ldots, s_{n-1})\in \C^{n-1}$, and  $\al=(\al_1,\al_2,\ldots,\al_n)=\C^{n-1}$ satisfies $\sum\limits_{i=1}^n \al_i=0$.

Suppose integers $m$ and $\delta$, with $1\le m\le n-1$ and $\delta\ge0$, are given.  The Ishii-Stade Conjecture expresses $\what{n}{\alpha}{s}$ as a finite linear combination, with coefficients that are rational functions of the $s_j$'s and $\al_k$'s, of shifted Mellin transforms $$\what{n}{\alpha}{s+\Sigma},$$where $\Sigma\in (\Z_{\ge0})^{n-1}$ and the $m$th coordinate of $\Sigma$ is $\ge \delta$.  In other words, for such $\delta$ and $m$, the conjecture expresses the Mellin transform $\what{n}{\alpha}{s}$ in terms of shifts of this Mellin transform  by at least $\delta$ units to the right in  the variable $s_m$.

Much as recurrence relations of the form $$\Gamma(s)=[(s+\delta-1)(s+\delta-2)\cdots (s+1)s]^{-1}\Gamma(s+\delta)$$for Euler's Gamma function imply concrete results concerning analytic continuation, poles, and residues of that function, so will the Ishii-Stade conjecture allow us to obtain explicit information about the behavior of $\what{n}{a}{s}$ beyond its original, {\it a priori} domain of definition.  This explicit information will be crucial to the analysis of our test function $h_T$, and consequently, to our derivation of an an asymptotic orthonality relation as in  Definition \ref{GLnOrthoDef}.

We have been able to prove the  Ishii-Stade Conjecture for  $\text{\rm GL}(n, \mathbb R)$ with $2\le n\le 5$.  See \S \ref{sec:IS-Conj} below.

\subsection{Lower bound conjecture for Rankin-Selberg L-functions}
 
Fix $n \ge 2$. Let $n = n_1+\cdots +n_r$ be a partition of $n$ with $n_i\in\mathbb Z_{>0}, (i=1,\ldots,r)$. The second conjecture we require for the proof of the asymptotic orthogonality relation for $\text{\rm GL}(n, \mathbb R)$  is a conjecture on the lower bound for Rankin-Selberg L-functions $L(s, \phi_k\times \phi_{k'})$ on the line $\text{\rm Re}(s) = 1,$ where  $\phi_k,\; \phi_{k'}$ (for  $1\le k<k'\le r$) are Hecke-Maass cusp forms for $\text{\rm SL}(n_k,\mathbb Z)$, $\text{\rm SL}(n_{k'},\mathbb Z)$, respectively.  For a Hecke-Maass cusp form $\phi$ with Langlands parameters $(\alpha_1,\ldots,\alpha_n)$, let
\begin{equation}\label{eq:analyticcond}
 c(\phi) = (1+ \lvert \alpha_1 \rvert)(1+ \lvert \alpha_2 \rvert)\cdots (1+ \lvert \alpha_n \rvert)
\end{equation}
denote the analytic conductor of $\phi$ as defined by Iwaniec and Sarnak \cite{MR1826269}. 
\begin{conjecture}[\bf Lower bounds for Rankin-Selberg L-functions]\label{RankinSelbergBound} Let $\varepsilon >0$ be fixed. Then we have the lower bound
$$
\left|L(1+it, \phi_k\times \phi_{k'})\right| \gg_\varepsilon \big(c(\phi_k)\cdot c(\phi_{k'})\big)^{-\varepsilon}\big(|t|+2\big)^{- \eps}.
$$
\end{conjecture}
\begin{remark} \label{KimShahidiReference}  Conjecture \ref{RankinSelbergBound} follows from Langlands' conjecture that $\phi_k\times\phi_{k'}$ is automorphic for $\text{\rm SL}(n_k\cdot n_{k'},\;\mathbb Z)$.  This can be proved via the method of de la Val\'ee Poussin as in Sarnak \cite{Sarnak2004}. Interestingly, Sarnak's approach can be extended to prove Conjecture \ref{RankinSelbergBound} if $\phi_{k'} $  is the dual of $\phi_k$ (see \cite{GLi2018}, \cite{HF2019}). Stronger bounds can also be obtained if one assumes the  Lindel\"of or Riemann hypothesis for Rankin-Selberg L-functions.  

If $n_k=n_k'=2$ it was proved by Ramakrishnan \cite{MR1792292} that $\phi_k\times\phi_{k'}$ is automorphic for $\text{\rm SL}(4,\mathbb Z),$ thus proving the lower bound conjecture for $n\le 4.$ Further, for $n_k=2 $ and $n_k'=3$, it was proved by Kim and Shahidi \cite{MR1923967} that $\phi_k\times\phi_{k'}$ is automorphic for $\text{\rm SL}(6,\mathbb Z),$ thus proving the lower bound conjecture for $n\le 5$.
\end{remark}

\subsection{Constructing the test functions}

Fix an integer $n\ge 2$. We now construct two complex-valued test functions on the space of Langlands parameters $$\Big\{\alpha=\left(\alpha_1,  \ldots, \alpha_n\right)\in\mathbb C^n \; \Big|\; \alpha_1+\cdots +\alpha_n=0\Big\}$$
that will be used in our proof of the orthogonality relation for $\text{\rm GL}(n,\mathbb R).$

We begin by introducing an auxiliary polynomial that is used in constructing the test functions.
\begin{definition}[\bf The polynomial $\mathcal F_R^{(n)}(\alpha)$]\label{def:FRn}
Let $R\in\mathbb Z_{>0}$ and let $\alpha = (\alpha_1, \ldots, \alpha_n)$ be a Langlands parameter. Then we define
$$\mathcal F_R^{(n)}(\alpha) := \prod_{j=1}^{n-2} \;\underset{\#K=\#L=j}{\prod_{K,L\,\subseteq\,(1,2,\ldots,n)}}\left(1+\sum_{k\in K}\alpha_k - \sum_{\ell\in L}\alpha_\ell     \right)^{\frac{R}{2}}.$$
\end{definition}
Note that if $\alpha\in (i\mathbb R)^n,$ then $\mathcal F_1^{(n)}(\alpha)$ is  {the square root of} a polynomial in $\alpha$ of degree  {$2D(n)$, where}
\begin{equation}\label{eq:Dn-explicit}
D(n) = \sum_{j=1}^{n-2} \frac12 \begin{pmatrix}  n\\ j \end{pmatrix} \left(\begin{pmatrix}  n\\ j \end{pmatrix}-1  \right) = \frac12 \begin{pmatrix}  2n\\ n \end{pmatrix} - \frac{n(n-1)}{2} - 2^{n-1}.
\end{equation}
 {By abuse of notation, we refer to $\mathcal{F}_R^{(n)}$ as a \emph{polynomial} although this is not strictly the case unless $R$ is even.}  For $\alpha$ with bounded real and imaginary parts, say, $|\text{\rm Re}(\alpha_j)|<R$ and $|\text{\rm Im}(\alpha_j)|<T^{1+\varepsilon}$ we have
\begin{equation} \label{FR-Bound}
\left|\mathcal F_R^{(n)}(\alpha)\right| \ll T^{\varepsilon+R\cdot D(n)}, \quad\qquad (T\to+\infty)
\end{equation}
with an implicit constant depending on $n, \varepsilon, R$.

\begin{definition}[\bf The test functions $p_{T,R}^{n,\#}(\alpha)$ and $h_{T,R}^{(n)}(\alpha)$]
\label{PTR-and-hTR}
 Let $R\in\mathbb Z_{>0}$ and $T\to+\infty$. Then for   a Langlands parameter $\alpha = (\alpha_1, \ldots, \alpha_n)$,  we define
 $$
 p_{T,R}^{n,\#}(\alpha) := e^{\frac{\alpha_1^2+\alpha_2^2+\cdots +\alpha_n^2}{2T^2}}\cdot \mathcal F_R^{\left(n \right)}\left(\tfrac{\alpha}{2}\right) \prod_{1\le j\ne k\le n} \Gamma\left(  \tfrac{1+2R+\alpha_j-\alpha_k}{4}  \right), 
 \quad 
  h_{T,R}^{(n)}(\alpha) := \frac{\left|p_{T,R}^{n,\#}(\alpha)\right|^2}{\prod\limits_{1\le j\ne k\le n} \Gamma\left(  \tfrac{1+\alpha_j-\alpha_k}{2}  \right)}.
$$
\end{definition}

We observe that, by Stirling's formula for the Gamma function and by \eqref{eq:Dn-explicit} and \eqref{FR-Bound}, we have\begin{equation}\label{hRgrowth}\big| h_{T,R}^{(n)}(\alpha)\big| \ll T^{R\cdot \left( \binom{2n}{n} -2^n\right)-\frac{n(n-1)}{2}}\end{equation}whenever $|\text{\rm Re}(\alpha_j)|$ is bounded and $|\text{\rm Im}(\alpha_j)|<T^{1+\varepsilon}$ for $1\le j\le n$.  The implied constant in \eqref{hRgrowth} depends on $n, \varepsilon,$ and $ R$.

 {
\begin{remark}[\bf Positivity of $h_{T,R}^{(n)}\,$]
Writing $\alpha=(\alpha_1,\alpha_2,\ldots,\alpha_n)$ with $\alpha_j=it_j$ and $t_j\in \R$ for each $j=1,2,\ldots,n$, the function $h_{T,R}^{(n)}(\alpha)$ is positive.  This is the case because $\Gamma(\frac{1+iu}{2})\Gamma(\frac{1-iu}{2}) = \lvert \Gamma(\frac{1+iu}{2})\rvert^2$ for $u\in \R$.
\end{remark}}

\begin{remark}[\bf Whittaker transform of the test function]
\label{SharpSymbol} The symbol $\#$ in the test function $
 p_{T,R}^{n,\#}$ means this function is the Whittaker transform of $
 p_{T,R}^{(n)}$. See \S~\ref{sec:IntRepOfpTR}. \end{remark}

\vskip 10pt
\subsection{The Main Theorem}

\begin{theorem} \label{Main-Theorem}
Fix $n\ge 2.$ Let $\{\phi_j\}_{j=1,2,\ldots}$ denote an orthogonal basis of Hecke-Maass cusp forms for $\SL(n, \Z)$ (assumed to be tempered at $\infty$)
  with associated Langlands parameter 
  \[ \alpha^{(j)} = \big(\alpha^{(j)}_1, \alpha^{(j)}_2,\ldots, \alpha^{(j)}_n\big)\in \left(i\cdot\mathbb R\right)^n \] 
 and L-function $L(s,\phi_j) :=  \sum\limits_{k=1}^\infty \lambda_j(k)\, k^{-s}$.  
 \vskip 5pt
 Fix positive integers $\ell, m$.  Then assuming the Ishii-Stade conjecture \ref{conj} and the lower bound conjecture for Rankin-Selberg L-functions  \ref{RankinSelbergBound}, we prove that for $T\to\infty$, 
\begin{multline*}
\sum\limits_{j=1}^\infty  \lambda_j(\ell) \, \overline{\lambda_j(m)} \,\frac{h_{T,R}^{(n)}\left(\alpha^{(j)}   \right)}{\mathcal L_j} 
  =    \delta_{\ell,m}\cdot 
\sum_{i=1}^{n-1} \mathfrak c_i \cdot T^{R\cdot \left( \binom{2n}{n} -2^n\right) + n-i}  
+\; \mathcal{O}_{\varepsilon,R,n}\left(  (\ell m)^{\frac{n^2+13}{4}}  \cdot T^{R\cdot \left(\binom{2n}{n} -2^n\right)\, +\,\varepsilon} \right)
\end{multline*}
where  $\delta_{\ell,m}$ is the Kronecker symbol,  $\mathcal L_j = L(1, \Ad\, \phi_j),$   and  $\mathfrak c_1,\ldots,\mathfrak c_{n-1}>0$ are absolute constants which depend at most on $R$ and $n$.

Because Conjectures   \ref{RankinSelbergBound} and \ref{conj} are known to be true for $2\le n\le5$ (see Remark \ref{KimShahidiReference} and \S \ref{sec:IS-Conj}), the above result is unconditional for such $n$.
\end{theorem} 
\begin{remark} Qiao Zhang \cite{2022arXiv221105747Z} recently proved the lower bound
\begin{equation} \label{QiaoZhang}
\left|L(1+it, \phi_k\times \phi_{k'})\right| \gg \big(c(\phi_k)\cdot c(\phi_{k'})\big)^{-\theta_{k,k'}}\big(|t|+2\big)^{-\frac{n_kn_{k'}}{2}\big( 1 - \frac{1}{n_k+n_{k'}}\big) - \eps}
\end{equation}
with $\theta_{k,k'} = n_k+n_{k'} +\varepsilon$. This improves on Brumley's bound who obtained nearly the same result but with the term $\frac{n_k n_{k'}}{2}$ replaced by $n_k n_{k'}$ (see \cite{Brumley2006} and the appendix of \cite{Lapid2013}). Assuming (\ref{QiaoZhang}) we can replace the error term in Theorem \ref{Main-Theorem} with
\[ \mathcal{O}_{\varepsilon,R,n,\ell,m}\Big(T^{R\cdot\left( \binom{2n}{n}-2^n \right) +n-1 \;+ \;\frac{n(n-2)}{6}\big(\theta_{k,k'} - \frac{8}{n^2}\big) } \Big).\] So if one could prove (\ref{QiaoZhang})  with $\theta_{k,k'}<\frac{8}{n^2}$ this would give a power savings error term in our main theorem and  would remove the assumption of the lower bound conjecture \ref{RankinSelbergBound}. In fact, the proof establishes a black box by which  improvements to bounds on Rankin-Selberg L-functions result in better power savings error terms for the continuous spectrum contribution to the asymptotic orthogonality relation.
\end{remark}

\begin{remark} \label{RemarkWeights} A variant of Theorem \ref{Main-Theorem} is obtained unconditionally in \cite{MR4324829}, \cite{MR4297181}, without the arithmetic weights $\mathcal L_j^{-1}$  and with different test functions, which are indicator functions of  $\alpha^{(j)} \in T\Omega$, where $\Omega$ is a Weyl group invariant bounded open subset of $i\cdot \mathfrak a^*$, where $\mathfrak a$ is the Lie algebra of the subgroup of diagonal matrices with positive entries. Additionally, the results of \cite{MR4324829}, \cite{MR4297181} do not entail the polynomial weights of size $T^{R\cdot \left( \binom{2n}{n} -2^n\right)-\frac{n(n-1)}{2}}$ coming from $h_{T,R}^{(n)}(\alpha)$ (cf. \eqref{hRgrowth}). 

The error term obtained in \cite{MR4297181},  in the present setting of $\SL(n,\Z)$, is $\ll T^{\frac{(n-1)(n+2)}{2}-1}$ as $T\to\infty$.  Here, $\frac{(n-1)(n+2)}{2}$ is the dimension of the generalized upper half-plane $\mathfrak h^n$, and the error term obtained by Finis-Matz has exponent equal to that dimension minus 1. By comparison, if one removes the polynomial weights $T^{R\cdot \left( \binom{2n}{n} -2^n\right)-\frac{n(n-1)}{2}}$ from the error term in Theorem \ref{Main-Theorem} above, then one obtains an error term that is $\ll T^{\frac{n(n-1)}{2}+\varepsilon}$.    
Also note that our main term is of a  stronger form than that of \cite{MR4324829}, \cite{MR4297181}, in that ours entails a sum of $n-1$  different high order asymptotics.

 More recently, Jana \cite{MR4297178} obtained a  proof of  the asymptotic orthogonality relation defined in \ref{GLnOrthoDef}, using the Kuznetsov trace formula and not the Selberg trace formula, with applications to the equidistribution of Satake parameters with respect to the Sato-Tate measure, second moment estimates of central values of L-functions as strong as Lindel\"of on average, and distribution of low lying zeros of automorphic L-functions in the analytic conductor aspect.  The paper of Jana does not contain a power saving error term. 
\end{remark}

\begin{remark} It is possible to remove the assumption 
of Ramanujan at the infinite place with more work which results in a weaker power savings error term in Theorem \ref{Main-Theorem}. For a  Maass form $\phi$ with Langlands parameter  $\alpha$, note that the test function $h_{T,R}(\alpha)$ is positive.  This is true because, even if $\alpha$ is a Langlands parameter   of an element in the complementary spectrum,  $-\alpha$ is a permutation of $\overline{\alpha}$.  A weaker version of Theorem \ref{Main-Theorem} can be proved if one assumes that almost all (except for a set of zero density) are tempered. Such results have been obtained in  \cite{MR4324829}, \cite{MR4297181}.
\end{remark}

\begin{proof}[ {Proof of Theorem~\ref{Main-Theorem}}]
Computing the inner product of certain Poincar{\'e} series in two ways (see the outline in \S\ref{sec:outline} below), we obtain a Kuznetsov trace formula relating the so-called geometric and spectral sides.  The geometric side consists of a main term $\mathcal{M}$ and a Kloosterman contribution $\mathcal{K}$.  The spectral side also consists of two components: a cuspidal (i.e., discrete) contribution $\mathcal{C}$ and an Eisenstein (i.e., continuous) contribution $\mathcal{E}$.

The left hand side of the theorem is precisely $\mathcal{C}$.  The first set of terms on the right hand side comes from the asymptotic formula for $\mathcal{M}$ given in Proposition~\ref{prop:mainterm}.  The power of $T$ in the error term comes from the bound for $\mathcal{E}$ given in  Theorem~\ref{th:EisensteinBound} (which also gives a factor of $\left( \ell m \right)^{\frac12-\frac{1}{n^2+1}}$).  A bound for $\mathcal{K}$, which is a (finite) sum of terms $\mathcal{I}_w$, with the same power of $T$ but with the given power of $\ell m$ follows as a consequence of Proposition~\ref{prop:Iwbounds}.
\end{proof}

\subsection{Outline of the key ideas in the proofs} \label{sec:outline}

Fix $n\ge 2.$ The  $\text{\rm GL}(n,\mathbb R)$ orthogonality relation appears directly in the spectral side of the Kuznetsov trace formula 
for $\text{\rm GL}(n,\mathbb R)$ which we now discuss. The Kuznetsov trace formula  is obtained by computing the inner product of two Poincar\'e series on $\text{\rm SL}(n, \mathbb Z)\backslash\mathfrak h^n$ in two different ways. The Poincar\'e series are constructed in a similar manner to Borel Eisenstein series by taking all $U_n(\mathbb Z)\backslash \text{\rm SL}(n, \mathbb Z)$ translates of a certain test function which  we choose to be the $p_{T,R}^{(n)}$ test function in Definition \ref{PTR-and-hTR}
 multiplied by a character and a  power function (see Definition \ref{def:PoincareSeries}).

The first way of computing the inner product of two Poincar\'e series is to replace one of the Poincar\'e series with its spectral expansion into cusp forms and Eisenstein series and then unravel the other Poincar\'e series with the Rankin-Selberg method. This gives the spectral contribution which has two parts: the cuspidal contribution and the Eisenstein contribution.  The second way of computing the inner product is to replace one of the Poincar\'e series with its Fourier Whittaker expansion and then unravel the other Poincar\'e series with the Rankin-Selberg method. This is called the geometric contribution to the trace formula, which also consists of two parts: a main term, and the so-called \emph{Kloosterman contribution}. The precise results of these computations are given in Theorems \ref{ThmSpectralDecomp} and \ref{GeomSide}, respectively. 

\vskip 24pt

\noindent
\underline{\bf Bounding the Eisenstein contribution}

\vskip 10pt
The key component of the Eisenstein contribution to the Kuznetsov trace formula 
is the inner product of an Eisenstein series
and the Poincar\'e series $P^M$ given in  Definition \ref{def:PoincareSeries}. By unraveling the Poincar\'e series in the inner product (see Proposition \ref{InnerProducts}) we essentially obtain     the $M^{th}$ Fourier coefficient of the Eisenstein series multiplied by the Whittaker transform of $p_{T,R}^{(n)}$. 
The explicit formula for the $M^{th}$ Fourier coefficient of the most general Langlands Eisenstein series given in Proposition  \ref{MthEisCoeff}  allows us to effectively bound all the terms in the integrals appearing in the Eisenstein contribution except for the product of adjoint L-functions
\begin{equation}
\label{ProductofAdjointLfunctions} \underset{n_k\ne1}{\prod_{k=1}^r} L^*\big(1, \text{\rm Ad}\; \phi_k\big)^{-\frac12 }
\end{equation}
  appearing in  that proposition. When considering the Eisenstein contribution to the Kuznetsov trace formula for $\text{\rm GL}(n,\mathbb R)$ all the adjoint L-functions in the above product are for cusp forms $\phi_k$ of lower rank $n_k <n$. Now in the special case that $\ell=m=1$, our Main Theorem  \ref{Main-Theorem} for $\text{\rm GL}(n,\mathbb R)$ gives a sharp bound for the sum of reciprocals of all adjoint L-functions of lower rank. This allows us to inductively prove  a power savings bound for the product (\ref{ProductofAdjointLfunctions}).

\vskip 24pt

\noindent
\underline{\bf Asymptotic formula for the geometric contribution}

\vskip 12pt

We prove that the geometric contribution is a sum of expressions $\mathcal I_w$ over elements $w$ in the Weyl group of $\text{\rm SL}(n,\mathbb Z)$. The  $\mathcal I_w$ are complicated multiple sums of multiple integrals weighted by Kloosterman sums (see (\ref{eq:IwDef})).
If $w_1$ is the trivial element of the Weyl group then we obtain an asymptotic formula for 
$\mathcal I_{w_1}$ (see Proposition~\ref{prop:mainterm})  while for all other Weyl group elements $\mathcal I_{w_i}$, with $i>1$, we obtain error terms with strong bounds for $|\mathcal  I_{w_i}|$ (see Proposition~\ref{prop:Iwbounds}) which are bounded by the final error term on the right side of our main theorem.

The key terms in \eqref{eq:IwDef}, the formula for $\mathcal I_w$, are the Kloosterman sums and two appearances of the test function $p_{T,R}^{(n)}$: one that is twisted by the Weyl group element $w$ and one that is not.  For the Kloosterman sums, we rely on bounds given by \cite{DR1998}.  The task of giving strong bounds for $p_{T,R}^{(n)}(y)$ occupies Sections~\ref{sec:IntRepOfpTR}, \ref{sec:ITRmBound} and \ref{sec:pTRnBound}.  We deal with the combinatorics of the twisted $p_{T,R}^{(n)}$-function, and we combine the bounds for it, the other $p_{T,R}^{(n)}$-function and the Kloosterman sums in Section~\ref{sec:GeometricSideBounds}.

The function $p_{T,R}^{(n)}$ is the inverse Whittaker transform of the test function $p_{T,R}^{n,\#}$ given in Definition~\ref{PTR-and-hTR} above.  Thanks to a formula of Goldfeld-Kontorovich \cite{MR2982424}, we can realize this as an integral of the product of $p_{T,R}^{n,\#}$, the Whittaker function $W_{\alpha}$ (see Definition~\ref{def:JacWhittFunction}), and certain additional  {gamma} factors.  We then write the Whittaker function as the inverse Mellin transform of its Mellin transform: $\what{n}{\al}{s}$.  This leads to the formula (valid for any $\varepsilon>0$):
\begin{multline*}
p_{T,R}^{(n)}(y) = \frac{1}{2^{n-1}} \int\limits_{\text{\rm Re}(\alpha_1)=0} \cdots
\int\limits_{\text{\rm Re}(\alpha_{n-1})=0}
e^{\frac{\alpha_1^2+\alpha_2^2+\cdots +\alpha_n^2}{T^2/2}} \, \mathcal F_R^{\left(n \right)}\left(\alpha \right)\prod_{1\le j\ne k\le n} \frac{\Gamma\left(  \tfrac{1+2R+\alpha_j-\alpha_k}{4}  \right)}{\Gamma\left( \frac{\alpha_j-\alpha_k}{2}  \right) }\\
\cdot \int\limits_{\text{Re}(s_1)=\varepsilon} \cdots \int\limits_{\text{Re}(s_{n-1})=\varepsilon} 
 \left(\,\prod_{j=1}^{n-1} y_j^{\frac{j(n-j)}{2}} (\pi y_j)^{-2s_j}\right)\widetilde{W}_{n, \alpha} \left(s\right) \, ds\, d\alpha.
\end{multline*}
To estimate the growth of $p_{T,R}^{(n)}(y)$ uniformly in $y$ and $T$ as $T\to+\infty$, we shift the line of integration in the $s$-integrals to $\text{\rm Re}(s) = -a$ with $a = (a_1,\ldots,a_{n-1})$ where $a_i>0$ for $i=1,\ldots,n-1$.  We remark that this is precisely where the Ishii-Stade Conjecture is required.  It is well known that
 \[ \what{2}{\alpha}{s} = \Gamma(s+\alpha)\Gamma(s-\alpha), \]
and hence understanding the values of $\what{2}{\alpha}{s}$ for $\re(s)<0$ is straightforward by applying the functional equation for the  {gamma} function or, equivalently, using an integral representation of the  {gamma} function valid for $\re(s)<0$.  A similar strategy can be used when $n=3$.  However, for $n\geq 4$,  {the analogous method seems intractable because the Mellin transform is not just a ratio of gamma functions, but an integral of such.  To overcome this difficulty, we apply}
the Ishii-Stade conjecture to describe the values of $\what{n}{\al}{s}$ in terms of sums of the Mellin transform of shifts of the $s$-variables.   {See also Remark \ref{rmk:mellincont} below.}

The Cauchy residue formula allows us to express $p_{T,R}^{(n)}$ as a sum of the shifted $s$-integral (termed the \emph{shifted $p_{T,R}^{(n)}$ term} and denoted $p_{T,R}^{(n)}(y;-a)$) and many residue terms.  The description of the shifted $p_{T,R}^{(n)}$ and residue terms is given in Section~\ref{sec:pTRnAsShiftPlusResidues}.  In order to bound $p_{T,R}^{(n)}(y;-a)$ it is convenient to introduce the function $\mathcal{I}_{T,R}(-a):=p_{T,R}^{(n)}(1;-a)$.

The next step is to use a result of Ishii-Stade (see Theorem~\ref{th:IS}) which allows us to write the Mellin transform $\what{n}{\al}{s}$ as an integral transformation of   $\what{n-1}{\beta}{z}$ against certain additional  {gamma} factors.  It is important to note that $\beta=(\beta_1,\ldots,\beta_{n-1})\in \C^{n-1}$ can be expressed in terms of $\alpha=(\alpha_1,\ldots,\alpha_n)$.  By carefullly teasing apart the portion of $\alpha$ which determines $\beta$ and that which doesn't, we are able to separate out the  {gamma} factors that don't depend on $\beta$ and bound $\mathcal{I}_{T,R}^{(n)}(-a)$ by the product of a power of $T$ and $\mathcal{I}_{T,R}^{(n-1)}(-b)$ for a certain $b=(b_1,\ldots,b_{n-1)}\in \R^{n-2}$.  This gives an inductive procedure, therefore, for bounding the shifted $p_{T,R}^{(n)}$ term.

In Section~\ref{sec:higherOrderResidues} we set notation for describing the $(r-1)$-fold shifted residue terms.  This requires generalizing a result of Stade (see Theorem~\ref{th:residue-description}) on the first set of residues of $\what{n}{\al}{s}$ (i.e., those that occur at $\re(s_i)=0$) to, first, higher order residues (i.e., taking the residue with respect to multiple values $s_i$), and second, to residues which occur along the lines $\re(s_i)=-k$ for $k\in \Z_{\geq 0}$.  This result, together with a teasing out of the variables similar to that described above, allows us to bound an $(r-1)$-fold residue term as the product of certain powers of $T$ and the variables $y_1,\ldots,y_{n-1}$ times 
 \[ \prod\limits_{j=1}^r \mathcal{I}_{T,R}^{(n_j)}(-a^{(j)}), \qquad\mbox{where $n=n_1+\cdots+n_r$.} \]
Applying the bounds on $\mathcal{I}_{T,R}^{(n_j)}$ that we inductively established for bounding the shifted $p_{T,R}^{(n)}$ term, and keeping careful track of all of the exponents and terms $a^{(j)}$, we eventually show that the bound for the shifted main term is in fact valid for every residue term as well.

\begin{remark}\label{rmk:GL3-GL4-compare}
In comparison to the results of \cite{GK2013} and \cite{GSW21}, we are using a slightly different normalization of the  {gamma} functions and the auxiliary polynomial $\mathcal{F}_R^{(n)}$ in the definition of the test functions $p_{T,R}^{n,\#}$ and $h_{T,R}^{(n)}$ (see Definition~\ref{PTR-and-hTR}).  Adjusting for this difference the results obtained here when applied to $n=3$ and $n=4$ recover the previously proven asymptotic formulas.
\end{remark}

\section{\bf\large Preliminaries}

\subsection{Notational conventions}

\begin{definition}[\bf{Hat notation for summation}] \label{sec:notation}
Suppose that $m\in \Z_+$ and $x=(x_1,\ldots,x_m)\in\C^{m}$.  For any $0\leq k \leq m$, define
 \[ \widehat{x}_k := x_1+\cdots + x_k. \]
Note that empty sums are assumed to be zero.
\end{definition}

\begin{definition}[\bf{Integration notation}]
Let $n\geq 2$.  We will often be working with $n$- and $(n-1)$-tuples of real or complex numbers.  We will denote such tuples without a subscript and use subscripts to refer to the components.  For example, we set $y=(y_1,\ldots,y_{n-1})\in \mathbb{R}_{>0}^{n-1}$, $s=(s_1,\ldots,s_{n-1})\in \C^{n-1}$ and $\alpha=(\alpha_1,\ldots,\alpha_n) \in \C^n$ such that
 \[ \alpha_1+\cdots+\alpha_n = 0. \]
In such cases, we denote integration over all such variables $x=(x_1,\ldots,x_k)$ subject to a condition(s) $\mathcal{C}=(\mathcal{C}_1,\ldots,\mathcal{C}_k)$ via
 \[ \int\limits_{\mathcal{C}} F(x) dx := \int\limits_{\mathcal{C}_1}\cdots \int\limits_{\mathcal{C}_k} F(x_1,\ldots,x_k)\; dx_1\, dx_2\cdots dx_k. \]
For example, given $\beta=(\beta_1,\ldots,\beta_{n-1})\in\C^{n-1}$ with $\hbeta_{n-1} = 0$, we denote integration over all such $\beta$ with $\re(\beta_j)=b_j$ for each $j=1,\ldots,n-2$ via
 \[ \int\limits_{\substack{\hbeta_{n-1}=0\\ \re(\beta)=b}} F(\beta)\, d\beta := 
 \int\limits_{\re(\beta_1)=b_1}
 \hskip -3pt \cdots \hskip -6pt
 \int\limits_{\re(\beta_{n-2})=b_{n-2}} F(\beta_1,\ldots,\beta_{n-2})\; d\beta_1\; d\beta_2\cdots d\beta_{n-2}. \]
We extend this notation liberally to integrals over $s$, $z$ and $\alpha$ and apply it also to integrals over the imaginary parts in the sequel.
\end{definition}

\begin{definition}[{\bf Polynomial notation}] \label{def:polynomialnotation}
Our analysis will often require us to bound certain polynomials in a trivial way.  Namely, for complex variables $x_j$ with $j=1,2,\ldots,k$, if $|x_j| \ll T^{1+\varepsilon}$ for each $j$ and $P(x_1,x_2,\ldots,x_k)$ is a polynomial, then $\lvert P(x_1,x_2,\ldots,x_k) \rvert \ll T^{\varepsilon+\deg{P}}$.  So, the relevant information about $P$ is its degree.  This being the case, we will use the notation $\mathcal{P}_d(x)$ (with $x=(x_1,\ldots,x_k)$) to represent an unspecified polynomial of degree less than or equal to $d$ in the variable(s) $x$.  Note that this notation agrees with the commonly employed practice (also used throughout these notes) of using $\varepsilon$ to represent an unspecified positive real number whose precise value is not specified and may differ from one usage to another.
\end{definition}

\begin{definition}[\bf{Vector or matrix notation depending on context}]
\label{VectorMatrix} Given a vector $a = (a_1,\ldots,a_{n-1})\in \mathbb R^{n-1}$, we shall  define the diagonal matrix 
 \[ t(a) := \mathrm{diag}(a_1a_2\cdots a_{n-1}, 
a_1 a_2\cdots a_{n-2} , \ldots, a_1,1).\] 
\end{definition}

\subsection{Structure of $\GL(n)$}

Suppose $n$ is a positive integer. Let $U_n(\R)\subseteq \GL(n,\R)$ denote the set of upper triangular unipotent matrices.

\begin{definition}[\bf{Character of $U_n(\R)$}]
Let $M=(m_1,\ldots,m_{n-1})\in \Z^{n-1}$.  For an element $x\in U_n(\R)$ of the form
\begin{equation}\label{eq:xmatrix-def}
 x = \left(\begin{smallmatrix} 
  1 & x_{1,2} & x_{1,3}& \cdots  & & x_{1,n}\\
  & 1& x_{2,3} &\cdots & & x_{2,n}\\
  & &\hskip 2pt \ddots & & & \vdots\\
  & && & 1& x_{n-1,n}\\
  & & & & &1\end{smallmatrix}\right),
\end{equation}
we define the character 
\begin{equation}\label{eq:psiMdef}
 \psi_M(x):= m_1x_{1,2}+m_2x_{2,3}+\cdots+m_{n-1} x_{n-1,n}.
\end{equation}
\end{definition}

\begin{definition}[\bf{Generalized upper half plane}]
We denote the set of (real) orthogonal matrices $\text{O}(n,\R)\subseteq \GL(n,\R)$, and we set 
 \[ \mathfrak h^n := \GL(n,\mathbb R)/\left(\text{O}(n,\R)\cdot\R^\times   \right). \]
Every element (via the Iwasawa decomposition of $\GL(n)$ \cite{Goldfeld2006}) of $\mathfrak h^n$ has a coset representative of the form $g=x y$, with $x$ as above and
\begin{equation}\label{eq:ymatrix-def}
 y =
    \left(\begin{smallmatrix} y_1y_2\cdots
    y_{n-1} & & & \\
    & \hskip -30pt y_1y_2\cdots y_{n-2} & & \\
    & \ddots &  & \\
    & & \hskip -5pt y_1 &\\
    & & &  1\end{smallmatrix}\right),
\end{equation}
where $y_i > 0$ for each $1 \le i \le n-1$.  The group $\GL(n,\R)$ acts as a group of transformations on $\mathfrak h^n$ by left multiplication.
\end{definition}

\begin{definition}[\bf{Weyl group and relevant elements}] \label{def:relevant}
Let $W_n \cong S_n$ denote the Weyl group of $\GL(n, \mathbb R).$  We consider it as the subgroup of $\GL(n,\R)$ consisting of permutation matrices, i.e., matrices that have exactly one $1$ in each row/column and all zeros otherwise.  An element $w\in W_m$ is called \emph{relevant} if
 \[ w = w_{(n_1,n_2,\ldots,n_r)} := \left( \begin{smallmatrix} & & I_{n_r} \\ & \rddots & \\ I_{n_1} & & \end{smallmatrix}\right),\]
where $I_{n_i}$ is the identity matrix of size $n_i\times n_i$ and $n=n_1+\cdots +n_r$ is a composition (a way of writing $n$ as a sum of  positive integers; see Section \ref{sec:pTRnAsShiftPlusResidues}).  The \emph{long element of $W_n$} is $\wlong:=w_{(1,1,\ldots,1)}$.
\end{definition}

\begin{definition}[{\bf Other subgroups of $\GL(n,\R)$}]
We define
 \[ \overline{U}_w := \Big(w^{-1}\cdot\hskip-8pt \phantom{U}^t U_n(\mathbb R)\cdot w\Big)\cap U_n(\R), \]
and
 \[ \Gamma_w := \Big(w^{-1}\cdot\hskip-8pt \phantom{U}^t U_n(\mathbb Z)\cdot w\Big)\cap U_n(\mathbb Z) = \SL(n,\Z) \cap \overline{U}_w, \]
where 
${}^t U_n$ denotes the transpose of $U_n$, i.e., the set of lower triangular unipotent matrices.
\end{definition}

\subsection{Basic functions on the generalized upper half plane $\mathfrak{h}^n$}

\begin{definition}[\bf{Power function}]
\label{def:Powerfunction}Let $\alpha=(\alpha_1,\ldots,\alpha_n)\in \C^n$ with $\halpha_n=0$.  Let $\rho=(\rho_1,\ldots,\rho_n)$, where $\rho_i=\frac{n+1}{2}-i$ for $i=1,2,\ldots,n$.   We define a power function on $xy\in \mathfrak h^n$  by
\begin{equation}\label{Ialpha}
I(xy,\alpha) =\prod_{i=1}^n d_i^{\alpha_i+\rho_i}= \prod_{i=1}^{n-1}y_i^{ {\halpha_{n-i}}+ {\widehat{\rho}_{n-i}}},
\end{equation}
where $d_i=\prod\limits_{j\le n-i}y_j$ is the $j$-th diagonal entry of the matrix $g=xy$ as above.
\end{definition}

\begin{definition}[\bf{Jacquet's Whittaker function}]\label{def:JacWhittFunction}
Let $g\in \GL(n,\R)$ with $n\geq 2$.  Let $\alpha = (\alpha_1,\alpha_2, \;\ldots, \;\alpha_n)\in \C^n$ with $\halpha_n=0$.  We define the completed    Whittaker function $W^{\pm}_{\alpha}: \GL(n,\mathbb R)\big/\left(
  \text{O}(n,\mathbb R)\cdot \mathbb R^\times\right) \to \mathbb C$   by the integral
$$W^{\pm}_{\alpha}(g) := \prod_{1\leq j< k \leq n} \frac{\Gamma\big(\frac{1 + \alpha_j - \alpha_k}{2}\big)}{\pi^{\frac{1+\alpha_j - \alpha_k}{2}} }\cdot \int\limits_{U_4(\mathbb R)} I(w_{\mathrm{long}} ug,\alpha) \,\overline{\psi_{1, \ldots, 1,\pm 1}(u)} \, du, $$
which converges absolutely if $\re(\al_i-\al_{i+1})>0$ for $1\le i\le n-1$ (cf.  \cite{GMW2021}), and has meromorphic continuation to all $\alpha\in \C^n$ satisfying $\halpha_n=0$.
\end{definition}

\begin{remark}
With the additional Gamma factors included in this definition (which can be considered as a ``\emph{completed}'' Whittaker function) there are $n!$ functional equations which is equivalent to the fact that the Whittaker function is invariant under all permutations of  $\alpha_1,\alpha_2, \ldots, \alpha_n$.  Moreover, even though the integral (without the normalizing factor) often vanishes identically as a function of $\al$, this normalization never does.

If $g$ is a diagonal matrix in $\GL(n,\mathbb R)$ then the value of $W^{\pm}_{n,\alpha}(g)$ is independent of sign, so we drop the $\pm$. We also drop the $\pm$ if the sign is $+1$.
\end{remark}

\begin{definition}[\bf{Whittaker transform and its inverse}]
Assume $n\ge 2$. Let $\alpha = (\alpha_1,\alpha_2, \ldots, \alpha_n)\in\C^n$ with $\halpha_n=0$.  Set $y := (y_1, y_2,\ldots y_{n-1})$ and $t(y)$ as in Definition~\ref{VectorMatrix}.  Let $f:\mathbb R_+^{n-1} \to \mathbb C$ be an integrable function. Then we define the Whittaker transform 
$f^\#: H^n \to \mathbb C$ (where $H^n:=\{ \alpha \in \C^n \mid \halpha_n=0\}$) by
 \begin{equation} \label{WhittakerTransform}
 f^\#(\alpha) := \int\limits_{y_1=0}^\infty \cdots  \int\limits_{y_{n-1}=0}^\infty   f(y)\,W_\alpha\big(t(y)\big) \prod_{k=1}^{n-1} \frac{dy_k}{y_k^{k(n-k)+1}  },
 \end{equation}
 provided the above integral converges absolutely and uniformly on compact subsets of $\mathbb R_+^{n-1}$. The inverse Whittaker transform  {\cite[Theorem~1.6]{MR2982424}} is
 \[ f(y) = \frac {1}{\pi^{n-1}} \int\limits_{\substack{ \halpha_n=0\\ \re(\al)=0}}  \frac{f^\#(\alpha) W_{-\alpha}\big(t(y)\big)}{\prod\limits_{1\le k\ne\ell\le n}\Gamma\left( \frac{\alpha_k-\alpha_\ell}{2}  \right)}\,d\al,\]
provided the above integral converges absolutely and uniformly on compact subsets of $(i\mathbb R)^n$.
\end{definition}

\begin{definition}[\bf{Normalized Poincar{\'e} series}]\label{def:PoincareSeries}
Let $M=(m_1,m_2,\ldots,m_{n-1})\in \Z^{n-1}$ with $m_i\neq 0$ for each $i=1,\ldots,n-1$.  As with $y$, we may think of $M$ as a matrix.  Let $g\in \mathfrak{h}^n$.  Then we define
\begin{equation}\label{eq:PoincareSeries-def}
 P^M(g,\alpha) := \frac{1}{\sqrt{\mathfrak{c}_n}}\cdot \prod_{k=1}^{n-1} m_k^{-\frac{k(n-k)}{2}}\sum_{\gamma \in U_n(\Z)\backslash \SL(n,\Z)} \psi_M(\gamma g) \cdot p_{T,R}^{(n)}(M\gamma g) \cdot I(\gamma g,\alpha),
\end{equation}
where $\mathfrak{c}_n$ is the (nonzero) constant given in Proposition~\ref{PropFirstCoeff}.  We extend the definition of $\psi_M$ and $p_{T,R}^{(n)}$ to all of $\mathfrak{h}^n$ by setting $\psi_M(xy):=\psi_M(x)$ and $p_{T,R}^{(n)}(xy):=p_{T,R}^{(n)}(y)$.
\end{definition}

\begin{remark}
This definition, up to the normalizing factor $\sqrt{\mathfrak c_n}\prod\limits_{k=1}^{n-1} m_k^{k(n-k)/2}$, of the Poincar{\'e} series agrees with that used in \cite{GSW21} with the minor caveat that $p_{T,R}$ takes on a slightly different normalization in terms of the polynomial $\mathcal{F}_R^{(n)}$ and in the Gamma factors appearing in Definition~\ref{PTR-and-hTR}.  The normalizing factor is inserted so that in the Kuznetsov trace formula the cuspidal term is precisely the orthogonality relation in Theorem~\ref{Main-Theorem}.
\end{remark}

\subsection{Fourier expansion of the Poincar\'e series}

\begin{definition}[\bf{Twisted Character}] Let 
$$V_n := \left\{v = \left.\left(\begin{smallmatrix}
v_1& & &\\
& v_2 & &\\
& & \ddots &\\
& & & v_n
\end{smallmatrix}\right) \right| \;v_1,\ldots,v_n\in\{\pm1\},  \; v_1\cdots v_n = 1\right\}.$$ Let $M = (m_1,\ldots,m_{n-1}) \in \mathbb Z^{n-1},$ and consider $\psi_M$ the additive character (see \eqref{eq:psiMdef}) of $U_n(\R)$. Then for $v \in V_n,$ we define the twisted character $\psi_M^v : U_n(\mathbb R) \to \mathbb C$ by 
$\psi_M^v(g) := \psi_M\left(v^{-1} g v\right).$
\end{definition}

\begin{definition}[\bf{Kloosterman Sum}] \label{KloostSum} Fix $L=(\ell_1,\ldots,\ell_{n-1}),\, M=(m_1,\ldots,m_{n-1})  \in\mathbb Z^{n-1}.$ Let $\psi_L, \psi_M$ be characters of $U_n(\mathbb R).$ Let $w \in W_n$ where $W_n$ is the Weyl group 
of $\GL(n).$ Let
\[ c = \left(\begin{smallmatrix} 1/c_{n-1} & & & & \\
& c_{n-1}/c_{n-2} & & & \\
& & \ddots & & \\
& & & c_2/c_1 & \\
& & & & c_1  \end{smallmatrix}\right) \] with $c_i\in  {\mathbb Z_{>0}}.$
 Then the Kloosterman sum is defined as
$$S_w(\psi_L, \psi_M, c) := \underset {\gamma = \beta_1cw\beta_2}{\sum\limits_{\gamma = U_n(\mathbb Z)\backslash \Gamma \cap G_w/\Gamma_w}} \psi_L(\beta_1) \,\psi_M(\beta_2),$$
with notation as in Definition~11.2.2 of \cite{Goldfeld2006}. The Kloosterman sum $S_w(\psi, \psi', c)$ is
well defined  {(i.e. independent of the choice of Bruhat decomposition for $\gamma$) if and only if} it satisfies the compatibility condition $\psi(cwuw^{-1}) = \psi'(u).$  
 {It is defined to be zero otherwise.  (See \cite{friedberg1987poincare}.)}
\end{definition}

\begin{proposition}[\bf{$M^{\mathrm{th}}$ Fourier coefficient of the Poincar\'e series $P^L$}]\label{FourierExp}
Let $L=(\ell_1,\ldots,\ell_{n-1})$ and $M=(m_1,\ldots,m_{n-1})  \in\mathbb Z^{n-1}$ satisfy $\prod\limits_{i=1}^{n-1} \ell_i\neq 0$ and $\prod\limits_{i=1}^{n-1} m_i\neq 0$.  If  {$\re(\alpha_k-\alpha_{k+1})$ is sufficiently large} for each $k=1,\ldots,n-1$, then
\begin{multline*}
\int\limits_{U_n(\mathbb Z)\backslash U_n(\mathbb R)} P^L\left(ug, \,\alpha \right) \cdot \overline{\psi_M(u)} \; d^*u  = \; \sum_{w\in W_n}\sum_{v\in V_n}\sum_{c_1=1}^\infty 
\cdots \sum_{c_{n-1}=1}^\infty  \frac{S_w(\psi_L,\psi_M^v,c) J_w(g; \alpha,\psi_L,\psi_M^v,c)}{ \sqrt{\mathfrak c_n} \prod\limits_{k=1}^{n-1}\left( \ell_k^{\frac{k(n-k)}{2}}c_k^{\alpha_k-\alpha_{k+1}+1} \right)  },
  \end{multline*}
where 
$$J_w(g; \alpha, \psi_L,\psi_M^v,c) = \int\limits_{U_w(\mathbb Z)\backslash {U}_w(\mathbb R)}\;
\int\limits_{\overline{U}_w(\mathbb R)} \psi_L(w u g)\, p_{T,R}^{(n)}\big(L c w u g\big)\, I(w ug,\alpha) \;\overline{\psi_M^v(u)} \; d^* u,$$
$${U}_w(\mathbb R) = \Big(w^{-1} \cdot  U_n(\mathbb R)\cdot w\Big) \cap 
U_n(\mathbb R), \qquad \overline{U}_w(\mathbb R) = \Big(w^{-1} \cdot \hskip-10pt\phantom{U}^t U_n(\mathbb R)\cdot w\Big) \cap U_n(\mathbb R),$$
and \hskip-8pt$\phantom{m}^t m$ denotes the transpose of a matrix $m$.
\end{proposition}
\begin{proof} 
See Theorem~11.5.4 of \cite{Goldfeld2006}.
\end{proof}

\section{\bf\large Spectral decomposition of  \boldmath $\mathcal{L}^2( \SL(n,\Z)\backslash \mathfrak{h}^n)$}

\subsection{Hecke-Maass cusp forms for $\SL(n,\Z)$}

\begin{definition}[\bf{Langlands parameters}]
Let $n\geq 2$.  A vector $\alpha=(\alpha_1,\ldots,\alpha_n)\in \C^n$ is termed a Langlands parameter if $\halpha_n=0$.
\end{definition}

\begin{definition}[\bf{Hecke-Maass cusp forms}]
 Fix $n\ge 2.$ A Hecke-Maass cusp form with Langlands parameter $\alpha\in \C^n$ for $\SL(n,\mathbb Z)$ is a smooth function $\phi: \mathfrak{h}^n \to \mathbb C$ which satisfies
 $\phi(\gamma g) = \phi(g)$
 for all $\gamma\in \SL(n,\Z)$, $g\in \mathfrak{h}^n$.  In addition $\phi$ is square integrable, is an eigenfunction of the algebra of Hecke operators on $\mathfrak{h}^n$, and is an eigenfunction of the algebra of $\GL(n,\R)$ invariant differential operators on $\mathfrak{h}^n$, with the same eigenvalues under this action as the power function $I(*,\alpha)$. The Laplace eigenvalue of $\phi$ is given by
   $$\frac{n^3 - n}{24} -\frac{\alpha_1^2+\alpha_2^2+\cdots+\alpha_n^2}{2}.$$
See  Section~6 in \cite{Miller_2002}. The Hecke-Maass cusp form $\phi$ is said to be tempered at $\infty$ if the Langlands parameters $\alpha_1,\ldots, \alpha_n$ are all pure imaginary.
\end{definition}

\begin{proposition}[\bf{Fourier expansion of Hecke-Maass cusp forms}] \label{FourierExpansion}
Assume $n\ge 2.$  Let 
$
\phi: \mathfrak h^n\to\mathbb C
$
be a Hecke-Maass cusp form for $\SL(n, \mathbb Z)$ with Langlands parameters $\alpha\in \C^n$. Then for $g\in \mathfrak{h}^n$,  we have the following Fourier-Whittaker expansion:

$$
\phi(g) = \sum_{\gamma\in{U}_{n-1}(\mathbb Z)\backslash \SL_{n-1}(\mathbb Z)} \;\sum_{m_1=1}^\infty\cdots  \sum_{m_{n-2}=1}^\infty\;\sum_{m_{n-1}\ne0} \, \frac{A_\phi(M)}
{\prod\limits_{k=1}^{n-1}  |m_k|^{\frac{k(n-k)}{2}}} \; W^{\text{\rm sgn}(m_{n-1})}_{\alpha}\left(t(M) \bigg(\begin{matrix} \gamma &0\\0&1\end{matrix}\bigg)g\right), $$
where $M = (m_1, m_2, \; \ldots, \; m_{n-1})$, $t(M)$ is the toric matrix in Definition~\ref{VectorMatrix}
and  $A_\phi(M)$ is  the $M^{th}$ Fourier coefficient of $\phi$.
\end{proposition}
\begin{proof}
See Section~9.1 of \cite{Goldfeld2015}. 
\end{proof}

\begin{definition}[{\bf L-function associated to a Hecke-Maass form $\pmb{\phi}$}]
Let $s\in\mathbb C$ with $\re(s)$  {sufficiently large}. Then the L-function  associated to a Hecke-Maass cusp form $\phi$ is defined as
$$L(s, \phi) := \sum_{m=1}^\infty \frac{A_\phi(m,1,\ldots,1)}{m^s}$$
and has holomorphic continuation to all $s\in \mathbb C$ and satisfies a functional equation $s\to 1-s.$
If $\phi$ is a simultaneous eigenfunction of all the Hecke operators then $L(s, \phi)$ has the following Euler product:
$$\begin{aligned} & L(s, \phi) = \prod_p \Bigg(1 - \frac{A(p,1,\ldots,1)}{p^{s}}
   + \frac{A(1,p,1,\ldots,1)}{p^{2s}}
    - \frac{A(1,1,p,\ldots,1)}{p^{3s}}\\
    &\hskip 170pt +\;\;\;\cdots  \;\;\; + \;(-1)^{n-1} \frac{A(1,,\ldots,1,p)}{p^{(n-1)s}} +\frac{(-1)^n}{p^{ns}}\Bigg)^{-1}.\end{aligned}$$
\end{definition}

\subsection{Langlands Eisenstein series for $\SL(n,\mathbb Z)$}

 \begin{definition}[\bf
 Parabolic Subgroup]\label{GLnParabolic} For $n\ge 2$ and $1\le r\le n,$ consider a partition of $n$ given by
$n = n_1+\cdots +n_r$ with positive integers $n_1,\cdots,n_r.$  We define the standard parabolic subgroup $$\mathcal P := \mathcal P_{n_1,n_2, \ldots,n_r} := \left\{\left(\begin{matrix} \GL(n_1) & * & \cdots &*\\
0 & \GL(n_2) & \cdots & *\\
\vdots & \vdots & \ddots & \vdots \\
0 & 0 &\cdots & \GL(n_r)\end{matrix}\right)\right\}.$$

Letting $I_r$ denote the $r\times r$ identity matrix, the subgroup
$$N^{\mathcal P} := \left\{\left(\begin{matrix} I_{n_1} & * & \cdots &*\\
0 & I_{n_2} & \cdots & *\\
\vdots & \vdots & \ddots & \vdots \\
0 & 0 &\cdots & I_{n_r}\end{matrix}\right)\right\}$$
is the unipotent radical of $\mathcal P$.  The subgroup
$$M^{\mathcal P} := \left\{\left(\begin{matrix} \GL(n_1) & 0 & \cdots &0\\
0 & \GL(n_2) & \cdots & 0\\
\vdots & \vdots & \ddots & \vdots \\
0 & 0 &\cdots & \GL(n_r)\end{matrix}\right)\right\}$$
 is the standard choice of Levi subgroup of $\mathcal P$.
\end{definition}

 \begin{definition}[\bf Hecke-Maass form $\Phi$ associated to a parabolic $\mathcal P$]\label{InducedCuspForm} \label{def:Hecke-Maassparabolic} Let $n\ge 2$. Consider a partition $n = n_1+\cdots +n_r$ with $1 < r < n$. Let  $\mathcal P := \mathcal P_{n_1,n_2, \ldots,n_r} \subset \GL(n,\mathbb R).$
  For $i = 1,2,\ldots, r$, let
$\phi_i:\GL(n_i,\mathbb R)\to\mathbb C$ be either the constant function 1 (if $n_i=1$) or a Hecke-Maass cusp form for $\SL(n_i,\mathbb Z)$ (if $n_i>1$).  The form $\Phi := \phi_1\otimes \cdots\otimes \phi_r$ is defined on $\GL(n,\mathbb R)=\mathcal P(\mathbb R)$  (where $K={\rm O}(n,\R)$) by the formula
$$\Phi(n m k ) := \prod_{i=1}^r \phi_i(m_i), \qquad (n\in  N^{\mathcal P}, m\in  M^{\mathcal P},k\in K)$$
where
  $m \in M^{\mathcal P}$ has the form
$m = \left(\begin{smallmatrix} m_1 & 0 & \cdots &0\\
0 & m_2 & \cdots & 0\\
\vdots & \vdots & \ddots & \vdots \\
0 & 0 &\cdots & m_r\end{smallmatrix}\right)$, with    $m_i\in \GL({n_i},\mathbb R).$  In fact, this construction works equally well if some or all of the $\phi_i$ are Eisenstein series.

 \end{definition}

\begin{definition}[\bf Character of a parabolic subgroup] \label{ParabolicNorm}  Let $n\ge 2.$ Fix a partition
$n = n_1+n_2 + \cdots +n_r$ with associated parabolic subgroup $\mathcal P := \mathcal P_{n_1,n_2, \ldots,n_r}.$  Define
\begin{equation}\label{rhoPj}
\rho_{_{\mathcal P}}(j)=\left\{
                            \begin{array}{ll}
                              \frac{n-n_1}{2}, & j=1 \\
                              \frac{n-n_j}{2}-n_1-\cdots-n_{j-1}, & j\ge 2.
                            \end{array}
                          \right.
\end{equation}
 Let $s = (s_1,s_2, \ldots, s_r) \in\mathbb C^r$ satisfy
$\sum\limits_{i=1}^r n_i s_i = 0.$
Consider the function  {(see Definition \ref{def:Powerfunction})}
$$| \cdot   |_{_{\mathcal P}}^s := I(\cdot,\alpha)$$
 on $\GL(n,\mathbb R)$, where
\begin{multline*}
\alpha=(\overbrace{s_1-\rho_{_{\mathcal P}}(1)+\smallf{1-n_1}{2}, \; s_1-\rho_{_{\mathcal P}}(1)+\smallf{3-n_1}{2}, \; \ldots \; ,s_1-\rho_{_{\mathcal P}}(1)+\smallf{n_1-1}{2}}^{n_1 \;\,\text{\rm terms}},\\
\overbrace{s_2-\rho_{_{\mathcal P}}(2)+\smallf{1-n_2}{2}, \; s_2-\rho_{_{\mathcal P}}(2)+\smallf{3-n_2}{2}, \; \ldots \; ,s_2-\rho_{_{\mathcal P}}(2)+\smallf{n_2-1}{2}}^{n_2 \;\,\text{\rm terms}},\\
\vdots
\\
\ldots \;\; ,\overbrace{s_r-\rho_{_{\mathcal P}}(r)+\smallf{1-n_r}{2}, \; s_r-\rho_{_{\mathcal P}}(r)+\smallf{3-n_r}{2}, \; \ldots \; ,s_r-\rho_{_{\mathcal P}}(r)+\smallf{n_r-1}{2}}^{n_r \;\,\text{\rm terms}}).
\end{multline*}
The conditions $\sum\limits_{i=1}^r n_i s_i = 0$ and $\sum\limits_{i=1}^r n_i \rho_{_{\mathcal P}}(i) = 0$
 guarantee that  $\alpha$'s entries sum to zero.  When $g\in \mathcal P$, with diagonal block entries $m_i\in \GL(n_i,\mathbb R)$, one has
 \vskip-10pt
 $$| g  |_{_{\mathcal P}}^s=\prod_{i=1}^r\left| \text{\rm det}(m_i)\right|^{s_i},$$
   so that $| \cdot   |_{_{\mathcal P}}^s$ restricts to a character of $\mathcal P$ which is trivial on $N^{\mathcal P}$.   
 \end{definition}

 \begin{definition}[\bf Langlands Eisenstein series twisted by Hecke-Maass forms of lower rank] \label{EisensteinSeries}
Let $\Gamma = \SL(n, \mathbb Z)$ with $n \ge 2.$  Consider a parabolic subgroup $\mathcal P = \mathcal P_{n_1,\ldots,n_r}$ of $\GL(n,\mathbb R)$ and functions $\Phi$  and $| \cdot |_{_{\mathcal P}}^s$   as given in Definitions \ref{InducedCuspForm} and  \ref{ParabolicNorm}, respectively. Let $$s = (s_1,s_2, \ldots,s_r)\in\mathbb C^r, \;\;\text{where}\;\; \sum_{i=1}^r n_is_i =0.$$
The Langlands Eisenstein series determined by this data is defined by
\begin{equation}\label{def:EPhi}
 E_{\mathcal P, \Phi}(g,s) := \sum_{\gamma\,\in\, (\mathcal P\, \cap \,\Gamma)\backslash \Gamma}  \Phi(\gamma g)\cdot |\gamma g|^{s+\rho_{_{\mathcal P}}}_{_{\mathcal P}}
 \end{equation}
as an absolutely convergent sum for $\re(s_i)$ sufficiently large, and extends to all $s\in \C^r$ by meromorphic continuation.
\end{definition}

For $k=1,2,\ldots,r,$ let $\alpha^{(k)} := (\alpha_{k,1},\ldots,\alpha_{k,n_k})$ denote the Langlands parameters of $\phi_k.$ We adopt the convention that if $n_k=1$ then $\alpha_{k,1} = 0.$
Then the Langlands parameters of $E_{\mathcal P, \Phi}(g,s)$ (denoted $\alpha_{_{\mathcal P,\Phi}}(s)$) are
\begin{multline}\label{langlandsparamsforPhi}
\bigg (\overbrace{\alpha_{1,1}+s_1, \;\ldots \;,\alpha_{1,n_1}+s_1}^{n_1 \;\,\text{\rm terms}}, \quad\overbrace{\alpha_{2,1}+s_2, \;\ldots \;,\alpha_{2,n_2}+s_2}^{n_2 \;\,\text{\rm terms}},
\quad\ldots \quad
,\overbrace{\alpha_{r,1}+s_r, \;\ldots \;,\alpha_{r,n_r}+s_r}^{n_r \;\,\text{\rm terms}}\bigg). \phantom{xxxxx}
\end{multline}

 {
\begin{definition} {\bf  \boldmath (The $M^{\rm{th}}$ Fourier coefficient of   $E_{\mathcal P, \Phi}$)}
 Let $s = (s_1,s_2, \ldots,s_r)\in\mathbb C^r,$  where $\sum\limits_{i=1}^r n_is_i =0.$
Consider  $E_{\mathcal P, \Phi}(*,s)$ with associated Langlands parameters $\alpha_{_{\mathcal P,\Phi}}(s)$ as defined in (\ref{langlandsparamsforPhi}).  Let $M = (m_1,m_2,\ldots, m_{n-1}) \in \mathbb Z_{>0}^{n-1}$.  Then the $M^{th}$ term in the Fourier-Whittaker expansion of $E_{\mathcal P, \Phi}$ is
\begin{align*}
\int\limits_{U_n(\mathbb Z)\backslash U_n(\mathbb R)} E_{\mathcal P, \Phi}(ug, s)\, \overline{\psi_M(u) }\; du \; = \; \frac{A_{E_{\mathcal P, \Phi}}(M,s)}
{\prod\limits_{k=1}^{n-1}  m_k^{k(n-k)/2}} \; W_{\alpha_{_{\mathcal P,\Phi}}(s)}\big(M g\big),
\end{align*}
\end{definition}
}

\subsection{Langlands spectral decomposition for $\SL(n,\mathbb Z)$}

\begin{definition}[\bf{Petersson inner product}]\label{def:PeterssonIP}
Let $n\geq 2$.  For $F,G\in \mathcal{L}^2( \SL(n,\Z)\backslash \mathfrak{h}^n)$ we define the Petersson inner product to be
 \[ \big\langle F,G \big\rangle := \int\limits_{ \SL(n,\Z)\backslash \mathfrak{h}^n} F(g)\overline{G(g)}\, dg. \]
For $g=xy\in \mathfrak{h}^n$, with
 \[ x = \left(\begin{smallmatrix} 
  1 & x_{1,2} & x_{1,3}& \cdots  & & x_{1,n}\\
  & 1& x_{2,3} &\cdots & & x_{2,n}\\
  & &\hskip 2pt \ddots & & & \vdots\\
  & && & 1& x_{n-1,n}\\
  & & & & &1\end{smallmatrix}\right), \qquad y =
    \left(\begin{smallmatrix} y_1y_2\cdots
    y_{n-1} & & & \\
    & \hskip -30pt y_1y_2\cdots y_{n-2} & & \\
    & \ddots &  & \\
    & & \hskip -5pt y_1 &\\
    & & &  1\end{smallmatrix}\right),\]
the measure $dg$ is given by $dx\, dy$, with 
 \[ dx = \prod_{1\leq i  {<} j \leq n} dx_{i,j}, \qquad dy = \prod_{k=1}^{n-1} \frac{dy_k}{y^{k(n-k)+1}}. \]
\end{definition}

The Langlands spectral decomposition for $\SL(n,\mathbb Z)$ states that
 \[ \boxed{\mathcal L^2(\SL(n,\mathbb Z)\backslash\mathfrak h^n) = \text{(Cuspidal spectrum)} \oplus \text{(Residual spectrum)}\oplus
\text{(Continuous spectrum)}}. \]
We shall be applying the Langlands spectral decomposition to Poincar\'e series which are orthogonal to the residual spectrum.

\begin{theorem}[\bf Langlands spectral decomposition for $\SL(n,\mathbb Z)$]
 \label{SpectralDecomp} Let $\phi_1, \phi_2, \phi_3, \ldots$\linebreak denote an orthogonal   basis of  Hecke-Maass forms for $\SL(n,\mathbb Z)$.   Assume that  $F, G \in \mathcal L^2(\SL(n,\mathbb Z)\backslash\mathfrak h^n)$ are orthogonal to the residual spectrum. Then for
$g\in \GL(n,\mathbb R)$ we have 
\begin{align*}
& F(g) = \sum_{j=1}^\infty \langle F, \phi_j\rangle \frac{\phi_j(g)}{\langle\phi_j, \phi_j\rangle} \; 
+ \; \sum_{\mathcal P } \sum_{\Phi}\; c_{\mathcal P} \hskip-10pt\underset{\text{\rm Re}(s)=0}{\int\limits_{n_1s_1+\cdots +n_rs_{r}=0}}   \hskip-14pt\Big\langle F, E_{\mathcal P,\Phi}(*\, ,s)\Big\rangle E_{\mathcal P,\Phi}(g\, ,s) \; ds,\nonumber
 \end{align*}
 
\begin{align*}
&\langle F, G\rangle =  \sum_{j=1}^\infty  \frac{\langle F, \phi_j\rangle \, \langle\phi_j, G\rangle}{\langle\phi_j, \phi_j\rangle} 
\; + \sum_{\mathcal P} \sum_{\Phi}\; c_{\mathcal P} \hskip-10pt \underset{\text{\rm Re}(s)=0}{\int\limits_{n_1s_1+\cdots +n_rs_{r}=0}} \hskip-10pt\Big\langle F, \, E_{\mathcal P,\Phi}(*\, ,s)\Big\rangle      \Big\langle E_{\mathcal P,\Phi}(*\, ,s), \,G \Big\rangle\; ds,
\nonumber
\nonumber\end{align*}
where the sum over $\mathcal P$ ranges over  parabolics  associated to partitions $n =  \sum\limits_{k=1}^r n_k$, while the sum over $\Phi$ (see Definition \ref{InducedCuspForm}) ranges over an orthonormal basis of Hecke-Maass forms associated to $\mathcal P$.  Furthermore,  $c_\mathcal P$ is a fixed non-zero constant.
\end{theorem}

\begin{proof} For proofs see \cite{Arthur1979}, \cite{Langlands1976}, \cite{MW1995}.
\end{proof}

\section{\large\bf Kuznetsov trace formula}\label{KuznetsovTraceFormula}

The Kuznetsov trace formula is derived by computing the inner product of two Poincar{\'e} series in two different ways.  More precisely, let $L=(\ell_1,\ldots,\ell_{m-1}),M=(m_1,\ldots,m_{n-1})\in \Z^{n-1}$ with $\prod\limits_{i=1}^{n-1}m_i\neq 0$ and $\prod\limits_{i=1}^{n-1}\ell_i\neq 0$, and consider the Petersson inner product $\big\langle P^L,P^M\big\rangle$.

 In particular since $P^L, P^M \in \mathcal L^2\left(\SL(n, \mathbb Z) \backslash\mathfrak h^n \right)$  {(see \cite{friedberg1987poincare})},  the inner product can be computed with the spectral expansion of the Poincar\'e series. The geometric approach utilizes the Fourier Whittaker expansion of the Poincar\'e series which involve Kloosterman sums.
\vskip 3pt

The trace formula takes the following form.
\begin{equation} \label{TraceFormula}
\boxed{\boxed{ \underset{\mbox{\scriptsize{spectral side}}}{\underbrace{\mathcal C \; + \; \mathcal E}} \;\; = \;\; \underset{\mbox{\scriptsize{geometric side}}}{\underbrace{\mathcal M \; + \; \mathcal K }}.}}
\end{equation}

Here $\mathcal C$ is the cuspidal contribution and $\mathcal E$ is the Eisenstein contribution.  See Theorem~\ref{ThmSpectralDecomp} for their precise definitions.  The geometric side consists of terms corresponding to elements of the Weyl group.  The identity element gives
the main term $\mathcal M$, and the Kloosterman contribution $\mathcal K$ is the sum of the remaining terms.  See Theorem~\ref{GeomSide} for their precise definitions.  The Kloosterman term $\mathcal{K}$ and the Eisenstein contribution $\mathcal{E}$ will be small with the special choice of the test function $ p_{T,R}$, and they constitute the error term in the main theorem.

\subsection{Spectral side of the Kuznetsov trace formula}

The first way to compute the inner product of the Poincar{\'e} series uses the spectral decomposition of the Poincar\'e series.  

Recall also the definition of the adjoint L-function:  $L(s, \Ad{\phi}) := L(s, \phi\times\overline{\phi})/\zeta(s)$ where $L(s, \phi\times\overline{\phi})$ is the Rankin-Selberg convolution L-function as in \S12.1 of 
\cite{Goldfeld2015}. 

\begin{theorem}[\bf Spectral decomposition for the inner product of  Poincar\'e series]\label{ThmSpectralDecomp} Fix $n\ge 2$ and $L=(\ell_1,\;\ldots\; ,\ell_{n-1}), \, M= (m_1,\;\ldots\;,m_{n-1}) \in\mathbb Z^{n-1}$ .  Then For $\alpha_0:= \big( -\tfrac{n-1}{2}+j-1\big)_{j=1,\ldots,n}$ we have
 \[ \boxed{\ 
 \Big\langle P^L(*, \alpha_0), \; P^M(*, \alpha_0)  \Big\rangle  = \mathcal{C} + \mathcal{E}.\ } \]
\vskip 5pt
With the notation of the Spectral Decomposition Theorem~\ref{SpectralDecomp}, the cuspidal contribution to the Kuznetsov trace formula is
\begin{align*}
  \mathcal{C} := \;\sum_{i=1}^\infty \frac{ \lambda_{\phi_i}(M)\overline{\lambda_{\phi_i}(L)}\cdot\left| p_{T,R}^{n,\#}\left(\, \alpha^{(i)}\,\right)\right|^2}{ L(1,\Ad{\phi_i}) \cdot \prod\limits_{1\le j \ne k\le n} \Gamma\left( \frac{1+\alpha_j^{(i)}-\alpha_k^{(i)}}{2}  \right)},
\end{align*}
and the Eisenstein contribution to the Kuznetsov trace formula is 
\begin{equation*}
\mathcal{E} :=  \sum_{\mathcal P} \sum_{\Phi}\; c_{\mathcal P}\hskip-4pt  \underset{\text{\rm Re}(s_j)=0}{\int\limits_{n_1s_1+\cdots +n_rs_{r}=0}}\hskip -11pt 
A_{E_{\mathcal P,\Phi}}(L, s) \, \overline{A_{E_{\mathcal P,\Phi}}(M, s)}\cdot\Big| p_{T,R}^{n,\#}\big(\alpha_{_{(\mathcal P,\Phi)}}(s)\big)\Big |^2\; ds,
\end{equation*}
for constants   $c_{\mathcal P} >0$.\end{theorem}

\begin{proof}
The proof follows from the Langlands Spectral Decomposition Theorem \ref {SpectralDecomp} with the choices  $F = P^L$ and $G=P^M$.  We have
\begin{align*}
&\langle P^L, P^M\rangle =  \sum_{j=1}^\infty  \frac{\langle P^L, \phi_j\rangle \, \langle\phi_j, P^M\rangle}{\langle\phi_j, \phi_j\rangle} 
\; + 
\sum_{\mathcal P} \sum_{\Phi}\; c_{\mathcal P} \hskip-10pt \underset{\text{\rm Re}(s_j)=0}{\int\limits_{n_1s_1+\cdots +n_rs_{r}=0}} \hskip-10pt\Big\langle F, \, E_{\mathcal P,\Phi}(*\, ,s)\Big\rangle      \Big\langle E_{\mathcal P,\Phi}(*\, ,s), \,G \Big\rangle\; ds.
\nonumber
\nonumber\end{align*}
We then insert the inner products given in Proposition \ref{InnerProducts} below.  Doing so, we see that the cuspidal spectrum is
\begin{align*}
 \sum_{i=1}^{\infty} \frac{\big\langle P^L, \phi_i \big\rangle \big\langle \phi_i, P^M \big\rangle}{\langle \phi_i,\phi_i \rangle} & =
 \sum_{i=1}^\infty \frac{A_{\phi_i}(M) \overline{A_{\phi_i}(L)}}{\mathfrak c_n \cdot \langle \phi_i,\phi_i \rangle} \left| p_{T,R}^{n,\#}(\alpha^{(i)}) \right|^2.
\end{align*}
From Proposition~\ref{PropFirstCoeff}, we see that
 \[ A_\phi(M) \overline{A_\phi(L)} = \lvert A_\phi(1)\rvert^2 \lambda_\phi(M)\overline{\lambda_\phi(L)} = \frac{\mathfrak c_n\cdot\big\langle \phi, \,\phi\big\rangle \cdot \lambda_\phi(M)\overline{\lambda_\phi(L)}}{L(1, \Ad \; \phi) \prod\limits_{1\le j\ne k\le n}\Gamma\big(\frac{1+\alpha_j-\alpha_k}{2}  \big)}. \]
The cuspidal part is now immediate.  The contributions from the Eisenstein series are computed in like manner using Proposition   \ref{MthEisCoeff}.
\end{proof}

\begin{proposition}[\bf The inner product of $P^M$ with an Eisenstein series or Hecke-Maass form]
\label{InnerProducts} Let $M=(m_1,m_2,\ldots,m_{n-1})$. Consider the Eisenstein series 
$E_{\mathcal P, \Phi}(*,s)$, with associated Langlands parameters    $\alpha_{_{\mathcal P,\Phi}}(s)$. Let $\phi$ denote a Hecke-Maass cusp form for $\SL(n,\mathbb Z)$ with Langlands parameter $\alpha$ and $M^{th}$ Fourier coefficient $A_{\phi}(M).$ Then for $\alpha_0 := \big( -\tfrac{n-1}{2}+j-1\big)_{j=1,\ldots,n}$,
\begin{align*}
 \Big\langle \phi , \; P^M(*, \alpha_0)  \Big\rangle  & =\;\frac{1}{\sqrt{\mathfrak c_n}} \;A_{\phi}(M)\cdot p_{T,R}^{n,\#}(\alpha),
\\
 \Big\langle E_{\mathcal P, \Phi}(*,s), \; P^M(*, \alpha_0)  \Big\rangle 
&
=\;\frac{ 1}{\sqrt{\mathfrak c_n}} \;A_{E_{\mathcal P, \Phi}}(M,s)\cdot p_{T,R}^{n,\#}(\alpha_{_{\mathcal P,\Phi}}(s)),
\end{align*}
where  {the inner products on the left are defined by analytic continuation and} $\mathfrak c_n$ is the nonzero constant (depending only on $n$) from Proposition~\ref{PropFirstCoeff}.
\end{proposition}
\begin{proof}
We outline the case of the Hecke-Maass forms. The series definition of the Poincar\'e series converges absolutely for  {sufficiently large} $\re(\alpha_i'-\alpha_{i+1}')$ ($1\leq i\leq n-1$).  It follows that  {for such $\alpha'$} we may unravel the Poincar\'e series $P^M(*,\alpha')$ in the inner product $\langle \phi , P^M \rangle$ with the Rankin-Selberg Method.
The inner product picks out the $M^{th}$ Fourier coefficient of $\phi$ multiplied by a certain Whittaker transform of $p_{T,R}^{(n)}\big(My\big)\cdot I(y,\alpha')$.  {This Whittaker transform has analytic continuation in $\alpha'$ to a region including $\alpha_0$.}   {For sufficiently large $\re(\alpha_i'-\alpha_{i+1}')$, we have} from \eqref{eq:PoincareSeries-def} that
 \begin{align} \label{InnerProductPoincare}
 \big\langle \phi, \; P^M(*, \al')  \big\rangle \; & = \;\frac{A_\phi(M)}
{\sqrt{\mathfrak c_n} \prod\limits_{k=1}^{n-1} m_k^{ {k(n-k)} }}\int\limits_{y_1=0}^\infty  \cdots
 \int\limits_{y_{n-1}=0}^\infty  \overline{p_{T,R}^{(n)} \big(My\big) {\cdot I(y,\alpha')}}\cdot W_\alpha\big(My\big)\; \prod_{k=1}^{n-1} \frac{dy_k}{y_k^{k(n-k)+1}}.
 \end{align}
Note that $I(y,\alpha_0)=1$.  {The integral in \eqref{InnerProductPoincare} converges (as a function of $\alpha'$) to a region which includes $\alpha_0$.  It follows that the analytic continuation in $\alpha'$ to $\alpha_0$ of the inner product satisfies
 \[ \big\langle \phi, \; P^M(*, \al_0)  \big\rangle \; = \; \frac{1}{\sqrt{\mathfrak c_n}} \cdot A_\phi(M)
  \cdot
p_{T,R}^{n,\#}(\alpha). \]
}  The proof for $E_{\mathcal P, \Phi}$ is the same.
\end{proof} 

For $n\ge 2$,  consider a  Hecke-Maass cusp form $\phi$ for $\SL(n,\mathbb Z)$ with Fourier Whittaker expansion given by Proposition \ref{FourierExpansion}. Assume $\phi$ is a Hecke eigenform.   Let $A_\phi(1) := A_\phi(1,1,\ldots,1)$ denote the first Fourier-Whittaker coefficient of $\phi.$ Then we have
$$A_\phi(M) = A_\phi(1) \cdot\lambda_\phi(M)$$
where $\lambda_\phi(M)$ is the Hecke eigenvalue (see Section~9.3 in \cite{Goldfeld2015}), and $\lambda_\phi(1) = 1$.

\begin{proposition}[{\bf First Fourier-Whittaker coefficient of  a Hecke-Maass cusp form}] \label{PropFirstCoeff} Assume $n\ge 2.$ Let $\phi$ be a Hecke-Maass cusp form for $\SL(n,\mathbb Z)$ with Langlands parameters $\alpha=(\alpha_1,\ldots,\alpha_n)$. Then the first coefficient $A_\phi(1)$ is given by 
$$|A_\phi(1)|^2 =  \frac{\mathfrak c_n\cdot\big\langle \phi, \,\phi\big\rangle}{L(1, \Ad{ \phi}) \prod\limits_{1\le j\ne k\le n}\Gamma\big(\frac{1+\alpha_j-\alpha_k}{2}  \big)}$$
where $\mathfrak c_n \ne 0$ is a constant depending on $n$ only.
\end{proposition}

\begin{proof}
See \cite{GMW2021}.
\end{proof}

\begin{proposition}[\bf  \boldmath The $M^{\rm{th}}$ Fourier coefficient of   $E_{\mathcal P, \Phi}$] \label{MthEisCoeff}
  
  Let $s = (s_1,s_2, \ldots,s_r)\in\mathbb C^r,$  where $\sum\limits_{i=1}^r n_is_i =0.$
Consider  $E_{\mathcal P, \Phi}(*,s)$ with associated Langlands parameters $\alpha_{_{\mathcal P,\Phi}}(s)$ as defined in (\ref{langlandsparamsforPhi}).   Assume that each Hecke-Maass form $\phi_k$ (with $1\le k\le r$) occurring in $\Phi$ has Langlands parameters
  $\alpha^{(k)} := (\alpha_{k,1},\ldots,\alpha_{k,n_k})$ with the convention that if $n_k=1$ then $\alpha_{k,1}=0.$ We also assume that each $\phi_k$ is normalized to have Petersson norm $\langle \phi_k, \phi_k\rangle = 1.$

Let $L^*(1+s_j-s_\ell, \phi_j\times\phi_\ell)$ denote the completed Rankin-Selberg L-function if $n_j\ne1\ne n_\ell$; otherwise define$$L^*(1+s_j-s_\ell, \;\phi_j\times\phi_\ell) = \begin{cases} L^*(1+s_j-s_\ell,\, \phi_j) & \text{if}\; n_\ell =1 \;\text{and}\;n_j\ne 1,\\
L^*(1+s_j-s_\ell,\, \phi_\ell) & \text{if}\; n_j=1 \;\text{and}\;n_\ell\ne 1,\\
\zeta^*(1+s_j-s_\ell) & \text{if}\; n_j=n_\ell=1,
\end{cases}$$ where $\zeta^*(w) = \pi^{-\frac{w}{2}}\Gamma\left( \frac{w}{2} \right)\zeta(w)$ is the completed Riemann $\zeta$-function.  Also define$$L^*(1,\,\text{\rm Ad}\;\phi_k) = L(1,\,\text{\rm Ad}\;\phi_k) \prod_{1\le i \ne j\le n_k}  \Gamma\left(\frac{1+\alpha_{k,i}-\alpha_{k,j}}{2}   \right),$$with the convention that $L^*(1,\,\text{\rm Ad}\;1)=1$.

  Let $M = (m_1,m_2,\ldots, m_{n-1}) \in \mathbb Z_{>0}^{n-1}$.  Per our convention (Definition~\ref{VectorMatrix}), we may think of $M$ as a vector or a diagonal matrix.  
 Then the $M^{th}$ term in the Fourier-Whittaker expansion of $E_{\mathcal P, \Phi}$ is
\begin{align*}
\int\limits_{U_n(\mathbb Z)\backslash U_n(\mathbb R)} E_{\mathcal P, \Phi}(ug, s)\, \overline{\psi_M(u) }\; du \; = \; \frac{A_{E_{\mathcal P, \Phi}}(M,s)}
{\prod\limits_{k=1}^{n-1}  m_k^{k(n-k)/2}} \; W_{\alpha_{_{\mathcal P,\Phi}}(s)}\big(M g\big),
\end{align*}
where $A_{E_{\mathcal P, \Phi}}(M,s) = A_{E_{\mathcal P, \Phi}}\big((1,\ldots,1),s\big) \cdot \lambda_{E_{\mathcal P, \Phi}}(M,s),$
 \begin{align}\label{HeckeCoeff}
    \lambda_{E_{P,\Phi}}\big((m,1,\ldots,1),s\big) & =   
    \underset{c_1c_2\cdots c_n=m} {\sum_{c_1,\ldots,c_n\in\mathbb Z_{>0}}}    \hskip-10pt \lambda_{\phi_1}(c_1)  \cdots \lambda_{\phi_r}(c_r)
    \cdot c_1^{s_1}  \cdots c_r^{s_r}
\end{align}
 is the $(m,1,\ldots,1)^{th}$ (or more informally  the $m^{th}$) Hecke eigenvalue of $E_{\mathcal P,\Phi}$, and  \begin{align*} 
 & A_{E_{\mathcal P, \Phi}}\big((1,\ldots,1),s\big) =  d_0\underset{n_k\ne1}{\prod_{k=1}^r} L^*\big(1, \text{\rm Ad}\; \phi_k\big)^{-\frac12 }
  \prod_{1\le j<\ell\le r}   L^*\big(1+s_j-s_{\ell}, \;\phi_j\times\phi_{\ell}\big)^{-1}
 \end{align*} for some constant 
 $d_0\ne0$ depending  only on $n$. 
\end{proposition}
\begin{proof} See  \cite{GSW24}.\end{proof}

\subsection{Geometric side of the Kuznetsov trace formula}

In this section, we obtain explicit descriptions of the terms $\mathcal{M}$ and $\mathcal{K}$ appearing on the geometric side of the Kuznetsov trace formula.  In order to do this, we introduce Kloosterman sums for $\SL(n,\Z)$, which appear in the Fourier expansion of the Poincar{\'e} series.  In the inner product $\big\langle P^L,P^M \big\rangle$, we replace $P^L$ with its Fourier expansion and unravel $P^M$ following the Rankin-Selberg method.

\begin{theorem} [\bf Geometric side of the trace formula]\label{GeomSide}
Fix $L= (\ell_1,\ldots,\ell_{n-1})$ and\linebreak $M = (m_1,\ldots,m_{n-1})\in \Z^{n-1}$  ($\mathfrak c_n$ is a nonzero constant; see Proposition~\ref{PropFirstCoeff}). It follows that for $\alpha_0 := \big( -\tfrac{n-1}{2}+j-1\big)_{j=1,\ldots,n}$,
 \[ \boxed{\ \big\langle P^L(*, \alpha_0), \; P^M(*, \alpha_0)  \big\rangle\;  = \; \mathcal{M}\; + \; \mathcal{K}.\ } \]
For $w_1$ the trivial element in the Weyl group $W_n$, we define
 \[ \mathcal{M} :=  \mathcal{I}_{w_1}, \qquad \mbox{and} \qquad 
 \mathcal{K} :=  \sum_{\substack{w\in W_n \\ w\neq w_1}} \mathcal{I}_w, \]
where
\begin{multline}\label{eq:IwDef}
\mathcal I_w :=  \;  \sum_{v\in V_n}\sum_{c_1=1}^\infty 
\cdots
\sum_{c_{n-1}=1}^\infty   \frac{S_w(\psi_L,\psi_M^v,c)}{\mathfrak c_n \cdot \prod\limits_{k=1}^{n-1} (m_k \ell_k)^{\frac{k(n-k)}{2}}}\int\limits_{y_1=0}^\infty 
\cdots
\int\limits_{y_{n-1}=0}^\infty\;\;\int\limits_{ {U}_w(\mathbb Z)\backslash  {U}_w(\mathbb R)}\;\int\limits_{\overline{{U}}_w(\mathbb R)}\\ 
\cdot
\psi_L(w u y)\,\overline{\psi_M^v(u)}\;
  p_{T,R}^{(n)}(Lc w u y)\,\overline{ p_{T,R}^{(n)}(My)}\;  d^* u
\,\frac{dy_1\cdots dy_{n-1}}{\prod\limits_{k=1}^{n-1}y_k^{k(n-k)+1}}.
\end{multline}
\end{theorem}

\begin{proof}
We  compute the  inner product
\begin{align*}&
\lim_{\alpha\to\alpha_0}\big\langle P^L\left(*, \,\alpha \right), \, P^M\left(*, \,\alpha \right) \big\rangle  = \lim_{\alpha\to\alpha_0} \int\limits_{\SL(n,\mathbb Z)\backslash\mathfrak h^n} P^L\left(g, \alpha \right) \cdot  \overline{P^M\left(g, \alpha \right)} \; dg
\\
  = & \frac{1}{\sqrt{\mathfrak c_n} \prod\limits_{k=1}^{n-1} m_k^{\frac{k(n-k)}{2}}} \lim_{\alpha\to\alpha_0} \int\limits_{U_n(\mathbb Z)\backslash\mathfrak h^n} P^L\left(g, \alpha \right) \cdot  \overline{\psi_M(g)\,  p_{T,R}^{(n)}(Mg)\, I(g,\alpha)} \; dg
  \\
 =& \frac{1}{\sqrt{\mathfrak c_n}} \left( \prod_{k=1}^{n-1} m_k^{-\frac{k(n-k)}{2}} \right) \lim_{\alpha\to\alpha_0} \hskip -1pt \int\limits_{\substack{y\in\R^{n-1}\\ y>0}} \hskip -1pt \left( \;\int\limits_{U_n(\mathbb Z)\backslash U_n(\mathbb R)} \hskip -3pt P^L\left(uy, \alpha \right) \cdot \overline{\psi_M(u)} \;du\right) \overline{ p_{T,R}^{(n)}(My)\, I(y,\alpha)}\;\, dy.
\end{align*}

Note that, as $\alpha\to\alpha_0$, the function $I(g,\alpha)\to 1$ (for any $g\in \mathfrak{h}^n$) and $\prod\limits_{k=1}^{n-1}c_k^{\alpha_k-\alpha_{k+1}+1}\to 1$.  It follows from this and Proposition~\ref{FourierExp} above that
\begin{align*}
\mathfrak c_n \cdot & \prod\limits_{k=1}^{n-1} (m_k \ell_k)^{k(n-k)/2} \cdot \lim_{\alpha\to\alpha_0}\big\langle P^L\left(*, \,\alpha \right), \, P^M\left(*, \,\alpha \right) \big\rangle
\\
& \quad = \lim_{\alpha\to\alpha_0} \sum_{w\in W_n}\sum_{v\in V_n}\sum_{c_1=1}^\infty \cdots \sum_{c_{n-1}=1}^\infty   \frac{S_w(\psi_L,\psi_M^v,c)}{\mathfrak c_n \cdot \prod\limits_{k=1}^{n-1} (m_i \ell_i)^{\frac{i(n-i)}{2}} \prod\limits_{k=1}^{n-1}c_k^{\alpha_k-\alpha_{k+1}+1}   }  \int\limits_{y_1=0}^\infty \cdots \int\limits_{y_{n-1}=0}^\infty\;\;\int\limits_{{U}_w(\mathbb Z)\backslash {{U}}_w(\mathbb R)}\;\int\limits_{\overline{{U}}_w(\mathbb R)}
\\
& \quad \qquad
\cdot
\psi_L(w u y)\,\overline{\psi_M^v(u)}\;
 p_{T,R}^{(n)}(Lc w u y)\,\overline{ p_{T,R}^{(n)}(My)}\; I(w u  y,\alpha) \;\overline{ I(y,\alpha)}\;
  d^* u
\,\frac{dy_1\cdots dy_{n-1}}{\prod\limits_{k=1}^{n-1}y_k^{k(n-k)+1}}\\
& \quad =  \; \sum_{w\in W_n}\sum_{v\in V_n}\sum_{c_1=1}^\infty \cdots \sum_{c_{n-1}=1}^\infty   \frac{S_w(\psi_L,\psi_M^v,c)}{\mathfrak c_n \cdot \prod\limits_{k=1}^{n-1} (m_k \ell_k)^{\frac{k(n-k)}{2}}}\int\limits_{y_1=0}^\infty \cdots \int\limits_{y_{n-1}=0}^\infty\;\;\int\limits_{ {U}_w(\mathbb Z)\backslash {{U}}_w(\mathbb R)}\;\int\limits_{\overline{{U}}_w(\mathbb R)}
\\
& \quad \qquad
\cdot
\psi_L(w u y)\,\overline{\psi_M^v(u)}\;
  p_{T,R}^{(n)}(Lc w u y)\,\overline{ p_{T,R}^{(n)}(My)}\;  du
\,\frac{dy_1\cdots dy_{n-1}}{\prod\limits_{k=1}^{n-1}y_k^{k(n-k)+1}}\\
& \quad = \sum_{w\in W_n} \mathcal I _w,
\end{align*}
as claimed.
\end{proof}

\section{\large\bf Asymptotic formula for the main term}

\begin{proposition}[\bf Main term in the trace formula]\label{prop:mainterm} 
Let $L=(\ell_1,\ldots,\ell_{n-1}),\, M=(m_1,\ldots,m_{n-1})  \in\mathbb Z^{n-1}$ satisfy $\prod\limits_{i=1}^{n-1} \ell_i\neq 0$ and $\prod\limits_{i=1}^{n-1} m_i\neq 0$.  There exist  fixed constants $\mathfrak c_1,\ldots, \mathfrak c_{n-1}>0$ (depending only on $R$ and $n$) such that the main term $\mathcal M$ in the Kuznetsov trace 
formula \eqref{TraceFormula} is given by
 $$\mathcal M =  \delta_{L,M}\cdot 
 \left(\left(\sum_{i=1}^{n-1} \mathfrak c_i \cdot T^{R(2\cdot D(n) + n(n-1)) + n-i}\right) + \mathcal O\left( T^{R(2\cdot D(n) + n(n-1))}  \right) \right),$$
where $D(n)=\frac{1}{2} \binom{2n}{n}-\frac{n(n-1)}{2}-2^{n-1}$ and $\delta_{L,M}$ is the Kronecker symbol (i.e., $\delta_{L,M}=0$ if $L\neq M$ and $\delta_{L,L}=1$).
  \end{proposition}
  
\begin{proof}
It follows from the definition  $\mathcal{M}=\mathcal{I}_{w_1}$, making the change of variables $y\mapsto M^{-1}y$, and noting that $U_{w_1}(\Z)=U_n(\Z)$ and $U_{w_1}(\R)=U_n(\R)$, that
\begin{align*}
\mathcal M & = \frac{1}{\mathfrak c_n}\cdot \int\limits_{y_1=0}^\infty \cdots \int\limits_{y_{n-1}=0}^\infty\;\left(\;\int\limits_{{U}_n(\mathbb Z)\backslash {{U}}_n(\mathbb R)} \hskip-5pt\psi_{L}(u)\overline{\psi_M(u)}\; d^*u\right) p_{T,R}(LM^{-1}y)\, \overline{p_{T,R}(y)} \;\frac{dy_1\cdots dy_{n-1}}{\prod\limits_{i=1}^{n-1} y_i^{i(n-i)+1}} \\& = \delta_{L,M}\cdot \mathfrak d_n\int\limits_{y_1=0}^\infty \cdots \int\limits_{y_{n-1}=0}^\infty | p_{T,R}^{(n)}(y)|^2 \,\frac{dy_1\cdots dy_{n-1}}{ \prod\limits_{i=1}^{n-1} y_i^{i(n-i)+1} } = \delta_{L,M}\cdot \mathfrak d_n\big\langle  p_{T,R},  p_{T,R}\big\rangle \\
 & = \; \delta_{L,M}\cdot \mathfrak d_n\int\limits_{\substack{\halpha_n=0 \\ \re(\alpha)=0}}  \frac{\big| p_{T,R}^{n,\#}(\alpha)\big|^2 }{\prod\limits_{1\le j\ne k\le n}\Gamma\left(\frac{\alpha_j-\alpha_k}{2}  \right) } \;d\alpha                     \\
 & = \delta_{L,M}\, \mathfrak d_n\cdot \big\langle  p_{T,R}^{n,\#}, \; p_{T,R}^{n,\#}\big\rangle,
 \end{align*}
where the representation in terms of the norm of $p_{T,R}^{n,\#}$ follows from the Plancherel formula in Corollary 1.9 of \cite{MR2982424} and $\mathfrak d_n$ is a nonzero constant depending only on $n$. 
Hence the main term for $\GL(n)$ is thus
 \[ \mathcal M = \delta_{L,M} \mathfrak d_n \cdot \int\limits_{\substack{\halpha_n=0 \\ \re(\alpha_j)=0}} \frac{ \Big|e^{\frac{\alpha_1^2+\cdots+\alpha_n^2}{2T^2}}
   \mathcal{F}_R^{(n)}\big(\frac{\alpha}{2}\big)  \prod\limits_{1\leq\, j \ne k\, \leq n}
      \Gamma\left(\textstyle{\frac{2R+1+\alpha_j - \alpha_k}{4}}\right)\Big|^2}{   \prod\limits_{1\le j\ne k\le n}\Gamma\left(\frac{\alpha_j-\alpha_k}{2}  \right)     } \; d\alpha. \]
Let $\alpha_j = i\tau_j$ with $\tau_j\in\mathbb R$.  It then follows from Stirling's asymptotic formula 
that
\begin{align*}
\mathcal M & \sim \delta_{L,M} \mathfrak d_n \cdot
\int\limits_{\widehat{\tau}_n=0} 
    e^{\frac{-\tau_1^2-\cdots-\tau_n^2}{T^2}}
    \Big(\mathcal{F}_R^{(n)}\big(\tfrac{i\tau}{2}\big) \Big)^2 \prod\limits_{1\leq\, j < k\, \leq n}
      \big( 1+\lvert \tau_j-\tau_k \rvert \big)^{2R}\ d\tau.
\end{align*}
If we now make the change of variables $\tau_j\to \tau_j T$ for each $j=1,\ldots,n$, and we use the fact that the degree of $\mathcal{F}_1^{(n)}$ is $D(n)$ (see Definition~\ref{def:FRn}) it follows that, if $L=M$, as $T\to\infty$ we have $\mathcal M  \sim \mathfrak{c} T^{R\cdot\big(2D(n)+n(n-1) \big)+n-1}$, where 
\begin{align*} 
 \mathfrak c & 
 = \mathfrak d_n \cdot \int\limits_{\widehat{\tau}_n=0} 
    e^{-\tau_1^2-\cdots-\tau_n^2}
    \Big(\mathcal{F}_R^{(n)}\big(\tfrac{i\tau}{2}\big) \Big)^2 \prod\limits_{1\leq\, j < k\, \leq n}
      \big( 1+\lvert \tau_j-\tau_k \rvert \big)^{2R}\ d\tau,
\end{align*}
and otherwise, the main term is zero.  This gives the $i=1$ term in the statement of the proposition.  The method of proof can be extended by using additional terms in Stirling's asymptotic expansion for the Gamma function to obtain the additional terms.
\end{proof}

\begin{remark}
Note that this doesn't agree with \cite{GSW21} in the case of $n=4$ because we have used a different normalization.  Namely, the linear factors of $\mathcal{F}_R^{(n)}$ agree with those defined previously, but we take a different power of each.  Also, the  {gamma} factors which appear in $p_{T,R}^{n,\#}$ have a different $R$: namely, what was $2+R$ in each  {gamma} factor previously has been replaced by $2R+1$ here.
\end{remark}

\section{\large\bf Bounding the geometric side}\label{sec:GeometricSideBounds}
\vskip 5pt

The goal of this section is to use the bound given in Theorem~\ref{th:pTRn-bound} to prove the following, i.e., to bound $\mathcal{K}$ the geometric side of the Kuznetsov trace formula.

\begin{proposition}\label{prop:Iwbounds}
Let $\mathcal{I}_w$ be as above.  Let $M=(m_1,\ldots,m_{n-1}), \;  L=(\ell_1,\ldots,\ell_{n-1}) \in (\mathbb Z_{>0})^{n-1}$.  Let $\rho \in \frac12+\Z$.   {Let $D(n) =  \frac12 \binom{2n}{n} - \frac{n(n-1)}{2} - 2^{n-1}$ as in \eqref{eq:Dn-explicit}.} Then for $R$ sufficiently large and any $\varepsilon>0$, we have 
 \[ \big| \mathcal{I}_w \big| \ll_{\varepsilon, R} \; T^{\varepsilon+R\big(2D(n)+n(n-1)\big) + \frac{(n-1)(n+4)}{2} - \left\lfloor \frac{n-1}{2} \right\rfloor - \rho n - \Phi(w)} \cdot \prod_{i=1}^{n-1}(\ell_im_i)^{2\rho+ \frac{n^2+1}{4}}, \]
where if $w=w_{(n_1,\ldots,n_r)}$ with $r\geq 2$,
 \[ \Phi(w) := \Phi(n_1,\ldots,n_r) := \frac12 \sum_{k=1}^{r-1} (n_k+n_{k+1})(n-\hn_k)\hn_k. \]
\end{proposition}

\begin{remark}\label{rmk:lm-bound}
Assuming the lower bound conjecture for Rankin-Selberg L-functions, the resulting bound for the Eisenstein series contribution to the Kuznetsov trace formula (see Theorem~\ref{th:EisensteinBound}) is of the magnitude $T$ to the power $R\big( 2D(n) + n(n-1) \big) + \eps$.  Therefore, given Proposition~\ref{prop:Iwbounds} and Lemma~\ref{lem:Phi-properties} (which says that $\Phi(w)\geq \Phi(1,n-1)=\frac{n(n-1)}{2}$), in order for the bound from the geometric side of the trace formula to be less than the Eisenstein series contribution, it suffices that 
 \[ \frac{(n-1)(n+4)}{2} - \left\lfloor \frac{n-1}{2} \right\rfloor - \rho n - \frac{n(n-1)}{2} \leq 0, \]
which simplifies to give
 \[ \rho \geq 
 \begin{cases}
   \frac32 - \frac{3}{2n} & \mbox{ if $n$ is odd}, 
   \\
   \frac32 -\frac{1}{n} & \mbox{ if $n$ is even}.
 \end{cases}  \]
Since we require that $\rho \in \frac12+\Z$, we find that it suffices to take $\rho=\frac32$ universally, meaning that the exponent of each term $\ell_i m_i$ can be taken to be $\frac{n^2+13}{4}$.  In particular, for the case of $n=4$, we see that this exponent is $\frac{29}{4}$ which is an improvement on the bound of $\frac{15}{2}$ obtained in \cite{GSW21}.
\end{remark}

As remarked above, the main result that we will need is Theorem~\ref{th:pTRn-bound}  {or, more specifically, Remark~\ref{rmk:pTRnBound},} which states that for any $0<\eps<\frac12$, and for $a=(a_1,a_2,\ldots,a_{n-1})$ satisfying $\lfloor a_j \rfloor+\eps < a_j < \lceil a_j \rceil -\eps$ for each $j=1,\ldots,n-1$, that
\begin{equation}\label{eq:final-pTRbound}
\boxed{
 \big| p_{T,R}^{(n)}(y) \big| \ll \delta^{-\frac12}(y) \cdot \big\lVert y \big\rVert^{2a} \cdot T^{\eps+\frac{(n+4)(n-1)}{4}+ R\cdot\big(D(n)+\frac{n(n-1)}{2} \big) - \sum\limits_{j=1}^{n-1}B( a_j)} }.
\end{equation}
(The terms $\delta^{-\frac12}(y)$ and $\lVert y \rVert^{2a}$  {are defined in Section~\ref{sec:Iw-notation}} below.  The function $B$ is defined in Theorem~\ref{th:ITRnaBound}.)

This bound for $p_{T,R}^{(n)}(y)$ is obtained via an integral representation denoted $p_{T,R}^{(n)}(y;b)$ (see \eqref{eq:pTRnstart}) over variables $s=(s_1,\ldots,s_n)$ valid for any $b=(b_1,\ldots,b_n)$ with $b_j>0$ for each $j=1,\ldots,n-1$.  The integral is taken over the lines $\re(s_j)=b_j$.  Essentially, the bound is then obtained by moving the lines of integration to $\re(s_j) = -a_j$ for some $a=(a_1, \ldots,a_{n_1})\in \big(\R_{>0}\big)^{n-1}$.

The strategy for proving Proposition~\ref{prop:Iwbounds} will be to, first, introduce notation to rewrite $\mathcal I_w$ in a simplified form.  We do this in Section~\ref{sec:Iw-notation}.  Then, in Section~\ref{sec:Iwbounds-proof} we give bounds for $\mathcal I_w$ obtained by applying \eqref{eq:final-pTRbound} to $\lvert p_{T,R}(Lcwuy)\rvert$ (with a parameter $a=(a_1,\ldots,a_{n-1})$) and to $\lvert p_{T,R}(My)\rvert$ (with a parameter $b=(b_1,\ldots,b_{n-1})$) for general $a,b\in (\R_{>0})^{n-1}$.  In particular, we establish \eqref{eq:newerIw}, bounding $\lvert \mathcal I_w\rvert$ in terms of the product of three independent quantities $K(c,w;a)$, $X(u,w;a)$ and $Y(y,w;a,b)$.  In Section~\ref{sec:abconditions}, we will show that $K(c,w;a)$ will converge provided that $a$ satisfies certain conditions (independent of $w$), and that for this choice of $a$, $X(u,w;a)$ also converges.  We then determine $b$ (dependent on by $w$ and $a$) for which $Y(y,w;a,b)$ is also convergent.  Finally, in Section~\ref{sec:Iwbounds-proof-final}, we complete the proof of Proposition~\ref{prop:Iwbounds} by simplifying the expression for the given choices of $a$ and $b$.

\subsection{Rewriting $\mathcal{I}_w$}\label{sec:Iw-notation}

Let $T_n(\R)$ and $U_n(\R)$ be the subgroups of $\GL_n(\R)$ consisting of diagonal matrices (with positive terms) and upper triangular unipotent matrices, respectively. Recall that if $t=\diag(t_1,\ldots,t_n)\in T_n(\R)$ and $u\in U_n(\R)$, the modular character $\delta:T_n(\R)\to \R$ is defined to satisfy $d(t^{-1}ut)=\delta(t)\, du$.  Explicitly, it is given by
 \[ { \delta(t) = \prod_{i=1}^n t_i^{2i-n-1}}. \]
More generally, if $a=(a_1,\ldots,a_{n-1})\in \R^{n-1}$, for
 \[ y=(y_1,\ldots,y_{n-1}) := \diag(y_1\cdots y_{n-2}y_{n-1},\ldots,y_1y_2,y_1,1) \]
with $y_1,\ldots,y_{n-1}>0$ we define
 \[ \lVert y\rVert^a := \prod_{k=1}^{n-1} y_k^{a_k}. \]
One checks that in the special case of $a_j=\frac{j(n-j)}{2}$ for $j=1,\ldots,n-1$,
\begin{equation}\label{eq:delta-Vert}
 \delta^{-\frac12}(y) = \lVert y \rVert^a.
\end{equation}

Similarly, if ${}^t U_n(\R)$ is the subgroup of $\GL_n(\R)$ consisting of lower triangular unipotent matrices and 
 \[ \overline{U}_w := \big(w^{-1}\; {}^t U_n(\R)\; w\big) \cap U_n(\R), \]
then we can consider the character $\delta_w$ on $T_n(\R)$ which satisfies $d(tut^{-1}) = \delta_w(t) du$ upon restricting the measure on $U_n(\R)$ to $\overline{U}_w$.  It can be checked that
\begin{equation}\label{eq:deltaw}
 \delta_w(y) = \delta^{\frac12}(y) \cdot \delta^{-\frac12}(wyw^{-1}).
\end{equation}

Recall from Theorem~\ref{GeomSide} that for $L=(\ell_1,\ldots,\ell_{n-1}),M=(m_1,\ldots,m_{n-1})\in (\Z_{>0})^{n-1}$ and
\[ c = \left(\begin{smallmatrix} 1/c_{n-1} & & & & \\
& c_{n-1}/c_{n-2} & & & \\
& & \ddots & & \\
& & & c_2/c_1 & \\
& & & & c_1  \end{smallmatrix}\right), \] 
where $c_i\in \mathbb \Z_{>0}$ for $i=1,\ldots,n-1$,
the Kloosterman contribution to the Kuznetsov trace formula is given by
 \[\mathcal K \; = \; \sum\limits_{\substack{w\in W_n\\ w\neq w_1}} \mathcal I_w, \]
where, using the notation defined above and letting $dy^\times$ denote the measure $\prod\limits_{k=1}^{n-1} \frac{dy_k}{y_k}$,
\begin{multline}\label{eq:newIw}
\mathcal I_w :=  \; \mathfrak c_n^{-1} \sum_{v\in V_n}\sum_{c_1=1}^\infty 
\cdots
\sum_{c_{n-1}=1}^\infty   S_w(\psi_L,\psi_M^v,c)\; \int\limits_{\substack{y=(y_1,\ldots,y_{n-1}) \\ y_1,\ldots,y_{n-1}>0}} \;\int\limits_{ {U}_w(\mathbb Z)\backslash  {U}_w(\mathbb R)}\;\int\limits_{\overline{{U}}_w(\mathbb R)}\\ 
\cdot
\delta^{\frac12}(LM) \cdot \delta(y) \cdot
\psi_L(w u y)\,\overline{\psi_M^v(u)}\;
  p_{T,R}^{(n)}(Lc w u y)\,\overline{ p_{T,R}^{(n)}(My)}\;  d^* u
\; dy^\times.
\end{multline}
We recall that by Friedberg \cite{friedberg1987poincare},  $\mathcal{I}_w$ is identically zero unless $w$ is \emph{relevant} (see Definition~\ref{def:relevant}).

\subsection{Bounds for $\mathcal I_w$ in terms of $a$ and $b$}\label{sec:Iwbounds-proof}

Since $p_{T,R}^{(n)}(g)$ is determined by the Iwasawa decomposition of $g$, we first make the change of variables $u\mapsto y^{-1}uy$.  Then \eqref{eq:newIw} implies that
\begin{multline}\label{eq:prenewerIw}
 \lvert \mathcal I_w \rvert \ll  \;  \sum_{v\in V_n}\sum_{c_1=1}^\infty 
\cdots
\sum_{c_{n-1}=1}^\infty   \lvert S_w(\psi_L,\psi_M^v,c) \rvert 
\int\limits_{\substack{ y=(y_1,\ldots,y_{n-1}) \\ y_1,\ldots,y_{n-1}>0 }}\;\int\limits_{ {U}_w(\mathbb Z)\backslash  {U}_w(\mathbb R)}\;\int\limits_{\overline{{U}}_w(\mathbb R)}\\ 
\cdot
  \, \delta^{\frac12}(M) \cdot \delta^{\frac12}(L) \cdot \delta_w(y) \cdot \delta(y)\cdot \lvert p_{T,R}^{(n)}(Lcwyu)\rvert \,\lvert  p_{T,R}^{(n)}(My)\rvert\;  d^* u
\; dy^\times.
\end{multline}
 {
For the purposes of our analysis, we break up the integral in the $y$-variables.  To this end, let 
 \[ I_0 := (0,1], \qquad I_1 = (1,\infty). \]
For $\tau = (\tau_1,\ldots,\tau_{n-1}) \in \{0,1\}^{n-1}$, define
 \[ I_\tau := I_{\tau_1} \times \cdots \times I_{\tau_{n-1}} .\]
Hence,
 \[ \int\limits_{\substack{y=(y_1,\ldots,y_{n-1}) \\ y_1,\ldots,y_{n-1}>0}} = \sum_{\tau \in \{0,1\}^{n-1}} \;
\int\limits_{I_\tau}, \]
and \eqref{eq:prenewerIw} becomes 
\[ \lvert \mathcal{I}_w\rvert \ll \sum\limits_\tau \lvert \mathcal{I}_w(\tau)\rvert \]
where
\begin{multline}\label{eq:newerIw} 
 \lvert \mathcal I_w(\tau) \rvert \;:=  \;  \sum_{v\in V_n}\sum_{c_1=1}^\infty 
\cdots
\sum_{c_{n-1}=1}^\infty   \lvert S_w(\psi_L,\psi_M^v,c) \rvert 
\int\limits_{y_1\in I_{\tau_1}} \int\limits_{y_2\in I_{\tau_2}} \cdots \int\limits_{y_{n-1}\in I_{\tau_{n-1}}}\;\\ \cdot \int\limits_{ {U}_w(\mathbb Z)\backslash  {U}_w(\mathbb R)}\; \int\limits_{\overline{{U}}_w(\mathbb R)}
  \, \delta^{\frac12}(M) \cdot \delta^{\frac12}(L) \cdot \delta_w(y) \cdot \delta(y)\cdot \lvert p_{T,R}^{(n)}(Lcwyu)\rvert \,\lvert  p_{T,R}^{(n)}(My)\rvert\;  d^* u
\; dy^\times.
\end{multline}
Our strategy is now to, for each choice of $\tau$, replace the terms with $p_{T,R}^{(n)}$ with the bound from \eqref{eq:final-pTRbound} (in the first instance using a choice of $a=(a_1,\ldots,a_{n-1})\in \R^{n-1}$, and in the second instance using $b=(b_1,\ldots,b_{n-1})\in \R^{n-1}$).  Then we need to find choices of $a$ and $b$ for which the corresponding
integrals converge and give good bounds.}

Recall that if $g=utk$ is the Iwasawa decomposition of an element $g\in \GL_n(\R)$, then $p_{T,R}^{(n)}(g)=p_{T,R}^{(n)}(t)$.  With this in mind, consider the Iwasawa decomposition $wu=u_0t k$, where $u_0\in U_n(\R)$, $t\in T_n(\R)$ and $k\in O(n,\R)$.  Then
 \[ Lcwyu = Lc(wyw^{-1}) u_0 t k = u_1 Lc(wyw^{-1}) t k \qquad \mbox{($u_1 = \big(Lcwyw^{-1}\big)^{-1} u_0 \big(Lcwyw^{-1}\big)$)} \]
is the Iwasawa form of $Lcwyu$, hence $ {\big\lvert p_{T,R}^{(n)}(Lcwyu)\big\rvert = \big\lvert p_{T,R}^{(n)}(Lc wyw^{-1}t)\big\rvert}$.  Recall that the Iwasawa form of $wu$ is assumed to be $u_0 t k$, meaning $wu=u_0 tk$ where $u_0\in U_n(\R)$, $t\in T_n(\R)$ and $k\in O(n,\R)$.  It can be shown \cite{Jacquet1967}  that 
 \begin{equation}\label{eq:t_inTermsOf_u}
  t^2 = \left(\begin{smallmatrix} 1/\xi_{n-1} & & & & \\
& \xi_{n-1}/\xi_{n-2} & & & \\
& & \ddots & & \\
& & & \xi_2/\xi_1 & \\
& & & & \xi_1  \end{smallmatrix}\right),
\end{equation}
where $\xi_i=\xi_i(wu) \geq 1$ for any $u\in \overline{U}_w$.  For example, in the case $n=4$ and $w=\wlong=w_{(1,1,1,1)}$, we find that $\overline{U}_{\wlong}=U_4(\R)$ and, for
\[ u   = \left(\begin{matrix}1 & x_{12} & x_{13} & x_{14}\\
0 & 1 & x_{23} & x_{24}\\
0&0&1 & x_{34}\\
0 & 0 & 0 & 1
\end{matrix}\right), \]
that
\begin{align*}
 \xi_1(\wlong u) & = 1+x_{12}^2+x_{13}^2 + x_{14}^2,
 \\
 \xi_2(\wlong u) & = 1+x_{23}^2+x_{24}^2 + \big( x_{12} x_{24} - x_{14} \big)^2 
+\big( x_{12} x_{23} - x_{13} \big)^2
+\big( x_{13} x_{24} - x_{14} x_{23}   \big)^2,
 \\
 \xi_3(\wlong u) & = 1+x_{34}^2 + \big(x_{23} x_{34} - x_{24}\big )^2  + \big( x_{12} x_{23} x_{34} - x_{13} x_{34} - x_{12} x_{24} + x_{14}\big )^2.
\end{align*}
In general, the values $\xi_i$ are always of the form $1$ plus a sum of squares of functions consisting of the entries of $u$.

From \eqref{eq:final-pTRbound}, replacing $a$ with $b$, we see that $\big| p_{T,R}^{(n)}(My) \big|$ is bounded by
 \[ \ll \delta(M)^{-\frac12} \cdot \big \lVert M \big\rVert^{2b} \cdot \delta(y)^{-\frac12} \cdot \big \lVert y \big\rVert^{2b} \cdot T^{\eps+\frac{(n+4)(n-1)}{4}+ R\cdot\big(D(n)+\frac{n(n-1)}{2} \big) - \sum\limits_{j=1}^{n-1}B( b_j)}. \]
To similarly bound $\lvert p_{T,R}^{(n)}(Lc wyw^{-1}t)\rvert$, we first remark that since
 \[ c = c_1 \left( \begin{smallmatrix} 
 d_1d_2\cdots
    d_{n-1} & & & \\
    & \hskip -30pt d_1d_2\cdots d_{n-2} & & \\
    & \ddots &  & \\
    & & \hskip -5pt d_1 &\\
    & & &  1\end{smallmatrix} \right) =: c_1 d \qquad \mbox{where $d_i = \frac{c_{i-1}c_{i+1}}{c_i^2}$}, \]
setting $c_0=c_n:=1$ (and $a_0=a_n:=0$ as usual), we see that
 \[ \delta^{-\frac12}(c) \cdot \lVert c \rVert^{2a} = \delta(c)^{-\frac12} \prod_{i=1}^{n-1} \left( \tfrac{c_{i-1}c_{i+1}}{c_i^2} \right)^{2a_i} = \prod_{k=1}^{n-1} c_k^{-1+2a_{i-1}-4a_i+2a_{i+1}}. 
 \]
Therefore, it follows that 
\begin{multline*}
 \big\lvert p_{T,R}^{(n)}(Lc wy w^{-1} t) \big\rvert \;
 \ll \; \frac{ \delta(L)^{-\frac12}\cdot \lVert L \rVert ^{2a} \cdot \delta(t)^{-\frac12}\cdot \lVert t \rVert^{2a}  }{ \prod\limits_{k=1}^{n-1} c_k^{1-2a_{i-1}+4a_i-2a_{i+1}} } 
 \cdot
 \delta(wyw^{-1})^{-\frac12} 
 \\
 \cdot \lVert wyw^{-1} \rVert^{2a} \cdot T^{\eps+\frac{(n+4)(n-1)}{4}+ R\cdot\big(D(n)+\frac{n(n-1)}{2} \big) - \sum\limits_{j=1}^{n-1}B( a_j)}.
\end{multline*}

Recall that if $t=t(u)$ is as in \eqref{eq:t_inTermsOf_u}, if we define,  {for $a=(a_1,\ldots,a_{n-1}),b=(b_1,\ldots,b_{n-1})\in \R^{n-1}$,}
\begin{align*}
 K(w;a) := & \sum_{v\in V_n}\sum_{c_1=1}^\infty 
\cdots
\sum_{c_{n-1}=1}^\infty   \frac{\lvert S_w(\psi_L,\psi_M^v,c) \rvert}{ \prod\limits_{i=1}^{n-1} c_i^{1-2a_{i-1}+4a_i-2a_{i+1}}}, \\
X(w;a) := & \int\limits_{ {U}_w(\mathbb Z)\backslash  {U}_w(\mathbb R)}\;\int\limits_{\overline{{U}}_w(\mathbb R)} \delta(t)^{-\frac12}\cdot \lVert t \rVert^{2a} 
\;  d^* u,
\end{align*}
 {
and for a given choice of $\tau=(\tau_1,\ldots,\tau_{n-1})\in \{0,1\}^{n-1}$
\begin{align*}
Y(\tau,w;a,b) := & \int\limits_{y_1\in I_{\tau_1}} \int\limits_{y_2\in I_{\tau_2}} \cdots \int\limits_{y_{n-1}\in I_{\tau_{n-1}}}
\lVert y \rVert^{2b} \cdot \lVert wyw^{-1} \rVert^{2a} \; dy^\times,
\end{align*}}
then the bound on $\lvert \mathcal I_w {(\tau)}\rvert$ given in \eqref{eq:newerIw} can be replaced by
\begin{multline}\label{eq:Iw-bound-KXY}
 \lvert \mathcal I_w {(\tau)} \rvert \; \ll \; 
 T^{\eps+\frac{(n+4)(n-1)}{2}+ R\cdot\big(2D(n)+n(n-1) \big) - \sum\limits_{j=1}^{n-1}\big(B( a_j)+B(b_j)\big)}
 \\
 \cdot K(w;a) \cdot X(w;a) \cdot Y( {\tau,}w;a,b) \cdot \lVert L \rVert^{2a} \cdot \lVert M \rVert ^{2b}.
\end{multline}
We remark that in simplifying/finding $Y( {\tau,w;a,b})$, we have used \eqref{eq:deltaw}.  The basic strategy to prove Proposition~\ref{prop:Iwbounds} is now clear: we first find $a$ such that both $K(w;a)$ and $X(w;a)$ converge; then given this choice of $a$, we determine a particular value of $b$ for which $Y( {\tau,}w;a,b)$ converges as well; finally, we work out the corresponding bounds on $\lVert L \rVert^{2a}$, $\lVert M \rVert ^{2b}$ and $\sum\limits_{j=1}^{n-1} \big(B(a_j)+B(b_j)\big)$.

\subsection{Restrictions on the parameters $a$ and $b$}\label{sec:abconditions}

The trivial bound  {(see \cite{DR1998})} for the Kloosterman sum is given by
 \[ S(1,1,c)\ll \delta^{\frac12}(c) = c_1c_2\cdots c_{n-1}. \] 
Hence $K(w;a)$ is convergent whenever $a$ is chosen such that 
 \[  {\lVert c \rVert^{2a} = \prod_{k=1}^{n-1} c_{ {k}}^{2a_{k-1}-4a_k+2a_{k+1}} \ll \delta^{-\frac12-\eps}(c)} , \]
From \eqref{eq:delta-Vert}, if we set
$a_j = \frac{j(n-j)}{2}(1+\eps)$, then  {$\lVert c \rVert^{2a}=\delta^{-1-\eps}(c) \ll\delta^{-\frac12-\eps}(c) $.}  More generally, $K(w;a)$ converges in the following case:
\begin{equation}\label{eq:choice-of-a}
 \boxed{ a_j := \rho + \frac{j(n-j)}{2}(1+\eps), \quad  {\rho   >0},\ j=1,\ldots,n-1. } 
\end{equation}
That this choice of $a$ makes $K(w;a)$ converge is a consequence of the easily verifiable fact that
 \[  {\lVert c \rVert^{2a} = (c_1 c_{n-1})^{-2\rho} \cdot \delta^{-1-\eps}(c).} \]
We assume henceforth that $a$ satisfies \eqref{eq:choice-of-a}.

We next consider the convergence of $X(w;a)$.  Recall that the Iwasawa form of $wu$ is assumed to be $u_0 t k$, meaning $wu=u_0 tk$ where $u_0\in U_n(\R)$, $t\in T_n(\R)$ and $k\in O(n,\R)$.  Indeed, $t$ is given by \eqref{eq:t_inTermsOf_u}.  Then
 \[ X(w;a) = \int\limits_{ {U}_w(\mathbb Z)\backslash  {U}_w(\mathbb R)}\;\int\limits_{\overline{{U}}_w(\mathbb R)} \delta(t)^{-\frac12}\cdot \lVert t \rVert^{2a} 
\;  d^* u \ll \int\limits_{ {U}_w(\mathbb Z)\backslash  {U}_w(\mathbb R)}\;\int\limits_{\overline{{U}}_w(\mathbb R)} \delta^{-\frac32-\eps}(t) \; d^*u. \]
The fact that the right hand side converges is a consequence of Jacquet \cite{Jacquet1967}.  

We now turn to the convergence of $Y(\tau,w;a,b)$.
Applying Lemma~\ref{lem:yIwasawa-of-yprime} (which describes $\lVert wyw^{-1}\rVert^{2a}$), we see that
\begin{align*}
 Y( {\tau,w;a,b}) & 
 = \int\limits_{I_\tau} 
\left( \prod_{i=1}^s \prod_{j=1}^{n_i} y_{n-\hn_i+j}^{2b_{(n-\hn_i+j)}-2\big(a_{\hn_{i-1}}-a_{(\hn_{i-1}+j)}+a_{\hn_i}\big)} \right) dy^\times 
=  \prod_{i=1}^s \prod_{j=1}^{n_i}Y_{n-\hn_i+j}( {\tau,w;a,b}),
\end{align*}
where
 \[ Y_{n-\hn_i+j}( {\tau,w;a,b}) := \int\limits_{I_{\tau_{n-\hn_i+j}}} y_{n-\hn_i+j}^{2b_{n-\hn_i+j}-2(a_{\hn_{i-1}}-a_{\hn_{i-1}+j}+a_{\hn_i})}\; \frac{dy_{n-\hn_i+j}}{y_{n-\hn_i+j}} . \]
Hence, in order to bound $Y( {\tau,w;a,b})$ (and thereby show that $\mathcal I_w {(\tau)}$ converges), we must choose $b=(b_1,\ldots,b_{n-1})$ such that $Y_{n-\hn_i+j}( {\tau,w;a,b})$ converges.  Clearly
\begin{equation}\label{eq:choice-of-b}
 \boxed{ b_{n-\hn_i+j} = a_{\hn_{i-1}}-a_{\hn_{i-1}+j}+a_{\hn_i} + (-1)^{\tau_{{n-\hn_i+j}}}\cdot \tfrac{\eps}{2}, \quad\quad \mbox{($i=1,\ldots,s,\ j=1,\ldots,n_i$)} }
\end{equation}
suffices, since making this choice implies that, for each $k=1,\ldots,n-1$, {
 \[ Y_k( {\tau,w;a,b}) = \begin{cases} \displaystyle{\int_0^1} y^{\eps} \; \frac{dy}{y} & \mbox{ if }\tau_k=0, \\ \\ \displaystyle{\int_1^\infty} y^{-\eps} \; \frac{dy}{y} & \mbox{ if }\tau_k=1, \end{cases} \]}
which converges (and gives the same value $\frac{1}{\eps}$) in either case.

\subsection{Proof of Proposition~\ref{prop:Iwbounds}}\label{sec:Iwbounds-proof-final}

We have now shown that if $w=w_{(n_1,\ldots,n_r)}$ and we choose $a$ as in \eqref{eq:choice-of-a} and $b$ via \eqref{eq:choice-of-b} accordingly, the right hand side of \eqref{eq:Iw-bound-KXY} converges, hence gives a bound for $\lvert \mathcal I_w \rvert$.  Therefore, in order to complete the proof of Proposition~\ref{prop:Iwbounds}, we need to first show that 
 \[ \lVert L \rVert^{2a} \cdot \lVert M \rVert^{2b} \ll \prod_{i=1}^{n-1} \big( \ell_i m_i \big)^{2\rho+ \frac{n^2+1}{4}}, \]
and second that the given choice of $a$ and $b$ gives the claimed bound for the power of $T$ appearing in \eqref{eq:Iw-bound-KXY}.

 {To complete the first of these tasks we note that, by \eqref{eq:choice-of-a}  and the fact that $j(n-j)$ is maximized (in $j$) when $j=n/2$,  we have
\begin{equation}\label{eq:a-bound} a_j =\rho+\frac{j( n- j)}{2}(1+\eps) \le\rho+\frac{n^2}{8}(1+\eps)<\rho+\frac{n^2+1}{8} \end{equation}for $\eps<1/n^2$ and $1\le j\le n-1$. Similarly, using \eqref{eq:choice-of-a} and \eqref{eq:choice-of-b} we compute that, for $1\le i\le s$ and $1\le j\le n_i$,
 \begin{equation} b_{n-\hn_i+j} =\rho+\frac12\big(j^2+j(2 \hn_{i-1}-n)+\hn_i(n-\hn_i)\big) +\eps\label{eq:b-equation}\end{equation}for $\eps$ sufficiently small. Note that the right hand side of \eqref{eq:b-equation}is a concave up parabola in $j$, and therefore, on the interval $1\le j\le n_i$, can  attain its maximum only at $j=1$ or $j=n_i$. So, if we can show that $b_{n-\hn_i+1}$ and $b_{n-\hn_i+ni}$ both satisfy a suitable upper bound, then the same bound will hold for all $1\le j\le n_i$.}

   {We consider first the endpoint $j=n_i$. Using \eqref{eq:b-equation} and the fact that $\hn_i-n_i=\hn_{i-1}$, we find that
\[ b_{n-\hn_i+n_i} =\rho+\frac{1}{2}\hn_{i-1} (n -\hn_{i-1} )+\eps.\]Again, $j(n-j)$ is maximized when $j=n/2$, so we conclude that\begin{equation}\label{eq:b-bound1} b_{n-\hn_i+n_i} \le\rho+\frac{n^2}{8} +\eps<\rho+\frac{n^2+1}{8} \end{equation}for $\eps$ sufficiently small.}

 {Next we consider the endpoint $j=1$.  From \eqref{eq:b-equation} we find that
\begin{align}\label{eq:another-b} b_{n-\hn_i+1}& =\rho+\frac{1}{2}\big(1-n+\hn_i(n-\hn_i)+2\hn_{i-1}\big)+\eps\\&\le\rho+\frac{1}{2}\big(-1-n+\hn_i(n-\hn_i)+2\hn_{i}\big)+\eps,\nonumber\end{align} the last step because $\hn_{i-1}=\hn_i-n_i\le \hn_i-1$.  We find using calculus that, as a function of $\hn_i$,  the right hand side of \eqref{eq:another-b}   is maximized when $\hn_i=(n+2)/2$.  So 
\begin{align}\label{eq:b-bound2}  b_{n-\hn_i+1}& \le\rho+\frac{1}{2}\bigg(-1-n+\frac{n+2}{2}\bigg(n-\frac{n+2}{2}\bigg)+n+2\bigg)+\eps\\&=\rho+\frac{n^2}{8}+\eps\le\rho+\frac{n^2+1}{8} \nonumber\end{align}for $\eps$ small enough.  From \eqref{eq:b-bound1} and \eqref{eq:b-bound2} it follows, again, that \[b_{n-\hn_i+j} \le\rho+\frac{n^2+1}{8}\]for all $1\le i\le s$ and $1\le j\le n_i$. This and \eqref{eq:a-bound} yield the desired bound on $\lVert L \rVert^{2a} \cdot \lVert M \rVert^{2b}$.}
 
  {The second task is accomplished using Lemma \ref{lem:aplusb-bound}. \qed}

\section{\large\bf Bounding the Eisenstein spectrum \boldmath $\mathcal{E}$}\label{sec:Eisenstein-bound}
\vskip 5pt

Recall that if $L=(\ell_1,\ldots,\ell_{n-1})$, $M=(m_1,\ldots,m_{n-1})\in \Z^{n-1}$ with $\prod\limits_{i=1}^{n-1} \ell_i m_i\neq 0$, then, by Theorem \ref{ThmSpectralDecomp}, the Eisenstein contribution to the Kuznetsov trace formula is given by
\begin{align*}
 \mathcal{E} =  \sum_{\mathcal P } \sum_{\Phi} \mathcal{E}_{\mathcal{P},\Phi},
\end{align*}
where
\begin{align*}
 \mathcal{E}_{\mathcal{P},\Phi} :=  c_{\mathcal P}\hskip-4pt  \underset{\text{\rm Re}(s_j)=0}{\int\limits_{n_1s_1+\cdots +n_rs_{r}=0}}\hskip -11pt 
A_{E_{\mathcal P,\Phi}}(L, s) \, \overline{A_{E_{\mathcal P,\Phi}}(M, s)}\cdot\Big| p_{T,R}^{n,\#}\big(\alpha_{_{(\mathcal P,\Phi)}}(s)\big)\Big |^2\; ds.
\end{align*}
In this section we give bounds for $\mathcal{E}$ in the case that $L=(\ell,1,\ldots,1)$ and $M=(m,1,\ldots,1)$ with $\ell,m\neq 0$.

\subsection{The Eisenstein contribution $\mathcal{E}$ to the Kuznetsov trace formula} The main result of this section is the following.

\begin{theorem}[\bf{Bounding the Eisenstein contribution $\mathcal{E}$}]\label{th:EisensteinBound}
Fix $n\ge 2$ and a  sufficiently large integer $R>0$.  Let $L=(\ell,1,\ldots,1)$, $M=(m,1,\ldots,1) \in \Z^{n-1}$ with $\ell,m\neq 0$.  Then, assuming the Lower bound conjecture for Rankin-Selberg L-functions (see \eqref{RankinSelbergBound}), for $T\to\infty$ we have the bound
 \[ \sum_{\mathcal P}\sum_{\Phi}\left|\mathcal E_{\mathcal P,\Phi}\right| 
 \; \ll \;
   (\ell m)^{\frac12-\frac{1}{n^2+1}+\varepsilon} \cdot T^{R\cdot\left(\binom{2n}{n}-2^n \right)\; \;+\;\varepsilon}. \qedhere \]      
\end{theorem}

\subsection{Proof of Theorem 7.1.1}

\begin{proof} We proceed by induction on $n$, beginning with the case $n=2$.  In this case, the only parabolic subgroup is the minimum parabolic, or Borel, subgroup $\mathcal B=\mathcal P_{1,1}$, and the only function $\Phi$ corresponding to $\mathcal B$ (see Definition \ref{def:Hecke-Maassparabolic}) is the constant function $\Phi=1$.   The Eisenstein contribution in this case, then, is simply the quantity $\mathcal E_{\mathcal B,1}$.

By Theorem \ref{ThmSpectralDecomp} in the case $n=2$, we have
\begin{equation*}
\mathcal{E}_{\mathcal B,1}  =  \; c_{\mathcal B} \int\limits_{\re{ s_1}=0} 
A_{E_{\mathcal B,1}}(\ell, s) \, \overline{A_{E_{\mathcal B,1}}(m, s)}\cdot\Big| p_{T,R}^{2,\#}\big(\alpha_{_{(\mathcal B,1)}}(s)\big)\Big |^2\; ds_1,
\end{equation*}where $s=(s_1,-s_1)$.  Now note that, by \eqref{langlandsparamsforPhi}, $\alpha_{_{(\mathcal B,1)}}(s)=s$.  Moreover, by Definition \ref{def:FRn}, we have $\FR{2}  \equiv 1$, so by Definition \ref{PTR-and-hTR}, we have 
\begin{align*}
p_{T,R}^{2,\#}\big(\alpha_{_{(\mathcal B,1)}}(s)\bigr)\; = \;e^{{s_1^2}/{T^2}}
      \Gamma\left(\textstyle{\frac{2R+1+2s_1}{4}}\right) \Gamma\left(\textstyle{\frac{2R+1-2s_1}{4}}\right).
\end{align*}
Furthermore, we see from Proposition \ref{MthEisCoeff}  that
\begin{align*}
\left |A_{E_{\mathcal B,1}}\big( \ell ,s\big)\right| &\ll  \big|\zeta^*\big(1 + 2s_1 \big)\big|^{-1}  \underset{c_1c_2 =\ell} {\sum_{c_1, c_2\in\mathbb Z_{>0}}} \left|c_1^{\alpha_1} c_2^{\alpha_2} \right| \ll\;  { \ell^\varepsilon  \left|\Gamma\bigl({\textstyle{ \frac{1 + 2s_1}{2}}  }\bigr)\zeta\big(1+2s_1\big)  \right|^{-1}.}
\end{align*}
Then
\begin{align*}
\big|\mathcal E_{\mathcal B,1}\big| & \; \ll \; (\ell m)^\varepsilon \int\limits_{  \re(s_1)=0 } e^{{s_1^2}/{T^2}}
\frac{ \left|       \Gamma\left(\textstyle{\frac{2R+1+2s_1}{4}}\right) \Gamma\left(\textstyle{\frac{2R+1-2s_1}{4}}\right)\right|^2}
{\left|\Gamma\bigl({\textstyle{ \frac{1 + 2s_1}{2}}  }\bigr)\zeta\big(1+2s_1\big) \right|^2}\;|ds_1|.\end{align*}We may restrict our integration to the domain $|\im(s)|\le T$, since $e^{{s_1^2}/{T^2}}$ decays exponentially otherwise.  On this domain, we  use Stirling's bound  \eqref{eq:Stirling}
for the Gamma function, as well as the Vinogradov bound
$$|\zeta(1+it)|^{-1} \ll (1+|t|)^\varepsilon, \qquad (t\in\mathbb R).$$We get 
$$
\big|\mathcal E_{\mathcal B,1}\big|  \; \ll \;  (\ell m)^\varepsilon
\int\limits_{\substack{ \re(s_1)=0\\ \im(s_1) \le T}}
 \big|1+s_1   \big|^{2R-1+\varepsilon} \;|ds_1|,
$$from which it follows immediately that  $\big|\mathcal E_{\mathcal B,1}\big|   \; \ll \; T^{2R+\varepsilon}.$ So our desired result holds in the case $n=2$.
 
We now proceed to the general case.  For $n>2$, in order to establish bounds for $\mathcal{E}_{\mathcal{P},\Phi}$, we need to know that our main theorem is true for all $k<n$.  The reason this is needed is because we have to bound Rankin-Selberg L-functions  $L(s,\phi_{k}\times \phi_{k'})$ with $2\leq k,k'<n$.  This will require knowing the Weyl law with harmonic weights (Theorem \ref{th:WeylLaw2}) for $2\leq k,k'<n$.  We may assume by induction, however, that this is indeed the case, i.e., the Weyl law with harmonic weights holds for all $2\leq k<n$.

Now recall that, for the parabolic $\mathcal{P}$ associated to a partition $n=n_1+\cdots + n_r$,  we have
 $$\mathcal E_{\mathcal P,\Phi}=\underset{\text{\rm Re}(s_j)=0}{\int\limits_{n_1s_1+\cdots +n_rs_{r}=0}} \hskip-6pt 
A_{E_{\mathcal P,\Phi}}(L, s) \, \overline{A_{E_{\mathcal P,\Phi}}(M, s)}\cdot\Big| p_{T,R}^{n,\#}\big(\alpha_{_{(\mathcal P,\Phi)}}(s)\big)\Big |^2\; ds$$
where
 $\alpha_{_{\mathcal P,\Phi}}(s)$ is  given by (see  (\ref{langlandsparamsforPhi}))
\begin{align*}
\bigg (\overbrace{\alpha_{1,1}+s_1, \;\ldots \;,\alpha_{1,n_1}+s_1}^{n_1 \;\,\text{\rm terms}}, \quad\overbrace{\alpha_{2,1}+s_2, \;\ldots \;,\alpha_{2,n_2}+s_2}^{n_2 \;\,\text{\rm terms}}.
\quad\ldots \quad
,\overbrace{\alpha_{r,1}+s_r, \;\ldots \;,\alpha_{r,n_r}+s_r}^{n_r \;\,\text{\rm terms}}\bigg), 
\end{align*}
Since $\sum\limits_{i=1}^{n_k} \alpha_{k,i} = 0$ for all $1\le k\le r$
we see that
$\sum\limits_{k=1}^{r} \sum\limits_{i=1}^{n_k} (\alpha_{k,i}+s_k)^2 = \sum\limits_{k=1}^{r} \sum\limits_{i=1}^{n_k} (\alpha_{k,i}^2+s_k^2).$

Now, for any $\beta = (\beta_1, \ldots,\beta_n)\in\left(i\mathbb R \right)^n$ where $\hbeta_n=0$ we have
\begin{align*}
p_{T,R}^{n,\#}\big(\beta\big) \; := \;\text{\rm exp}\left(\frac{\beta_1^2+\beta_2^2+\cdots+\beta_n^2}{2T^2}\right)
\cdot \FR{n}(\tfrac{\beta}{2}) 
      \prod_{1\leq\, i < j\, \leq n}
     \left| \Gamma\left(\textstyle{\frac{2R+1+\beta_i - \beta_j}{4}}\right)\right|^2.
\end{align*}
It follows that
\begin{multline*}
p_{T,R}^{n,\#}\big(\alpha_{_{(\mathcal P,\Phi)}}(s)\big) \; =
\text{\rm exp}\left(\frac{\sum\limits_{k=1}^{r} \sum\limits_{i=1}^{n_k} \left(\alpha_{k,i}^2+s_k^2\right)}{2T^2}\right)
\FR{n}\left(\tfrac{\alpha_{_{(\mathcal P,\Phi)}}(s)}{2}\right) 
\\ \cdot
\underset{n_k\ne1}
{\prod_{k=1}^r} \prod_{1\leq\, i < j\, \leq n_k}
      \left|\Gamma\left(\textstyle{\frac{2R+1+\alpha_{k,i }- \alpha_{k,j}}{4}}\right)\right|^2
       \cdot \prod_{1\le k< k'\le r}\; \prod_{i=1}^{n_k} \;\prod_{j=1}^{n_{k'}} \left|\Gamma\left(\textstyle{\frac{2R+1+s_k-s_{k'}+\alpha_{k,i }- \alpha_{k',j}}{4}}\right)\right|^2.
\end{multline*}

By Proposition  \ref{MthEisCoeff}, the $m^{th}$ coefficient of $E_{\mathcal P, \Phi}$ is given by
 \begin{align*} 
  A_{E_{\mathcal P, \Phi}}\big((m,1,\ldots,1),s\big) & =  \underset{n_k\ne1}{\prod_{k=1}^r} L^*\big(1, \text{\rm Ad}\; \phi_k\big)^{-\frac12 }
  \prod_{1\le i<j\le r}   L^*\big(1+s_i-s_{j}, \;\phi_i\times\phi_{j}\big)^{-1}\\  
  &
  \hskip 50pt
   \cdot \underset {c_1c_2\cdots c_r = m } {\sum_{    1 \le c_1, c_2, \ldots, c_r \, \in \, \mathbb Z}} \hskip-10pt \lambda_{\phi_1}(c_1)  \cdots \lambda_{\phi_r}(c_r)
    \cdot c_1^{s_1}  \cdots c_r^{s_r}
 \end{align*} 
 up to a non-zero constant factor with absolute value depending  only on $n$. To bound the divisor sum above we will use the bound of Luo-Rudnick-Sarnak  \cite{LRS1999} for the $m^{th}$ Hecke Fourier coefficient of a $\GL(\kappa)$ (for $\kappa\ge 2)$ Hecke-Maass cusp form $\phi$ given by
 $$\left|\lambda_\phi(m,1,\ldots,1)\right| \le m^{\frac12-\frac{1}{\kappa^2+1}+\eps}.$$
 {(A slightly stronger result has been obtained by Kim and Sarnak \cite{MR1937203}.  However, the stated result above is sufficient for our purposes.)}  We immediately obtain the following bound for the divisor sum
$$\underset {c_1c_2\cdots c_r = m } {\sum_{    1 \le c_1, c_2, \ldots, c_r \, \in \, \mathbb Z}} \hskip-10pt \left|\lambda_{\phi_1}(c_1)  \cdots \lambda_{\phi_r}(c_r)
    \cdot c_1^{s_1}  \cdots c_r^{s_r}\right| \ll m^{\frac12-\frac{1}{n^2+1}+\varepsilon}.$$
    It follows that
  \begin{align*}
  \left|\mathcal E_{\mathcal P,\Phi}\right| \;
  & 
  \ll
 (m\ell)^{\frac12-\frac{1}{n^2+1}+\varepsilon} \cdot\text{\rm exp}\left(\frac{\sum\limits_{k=1}^{r} \sum\limits_{i=1}^{n_k} \alpha_{k,i}^2}{T^2}\right) 
\underset{\text{\rm Re}(s_j)=0, \;\text{\rm Im}(s_j)\ll T}{\int\limits_{n_1s_1+\cdots +n_rs_{r}=0}}   \left| \FR{n}\left(\tfrac{\alpha_{_{(\mathcal P,\Phi)}}(s)}{2}\right)\right|^2 
\\
 &
 \hskip -30pt 
\cdot \left(\underset{n_k\ne1}
{\prod_{k=1}^r} \prod_{1\leq\, i < j\, \leq n_k}
      \left|\Gamma\left(\textstyle{\frac{2R+1+\alpha_{k,i }- \alpha_{k,j}}{4}}\right)\right|^4\right)
       \cdot 
      \left( \prod_{1\le k< k'\le r}\; \prod_{i=1}^{n_k} \;\prod_{j=1}^{n_{k'}} \left|\Gamma\left(\textstyle{\frac{2R+1+s_k-s_{k'}+\alpha_{k,i }- \alpha_{k',j}}{4}}\right)\right|^4\right) \\
  &
  \hskip 17pt
 \cdot \left(\underset{n_k\ne1}{\prod_{k=1}^r} \big|L^*\big(1, \text{\rm Ad}\; \phi_k\big)\big|^{-1}\right) 
 \cdot  
 \left( \prod_{1\le k<k'\le r}  \big| L^*\big(1+s_k-s_{k'}, \;\phi_k\times\phi_{k'}\big)\big|^{-2} \;|ds|\right)\\
 &
 \\
 &
  \ll \;
 (m\ell)^{\frac12-\frac{1}{n^2+1}+\varepsilon}  \; 
 \underset{n_k\ne1} {\prod_{k=1}^r} \text{\rm exp}\left(\frac{  \alpha_{k,1}^2+\cdots+  \alpha_{k,n_k}^2}{T^2}\right) 
 \\
 &
 \cdot \hskip-16pt\underset{\text{\rm Re}(s_j)=0, \;\text{\rm Im}(s_j)\ll T}{\int\limits_{n_1s_1+\cdots +n_rs_{r}=0}}
 \hskip -10pt  \left| \FR{n}\left(\tfrac{\alpha_{_{(\mathcal P,\Phi)}}(s)}{2}\right)\right|^2 \;
 \underset{n_k\ne1}
{\prod_{k=1}^r} \frac{1}{ \left|L(1,\,\text{\rm Ad}\;\phi_k)\right|} \;\;
\prod_{1\leq\, i < j\, \leq n_k}
     \frac{ \left|\Gamma\left(\textstyle{\frac{2R+1+\alpha_{k,i }- \alpha_{k,j}}{4}}\right)\right|^4}{   \left| \Gamma\left(\frac{1+\alpha_{k,i}-\alpha_{k,j}}{2}   \right)\right|^2  }\\
     &
     \hskip 20pt
      \cdot 
      \prod_{1\le k< k'\le r}\frac{1}{\left|L\big(1+s_k-s_{k'}, \;\phi_k\times\phi_{k'}\big)\right|^2}\; \prod_{i=1}^{n_k} \;\prod_{j=1}^{n_{k'}} \frac{\left|\Gamma\left(\textstyle{\frac{2R+1+s_k-s_{k'}+\alpha_{k,i }- \alpha_{k',j}}{4}}\right)\right|^4}{       \left|\Gamma\left(\textstyle{\frac{1+s_k-s_{k'}+\alpha_{k,i }- \alpha_{k',j}}{2}}\right)\right|^2} \;|ds|.
  \end{align*}

\noindent  
\begin{lemma}\label{Lemma:mathcalF}  Assume $|s_k|\ll T^{1+\varepsilon}$ and $|\alpha_{k,j}|\ll T^{1+\varepsilon}$ for $1\le k\le r$ and $1\le j \le n_k$.  Then for $\alpha:=\alpha_{_{(\mathcal P,\Phi)}}(s)$ and $\alpha^{(k)}$  as in Definition \ref{partitioning-alpha}, we have
 $$\left| \FR{n}\big(\alpha_{_{(\mathcal P,\Phi)}}(s)\big)\right|^2 \ll \Bigg(\underset{n_k\ne1}{\prod_{k=1}^r} \left| \FR{n_k}\big(\alpha^{(k)}\big)\right|^2\Bigg) \cdot T^{R\cdot  B(n)+\varepsilon}$$
 where $B(n) = 2D(n) -2 \underset{n_k\ne 1}{\sum\limits_{k=1}^r} D(n_k).$
 \end{lemma}
 \begin{proof} This follows immediately from Lemma
 \ref{lem:FRnkdecomp}.
 \end{proof}

It follows from  Lemma \ref{Lemma:mathcalF} that  for  $|\alpha^{(k)}|^2 = \alpha_{k,1}^2+\cdots+  \alpha_{k,n_k}^2$, we have
 {\small \begin{align*} \left|\mathcal E_{\mathcal P,\Phi}\right| & \ll (m\ell)^{\frac12-\frac{1}{n^2+1}+\varepsilon}  \; 
   T^{R\cdot B(n)+\varepsilon}
 \underset{n_k\ne1} {\prod_{k=1}^r} \frac{\text{\rm exp}\left(\frac{ |\alpha^{(k)}|^2}{T^2}\right)
 \left| \FR{n_k}\left(\tfrac{\alpha^{(k)}}{2}\right)\right|^2
\prod\limits_{1\leq\, i < j\, \leq n_k}
     \frac{ \left|\Gamma\left(\textstyle{\frac{2R+1+\alpha_{k,i }- \alpha_{k,j}}{4}}\right)\right|^4}{   \left| \Gamma\left(\frac{1+\alpha_{k,i}-\alpha_{k,j}}{2}   \right)\right|^2  }}{\left|L(1,\,\text{\rm Ad}\;\phi_k)\right|}\\
     &
     \hskip-10pt
     \cdot \hskip-16pt \underset{\text{\rm Re}(s_j)=0, \;|\text{\rm Im}(s_j)|\ll T}{\int\limits_{n_1s_1+\cdots +n_rs_{r}=0}}\prod_{1\le k< k'\le r}\frac{1}{\left|L\big(1+s_k-s_{k'}, \;\phi_k\times\phi_{k'}\big)\right|^2}\; \prod_{i=1}^{n_k} \;\prod_{j=1}^{n_{k'}} \frac{\left|\Gamma\left(\textstyle{\frac{2R+1+s_k-s_{k'}+\alpha_{k,i }- \alpha_{k',j}}{4}}\right)\right|^4}{       \left|\Gamma\left(\textstyle{\frac{1+s_k-s_{k'}+\alpha_{k,i }- \alpha_{k',j}}{2}}\right)\right|^2} \;|ds|. 
     \\
     &
     \\
     &
     \\
     &
     \ll \; (m\ell)^{\frac12-\frac{1}{n^2+1}+\varepsilon} \cdot T^{R\cdot B(n)+\varepsilon} \;\; 
 \underset{n_k\ne1} {\prod_{k=1}^r}  \;\;\frac{ \left|h_{T,R}^{(n_k)}(\alpha^{(k)})\right|}{\left|L(1,\,\text{\rm Ad}\;\phi_k)\right| } \\
 &
 \hskip-20pt
     \cdot  \hskip-10pt\underset{\text{\rm Re}(s_j)=0, \;|\text{\rm Im}(s_j)|\ll T}{\int\limits_{n_1s_1+\cdots +n_rs_{r}=0}}\prod_{1\le k< k'\le r}\frac{1}{\left|L\big(1+s_k-s_{k'}, \;\phi_k\times\phi_{k'}\big)\right|^2} \;\,\prod_{i=1}^{n_k} \;\prod_{j=1}^{n_{k'}} \frac{\left|\Gamma\left(\textstyle{\frac{2R+1+s_k-s_{k'}+\alpha_{k,i }- \alpha_{k',j}}{4}}\right)\right|^4}{       \left|\Gamma\left(\textstyle{\frac{1+s_k-s_{k'}+\alpha_{k,i }- \alpha_{k',j}}{2}}\right)\right|^2} \;|ds|. 
  \end{align*}   }
  \vskip 10pt\noindent
  where
  \vskip -5pt
  $$h_{T,R}^{(n_k)}\left(\alpha^{(k)}\right) = \text{\rm exp}\left(\frac{|\alpha^{(k)}|^2}{T^2}\right) \FR{n_k}\big(\tfrac{\alpha^{(k)}}{2}\big)^2\prod_{1\leq\, i \ne j\, \leq n_k}
     \frac{\Gamma\left(\textstyle{\frac{2R+1+\alpha_{k,i }- \alpha_{k,j}}{4}}\right)^2}{    \Gamma\left(\frac{1+\alpha_{k,i}-\alpha_{k,j}}{2}   \right)  }.$$
 
 Next
 $$\prod_{1\le k< k'\le r}\; \prod_{i=1}^{n_k} \;\prod_{j=1}^{n_{k'}} \frac{\left|\Gamma\left(\textstyle{\frac{2R+1+s_k-s_{k'}+\alpha_{k,i }- \alpha_{k',j}}{4}}\right)\right|^4}{       \left|\Gamma\left(\textstyle{\frac{1+s_k-s_{k'}+\alpha_{k,i }- \alpha_{k',j}}{2}}\right)\right|^2} 
\; \ll \;
 T^{(2R-1)\hskip-5pt\sum\limits_{1\le k<k'\le r} n_k\cdot n_{k'}}.$$

 We obtain the bound
  \begin{align*} \left|\mathcal E_{\mathcal P,\Phi}\right| & 
   \ll 
   (m\ell)^{\frac12-\frac{1}{n^2+1}+\varepsilon} \cdot T^{R\cdot B(n)+\varepsilon + (2R-1)\;\cdot\hskip-5pt\sum\limits_{1\le k<k'\le r} \hskip-5pt n_k\cdot n_{k'}} \; 
\cdot  \prod_{k=1}^r  \frac{ \left|h_{T,R}^{(n_k)}\left(\alpha^{(k)}\right)\right|}{\left|L(1,\,\text{\rm Ad}\;\phi_k)\right| } 
 \\
     &
     \hskip 130pt
     \cdot
      \hskip-16pt \underset{\text{\rm Re}(s_j)=0, \;|\text{\rm Im}(s_j)|\ll T}{\int\limits_{n_1s_1+\cdots +n_rs_{r}=0}}\prod_{1\le k< k'\le r}
      \frac{|ds|}{\left|L\big(1+s_k-s_{k'}, \;\phi_k\times\phi_{k'}\big)\right|^2}.
  \end{align*}
  Next, we bound the $s$-integral above. It follows from Langlands conjecture  (see \ref{RankinSelbergBound})
 that  for $|\text{\rm Im}(s_k)|,|\text{\rm Im}(s_{k'})| \ll T$ we have the bound   $$\left|L\big(1+s_k-s_{k'}, \;\phi_k\times\phi_{k'}\big)\right|^{-2} \ll T^\varepsilon.$$ 
  This together with the bound 
   \[ \int\limits_{\substack{n_1s_1+\cdots +n_rs_r=0 \\ \re(s_j)=0,\ \lvert \im(s_j)\rvert \ll T}} \hskip -14pt |ds| \ \ll \ T^{r-1}, \]
   implies that
   \begin{equation} \label{EP-Phi-Bound2}
 \left|\mathcal E_{\mathcal P,\Phi}\right|  
   \ll 
   (m\ell)^{\frac12-\frac{1}{n^2+1}+\varepsilon}  \,T^{R\cdot B(n) \,+\, (2R-1)\hskip-8pt\sum\limits_{1\le k<k'\le r} \hskip -5pt n_k\cdot n_{k'} \; +(r-1)\, +\,\varepsilon} \; 
  \cdot\left(\prod_{k=1}^r  \frac{ \left|h_{T,R}^{(n_k)}\left(\alpha^{(k)}\right)\right|}{\left|L(1,\,\text{\rm Ad}\;\phi_k)\right| }\right). 
\end{equation} 
\vskip 5pt
Since each $n_k<n$ (for $k=1,2,\ldots,r)$,  we can apply our inductive procedure together with the following theorem to bound  $\sum_{\Phi}\left|\mathcal E_{\mathcal P,\Phi}\right| .$

 \vskip 5pt
\begin{theorem}[\bf{Weyl law with harmonic weights for \boldmath $ {\GL(n_k)}$ with ${n_k<n}$}]\label{th:WeylLaw2}
Suppose $n_k\in\mathbb Z$ with $2\leq n_k <n$.  Let $\{\phi_1,\phi_2,\ldots\}$ be an orthogonal basis of Hecke-Maass cusp forms for $\GL(n_k)$ ordered by eigenvalue.  If $\alpha^{(j)}$ are the Langlands parameters of $\phi_j$, then
\begin{equation}\label{eq:WeylLaw2}
 \sum_{j=1}^\infty  \frac{ h_{T,R}^{(n_k)}\left(\overline{\alpha^{(j)}}\right)}{\mathcal{L}_j }\ll_n \; T^{2R  \cdot\big(D(k)+\frac{n_k(n_k-1)}{2}\big)+ n_k-1}.
\end{equation}
\end{theorem}
\vskip 2pt
\begin{proof}
In \cite{GSW21}, all that was needed to prove this statement for $n=4$ was to have it be true for $n_k=2$ and $n_k=3$, which was already known.  A similar induction argument works in general.
\end{proof}
It immediately follows from the bound \eqref{EP-Phi-Bound2} and \eqref{eq:WeylLaw2} that
$$
 \sum_{\Phi}\left|\mathcal E_{\mathcal P,\Phi}\right|  
   \ll 
   (m\ell)^{\frac12-\frac{1}{n^2+1}+\varepsilon}  \,T^{R\cdot B(n) \,+\, (2R-1)\hskip-8pt\sum\limits_{1\le k<k'\le r} \hskip -5pt n_k\cdot n_{k'} \; +(r-1)\, +\,\varepsilon} 
   \hskip-5pt 
  \cdot  \;T^{\sum\limits_{k=1}^r\left(2R  \cdot\big(D(k)+\frac{n_k(n_k-1)}{2}\big)+ n_k-1\right)}. 
$$
Recall that
$B(n) = 2D(n) -2\sum\limits_{k=1}^r D(n_k),$ which implies that
$$ 
 \sum_{\Phi}\left|\mathcal E_{\mathcal P,\Phi}\right|  
   \ll 
   (m\ell)^{\frac12-\frac{1}{n^2+1}+\varepsilon}  \,T^{2R\cdot D(n) \,+\, 2R\left(\sum\limits_{1\le k<k'\le r} \hskip -5pt n_k\cdot n_{k'} + \sum\limits_{k=1}^r   \frac{n_k(n_k-1)}{2}\right)+ \sum\limits_{k=1}^rn_k\;\, - 
   \hskip-4pt\sum\limits_{1\le k<k'\le r} \hskip -5pt n_k\cdot n_{k'} \; -1\, +\,\varepsilon}. 
$$
Next, $\sum\limits_{1\le k <k'\le r} \hskip -5pt n_k\cdot n_{k'} \;+ \sum\limits_{k=1}^r \frac{n_k(n_k-1)}{2}  =  \frac{n(n-1)}{2}$ by Lemma \ref{combinatorial} and $\sum\limits_{k=1}^r n_k =n.$ 
It follows that
$$ 
 \sum_{\Phi}\left|\mathcal E_{\mathcal P,\Phi}\right|  
   \ll 
   (m\ell)^{\frac12-\frac{1}{n^2+1}+\varepsilon}  \,T^{2R\cdot \left(D(n) \,+\,\frac{n(n-1)}{2}\right) + n-1\;- \hskip-5pt\sum\limits_{1\le k<k'\le r} \hskip -5pt n_k\cdot n_{k'} \, +\,\varepsilon}. 
$$
To complete the proof, we need to sum over all parabolics $\mathcal P$.  It suffices, therefore, to consider the ``worst case scenario'' among the possible partitions $n=n_1+\cdots+n_r$ for which the expression $$\sum\limits_{1\le k<k'\le r}n_k n_{k'}$$ is minimized. It is easy to see that this occurs when $r=2$ and $\{n_1,n_2\} = \{ n-1,1\}$, giving the bound $n-1$.
It follows that
 \[ \sum_{\mathcal P}\sum_{\Phi}\left|\mathcal E_{\mathcal P,\Phi}\right| 
\ll 
   (m\ell)^{\frac12-\frac{1}{n^2+1}+\varepsilon}  \,T^{2R\cdot \left(D(n) \,+\,\frac{n(n-1)}{2}\right)\, +\,\varepsilon}. \]
Using \eqref{eq:Dn-explicit}, this immediately implies the desired result.
\end{proof}

 {
\begin{remark}
In \cite{2019arXiv191101880J} Jana and Nelson prove the bound
\begin{equation}\label{eq:JanaNelsonBound}
 \sum_{c\left(\phi_j \right)\le T^n} \frac{1}{\mathcal L_j} \; \ll \; T^{n^2-n},
\end{equation}
where $c(\phi)$ is the analytic conductor given in \eqref{eq:analyticcond}.  This is an unsmoothed version of Theorem~\ref{th:WeylLaw2}.  Our result is a smoothed version, and it doesn't seem possible to derive a bound as in \eqref{eq:JanaNelsonBound} with a sharp cutoff without using a different approach.
\end{remark}}

\section{\large\bf An integral representation of \boldmath $p_{T,R}^{(n)}(y)$} \label{sec:IntRepOfpTR}

Recall (see \eqref{PTR-and-hTR}) that
\[
p_{T,R}^{n,\#}(\alpha) \; := \;e^{\frac{\alpha_1^2+\alpha_2^2+\cdots+\alpha_n^2}{2T^2}}
\cdot \FR{n}(\tfrac{\alpha}{2}) 
      \prod_{1\leq\, j \ne k\, \leq n}
      \Gamma\left(\textstyle{\frac{1+2R+\alpha_j - \alpha_k}{4}}\right).
\]
Using the formula for the inverse Lebedev-Whittaker transform given in \cite{MR2982424}, it follows that
\begin{align*}
 p_{T,R}^{(n)}(y) & := \frac{1}{\pi^{n-1}}
  \int\limits_{\substack{\halpha_n=0\\\re(\alpha)=0}} 
     \frac{ p^{n,\#}_{T,R}(\alpha) \;\overline{W_{n,\alpha}(y)}
         }{\prod\limits_{1\leq\, j \ne k\, \leq n} \Gamma\left(\frac{\alpha_j-\alpha_k}{2}\right) } \; d\al
    \\ & =
 \frac{1}{\pi^{n-1}}
  \int\limits_{\substack{\halpha_n=0\\\re(\alpha)=0}} e^{\frac{\alpha_1^2+\alpha_2^2+\cdots+\alpha_n^2}{2T^2}}
\cdot \FR{n}(\tfrac{\alpha}{2}) 
      \prod_{1\leq\, j \ne k\, \leq n}
      \Gamma_R\big(\tfrac{\alpha_j-\alpha_k}{2}\big) \;\overline{W_{n,\alpha}(y)}
     \ d\alpha,
\end{align*}
where $\Gamma_R(z):= \frac{\Gamma\left(\frac{\frac12+R+z}{2}\right)}{\Gamma(z)}$.

The strategy in this section for giving a representation of $p_{T,R}^{(n)}(y)$ follows the same general outline as was used to obtain the results for $\GL(3)$ and $\GL(4)$ given in the papers \cite{GK2013} and \cite{GSW21}, respectively. 
As in the prior works, we express the Whittaker function as the inverse Mellin transform of its Mellin transform.  (See Section~\ref{sec:NormalizedWhitt}.)  Plugging this into the above formula gives an integral representation of $p_{T,R}^{(n)}(y)$ in terms of an additional variable $s=(s_1,\ldots,s_{n-1})$.  

 \vskip 8pt

\subsection{Normalized Mellin transform of Whittaker function}\label{sec:NormalizedWhitt}

We introduce (as in \cite{IS2007}) the following Mellin transform and its inverse.
\begin{definition}[\bf{Normalized Mellin transform of Whittaker function}]
Let $n\in \Z_+$ and $\al=(\alpha_1,\ldots,\alpha_n)\in \C^n$ such that $\halpha_n=0$.  Let $W_{n,\alpha}(y)$ be the Whittaker function of Definition~\ref{def:JacWhittFunction}.  The \emph{Mellin transform} is
\begin{equation}\label{eq:MellinOfGLnWhit}
 \what{n}{\al}{s}: = 2^{n-1}\int_0^\infty\cdots \int_0^\infty W_{n,2\alpha}(y) \prod_{j=1}^{n-1}(\pi y_j)^{2s_j} \frac{dy_j}{y_j^{1+\frac{j(n-j)}{2}}}, 
\end{equation}
and the \emph{inverse Mellin transform} is given by
\begin{equation}\label{eq:WhitAsMellinTransform}
 W_{n,\alpha}(y) = \frac{1}{2^{n-1}} \int\limits_{\substack{s=(s_1,\ldots,s_{n-1})\\\re(s)=2b}} 
 \bigg(\prod_{j=1}^{n-1}
 y_j^{\frac{j(n-j)}{2}}(\pi y_j)^{-s_j} \bigg)
 \widetilde{W}_{n,\frac{\alpha}{2}}(\tfrac{s}{2})\, ds.
\end{equation}
\end{definition}

As a consequence of this definition, we have
\begin{align}\label{eq:pTRnstart}
 p_{T,R}^{(n)}(y) = &
\frac{1}{(2\pi)^{n-1}} \int\limits_{\substack{\halpha_n=0\\ \re(\alpha)=0}} 
 e^{\frac{\alpha_1^2+\cdots+\alpha_n^2}{T^2/2}} 
 \cdot 
 \FR{n}(\alpha)
 \bigg(\prod_{1\leq j\neq k \leq n} \hskip -4pt
 \Gamma_R\big(\alpha_j-\alpha_k\big) \bigg)
 \\ \nn & \hskip 96pt \cdot \hskip -8pt
 \int\limits_{\substack{s=(s_1,\ldots,s_{n-1})\\\re(s)=b}} \hskip -8pt
 \bigg(\prod_{j=1}^{n-1}
 y_j^{\frac{j(n-j)}{2}}(\pi y_j)^{-2s_j} \bigg)
 \widetilde{W}_{n,\alpha}(s)\, ds\, d\alpha,
\end{align}
where $b=(b_1,\ldots,b_{n-1})$ with each $b_j>0$.

We use the following theorem to make \eqref{eq:pTRnstart} explicit and to begin setting up an inductive method to bound $p_{T,R}^{(n)}(y)$ for all $n\geq 2$.

\begin{theorem}[{\bf Ishii-Stade}]\label{th:IS}
Let $m\geq 2$ and $\eps>0$.  Fix a Langlands parameter $\alpha\in\C^m$. Let $s\in \C^{m-1}$ with $\re(s)>\eps$.  Then
\begin{multline}\label{eq:ISrecursion}
 \widetilde{W}_{m,\al}(s) 
 =\int\limits_{\substack{z=(z_1,\ldots,z_{m-2})\\ \re(z)=\eps }} \left( \prod_{j=1}^{m-1} \Gamma\Big(s_j-z_{j-1}+\frac{(m-j)\alpha_m}{m-1}\Big)\Gamma\Big(s_j-z_j-\frac{j\alpha_m}{m-1}\Big)\right)
 \cdot
\frac{ \widetilde{W}_{m-1,\beta}\left(z \right)}{(2\pi i)^{m-2} }\, dz,
\end{multline}
 where
 \[ z_0:=-0+\frac{0\cdot \alpha_m}{m-1}=0, \quad z_{m-1}:=\alpha_m-\frac{(m-1)\alpha_m}{m-1}=0, \]
and
 \[ \beta = \big(\beta_1,\ldots,\beta_{m-1}\big) := \left( \alpha_1+\frac{\alpha_m}{m-1},\ldots,\alpha_{m-1} + \frac{\alpha_m}{m-1}\right) . \]
\end{theorem}

\subsection{A shifted $p_{T,R}^{(n)}$ term and the Ishii-Stade Conjecture}\label{sec:IS-Conj}

Our goal is to insert \eqref{eq:ISrecursion} into \eqref{eq:pTRnstart} and then shift the lines of integration in $s$ to $\re(s)=-a$, to the left of some of the poles of $\what{n}{\al}{s}$, which (see Theorem~\ref{th:residue-description}) occur at $\re(s_i)=-\delta$ for every $1\leq i \leq n-1$ and $\delta\in \Z_{\geq 0}$.  By Cauchy's residue formula, this allows us to describe $p_{T,R}^{(n)}(y)$ in terms of a the sum of a \emph{shifted $p_{T,R}^{(n)}$ term} and finitely many \emph{shifted residue terms}.

\begin{definition}[\bf{shifted $p_{T,R}^{(n)}$ term}]
\label{def:shiftedpTR} Let $n \geq 2$ be an integer and $a=(a_1,\ldots,a_{n-1})\in \R^{n-1}$.  The \emph{shifted $p_{T,R}^{(n)}$ term} is given by the same formula as \eqref{eq:pTRnstart} but with $b$ replaced by $-a$:
\begin{align}\label{eq:shifted-pTR}
 p_{T,R}^{(n)}(y;-a) := &
\frac{1}{(2\pi)^{n-1}} \int\limits_{\substack{\halpha_n=0\\ \re(\alpha)=0}} 
 e^{\frac{\alpha_1^2+\cdots+\alpha_n^2}{T^2/2}} 
 \cdot 
 \FR{n}(\alpha)
 \bigg(\prod_{1\leq j\neq k \leq n} \hskip -4pt
 \Gamma_R\big(\alpha_j-\alpha_k\big) \bigg)
 \\ \nn & \hskip 96pt \cdot \hskip -8pt
 \int\limits_{\substack{s=(s_1,\ldots,s_{n-1})\\\re(s)=-a}} \hskip -8pt
 \bigg(\prod_{j=1}^{n-1}
 y_j^{\frac{j(n-j)}{2}}(\pi y_j)^{-2s_j} \bigg)
 \widetilde{W}_{n,\alpha}(s)\, ds\, d\alpha.
\end{align}
\end{definition}

One might be tempted to insert \eqref{eq:ISrecursion} into \eqref{eq:shifted-pTR}, but this is invalid if $n>3$, because Theorem~\ref{th:IS} requires that $\re(s_i)>\eps$ for each $i=1,\ldots,n-1$.  To overcome this difficulty, we use shift equations as given in the following conjecture.  This allows us to evaluate $\what{n}{\al}{s}$ for $\re(s)<0$.

\begin{conjecture}[{\bf Ishii-Stade}]\label{conj}
Let $m,n\in \Z$ with $1\le m\le n-1$; let $\delta\in \Z_{\ge0}$.  Let $(x)_n:=\frac{\Gamma(x+n)}{\Gamma(x)}=x(x+1)\cdots(x+n-1)$.  Then there exists a positive integer $r$ and, for  each $i$ with $1\le i\le r$,  a polynomial $P_i(s,\al)$ and an $(n-1)$-tuple $\Sigma_i\in(\Z_{\ge0})^{n-1}$, such that
\begin{align}&\what{n}{\al}{s} =\biggl[\prod_{1\le j_1<j_2<\ldots<j_m\le n}(s_m+\al_{j_1}+\al_{j_2}+\cdots+\al_{j_m})_{_{\scriptstyle{\delta}}}\biggr]^{-1}\sum_{i=1}^r P_i(s,\al)\what{n}{\al}{s+\Sigma_i}, \label{s1-shift-n1}\end{align}where the $m$th coordinate of each $\Sigma_i$ is at least $\delta$.  Moreover, for each $i$, we have
 \[ \deg(P_i(s,\al))+2\big|\Sigma_i\big|=\delta \binom{n}{m}. \]
\label{thislemma}
\end{conjecture}

\begin{proof}[Proof of conjecture for $2\le n\le 5$]
Note that the case $\delta=0$ of the conjecture is trivial. Moreover, for a given $m$ and $n$ with $1\le m\le n-1$, it's enough to prove the conjecture for $\delta=1$.  The case   $\delta>1$ then follows by applying the case $\delta=1$ to itself iteratively.

 For $\delta=1$ and $n=2$ or $n=3$, the conjecture follows immediately from the explicit formulae
\begin{align*}\what{2}{\al}{s}&=\Gamma(s+\al)\Gamma(s-\al);\\\what{3}{\al}{s}&=\frac{\Gamma(s_1+\al_1)\Gamma(s_1+\al_2)\Gamma(s_1+\al_3)\Gamma(s_2-\al_1)\Gamma(s_2-\al_2)\Gamma(s_2-\al_3)}{\Gamma(s_1+s_2)} ,\end{align*}respectively, together with the functional equation $\Gamma(s+1)=s\Gamma(s).$ The case $\delta=1$ and $n=4$ is a consequence of \cite[equations (21), (29), and  
(31)]{MR4166593}.

We now consider the case $\delta=1$ and $n=5$. Note that it suffices to derive the appropriate recurrence relations for $m=1$ and $m=2$ (that is, for the variables $s_1$ and $s_2$); the cases $m=3$ and $m=4$ then follow from the invariance of $\what{5}{\al}{s}$ under the involution$$(s_1,s_2,s_3,s_4,\al_1,\al_2,\al_3,\al_4,\al_5)\to(s_4,s_3,s_2,s_1,-\al_1,-\al_2,-\al_3,-\al_4,-\al_5).$$

We follow an approach developed by Taku Ishii (personal correspondence).  First, consider the case $m=1$: we wish to show that

\begin{align}&\biggl[\prod_{i=1}^{5}(s_1+\al_i)\biggr]\what{5}{\al}{s}  \label{s1-shift-5}\end{align}is equal to a finite sum of terms  $P_i(s,\al)\what{n}{\al}{s+\Sigma_i},$ where the first coordinate of each $\Sigma_i\in (\Z_{\ge0})^{4}$ is at least one, and  $\deg(P_i(s,\al))+2\big|\Sigma_i\big|=5$ for each $i$.
To this end, let \begin{equation}\sigma=(\sigma_1,\sigma_2,\sigma_3,\sigma_4,\sigma_5)=(-s_1,s_1-s_2,s_2-s_3,s_3-s_4,s_4);\label{sigdef}\end{equation}note that $\sum_i \sigma_i=0$.  Since $s_1+\sigma_1=0$, we have
\begin{equation}  \biggl[\prod_{i=1}^{5}(s_1+\al_i)\biggr]\what{5}{\al}{s} =\biggl[\prod_{i=1}^{5}(s_1+\al_i)-\prod_{i=1}^{5}(s_1+\sigma_i)\biggr]\what{5}{\al}{s}. \label{addzero}\end{equation}

But for indeterminates $T,x_1,x_2,x_3,x_4,x_5$, we have
  
\begin{align}  \prod_{i=1}^{5}(T+x_i)=T^5+T^4P_1(x)+T^3P_2(x)+T^2 P_3(x)+TP_4(x)+P_5(x),\label{sym-decomp41}\end{align}where $P_k(x)$ is the elementary symmetric polynomial of degree $k$ in $x_1,x_2,x_3,x_4,x_5$.  So by equation \eqref{addzero}  above, we have
     
      \begin{align}\label{first41} \biggl[\prod_{i=1}^5(s_1+\al_{i} )\biggr] \what{5}{\al}{s}&= \bigl[s_1^5+s_1^4P_1(\al)+s_1^3P_2(\al)+s_1^2  P_3(\al)  + s  P_4(\al) +  P_5(\al)  \bigr]\what{5}{\al}{s}\\&- \bigl[s_1^5+s_1^4P_1(\sigma)+s_1^3P_2(\sigma)+s_1^2  P_3(\sigma)  + s  P_4(\sigma) +  P_5(\sigma)  \bigr]\what{5}{\al}{s}\nr&= \bigl[s_1^3\{P_2(\al)-P_2(\sigma) \} +s_1^2\{P_3(\al) - P_3(\sigma)\}\nr&\hskip17pt +s_1 \{ P_4(\al) - P_4(\sigma)\}+ \{ P_5(\al) - P_5(\sigma)\}\bigr]\what{5}{\al}{s},\nn  \end{align}since $P_1(\al)=P_1(\sigma)=0$.  
      
Now let $e_k$, for $1\le k\le 4$, be the four-tuple with a one in the $k$th place and zeroes elsewhere.  By \cite[Proposition 3.6]{MR3146819}, we have
$$P_k(\al)-P_k(\sigma)=Z_k-C_k$$(as operators acting on functions in the variable $s=(s_1,s_2,s_3,s_4)$), where the ``Capelli elements'' $C_k$ annihilate $\what{5}{\al}{s}$, and

\begin{align*}Z_2f(s)&= f(s+e_1)+f(s+e_2)+f(s+e_3)+f(s+e_4);\\
Z_3f(s)&=P_1(\sigma_3,\sigma_4,\sigma_5) f(s+e_1)+P_1(\sigma_1,\sigma_4,\sigma_5) f(s+e_2)+P_1(\sigma_1,\sigma_2,\sigma_5) f(s+e_3)\nr&+P_1(\sigma_1,\sigma_2,\sigma_3) f(s+e_4);\nr
Z_4f(s)&=P_2(\sigma_3,\sigma_4,\sigma_5) f(s+e_1)+P_2(\sigma_1,\sigma_4,\sigma_5) f(s+e_2)+P_2(\sigma_1,\sigma_2,\sigma_5) f(s+e_3)\nr&+P_2(\sigma_1,\sigma_2,\sigma_3) f(s+e_4)+f(s+e_1+e_3)+f(s+e_1+e_4)+f(s+e_2+e_4)\nr
Z_5f(s)&=P_3(\sigma_3,\sigma_4,\sigma_5) f(s+e_1)+P_3(\sigma_1,\sigma_4,\sigma_5) f(s+e_2)+P_3(\sigma_1,\sigma_2,\sigma_5) f(s+e_3)\nr&+P_3(\sigma_1,\sigma_2,\sigma_3) f(s+e_4)+P_1(\sigma_5)f(s+e_1+e_3)+P_1(\sigma_3)f(s+e_1+e_4)\nr&+P_1(\sigma_1)f(s+e_2+e_4);\nr.\nn \end{align*}So by \eqref{first41},
      \begin{align} \label{second-41}\biggl[\prod_{i=1}^{5}(s_1+\al_{i} )\biggr] \what{5}{\al}{s}&= \bigl[s_1^3 Z_2+s_1^2Z_3+s_1 Z_4+Z_5 \bigr]\what{5}{\al}{s}
    \\ &= \bigl[s_1^3 +s_1^2 P_1(\sigma_3,\sigma_4,\sigma_5)+ s_1P_2(\sigma_3,\sigma_4,\sigma_5) +P_3(\sigma_3,\sigma_4,\sigma_5)  \bigr]\what{5}{\al}{s+e_1}
      \nr  &+ \bigl[s_1^3 +s_1^2 P_1(\sigma_1,\sigma_4,\sigma_5)+ s_1P_2(\sigma_1,\sigma_4,\sigma_5) +P_3(\sigma_1,\sigma_4,\sigma_5)  \bigr]\what{5}{\al}{s+e_2}
        \nr  &+\bigl[s_1^3 +s_1^2 P_1(\sigma_1,\sigma_2,\sigma_5)+ s_1P_2(\sigma_1,\sigma_2 ,\sigma_5) +P_3(\sigma_1,\sigma_2,\sigma_5)  \bigr]\what{5}{\al}{s+e_3}
          \nr  &+\bigl[s_1^3 +s_1^2 P_1(\sigma_1,\sigma_2,\sigma_3)+ s_1P_2(\sigma_1,\sigma_2,\sigma_3) +P_3(\sigma_1,\sigma_2,\sigma_3)  \bigr]\what{5}{\al}{s+e_4}
          \nr  &+\bigl[ s_1 +P_1(\sigma_5)  \bigr]\what{5}{\al}{s+e_1+e_3}
              \nr  &+\bigl[ s_1 +P_1(\sigma_3)  \bigr]\what{5}{\al}{s+e_1+e_4}
                  \nr  &+\bigl[ s_1 +P_1(\sigma_1)  \bigr]\what{5}{\al}{s+e_2+e_4}.\nn\end{align}Recalling that the $P_k$'s are the elementary symmetric polynomials of degree $k$ in their arguments, we see that
          $$s_1^3 +s_1^2 P_1(a,b,c)+ s_1P_2(a,b,c) +P_3(a,b,c) =(s_1+a)(s_1+b)(s_1+c),$$for indeterminates $a, b,c$.  So \eqref{second-41} gives 
                \begin{align*}& \biggl[\prod_{i=1}^{5}(s_1+\al_{i} )\biggr] \what{5}{\al}{s}\nr&= (s_1+\sigma_3)(s_1+\sigma_4)(s_1+\sigma_5)\what{5}{\al}{s+e_1}
 + (s_1+\sigma_1)(s_1+\sigma_4)(s_1+\sigma_5)\what{5}{\al}{s+e_2}
  \nr&+ (s_1+\sigma_1)(s_1+\sigma_2)(s_1+\sigma_5)\what{5}{\al}{s+e_3}\nr&
 + (s_1+\sigma_1)(s_1+\sigma_2)(s_1+\sigma_3)\what{5}{\al}{s+e_4}
  \nr&+(s_1+\sigma_5) \what{4}{\al}{s+e_1+e_3}+(s_1+\sigma_3) \what{5}{\al}{s+e_1+e_4}+(s_1+\sigma_1) \what{5}{\al}{s+e_2+e_4}\nr&= (s_1+s_2-s_3)(s_1+s_3-s_4)(s_1+s_4)\what{5}{\al}{s+e_1}
  \nr&+(s_1+s_4) \what{5}{\al}{s+e_1+e_3}+(s_1+s_2-s_3) \what{5}{\al}{s+e_1+e_4}, \end{align*}the last step by the definition \eqref{sigdef} of the $\sigma_i$'s.  This is our desired shift equation in $s_1$.
  
  The shift equation in $s_2$ is derived analogously.  A fundamental difference in this derivation is that, in place of  \eqref{sym-decomp41}, we use the following expression involving  {\it Schur polynomials}   $s_\mu $   (see \cite[\S I.3]{MR553598}, especially Exercise 10 of that section):
  \begin{align*}  \prod_{1\le i<j\le5}(T+x_i+x_j)=  \sum_{\mu=(\mu_1,\mu_2,\ldots,\mu_5)\in S} \biggl(\frac{T}{2}\biggr)^{10-(\mu_1+\mu_2+\cdots+\mu_5)}d_{ \mu}\,s_\mu (x_1,x_2,\ldots, x_5).\end{align*}Here, $$S=\bigl\{ (\mu_1,\mu_2,\ldots, \mu_5)\in (\Z_{\ge0})^5: \mu_i\le 5-i\  (1\le i\le 5)\hbox{ and }  \mu_1\ge \mu_2\ge \cdots\ge \mu_5\bigr\},$$ and $d_\mu$ is the determinant of the matrix$$\biggl(\binom{2(5-i)}{\mu_j+5-j}\biggr)_{1\le i,j\le 5}.$$The Schur polynomials are symmetric polynomials in the $x_k$'s, and are therefore expressible in terms of the elementary symmetric polynomials in the $x_k$'s.  Techniques like those employed above, in the case $m=1$, therefore apply.  We omit the details.  
  \end{proof}  

  \begin{remark}\label{rmk:mellincont}
The above proof, in the case $m=1$ (that is, for the variable $s_1$---and therefore also for the variable $s_{n-1}$), generalizes to the case of $\GL(n,\R)$, for any $n\ge 2$.  For $2\le m\le n-2$, we do not yet have a proof that works for all $n\ge2$, though we expect that the above ideas and techniques should prove relevant.  Indeed, using the above methods, and applying  Mathematica to help with the more arduous calculations, we have been able to verify Conjecture \ref{conj} in full generality for $n\le7$.

 {We further note that, alternatively, one might   continue   $\widetilde{W}_{n,\al}(s)$ in the $s_j$'s by shifting or deforming the lines of integration in \eqref{eq:ISrecursion}.  Unfortunately such an approach has, thus far, failed to yield suitable results.  In particular, the residues that one obtains in moving these lines of integration past poles of the integrand are quite complicated, and do not seem to lend themselves to bounds of the type required to estimate $p_{T,R}^{(n)}(y)$.}
  \end{remark}
\subsection{\boldmath $ {p_{T,R}^{(n)}(y)}$ is a sum of a shifted term and residues}\label{sec:pTRnAsShiftPlusResidues}

Besides the shifted $p_{T,R}^{(n)}$ term (because we cross poles of $\what{n}{\al}{s}$ upon shifting the lines of integration) there are also many residue terms.  The residue terms will be parameterized by compositions of $n$.  Recall that a \emph{composition of length $r$} of a positive integer $n$ is a way of writing $n=n_1+\cdots+n_r$ as a sum of strictly positive integers.  Two sums that differ in the order define different compositions.  Compare this, on the other hand with \emph{partitions} which are compositions of $n$ for which the order doesn't matter.

\begin{definition}\label{def:compositions}
{\bf ($\pmb{a}$-admissible compositions)}  Let $a=(a_1,\ldots,a_{n-1})\in \R^{n-1}$. A composition $n=n_1+\cdots +n_r$ is termed \emph{$a$-admissible} if 
 \[  a_{\hn_i}>0 \mbox{ for all }i=1,\ldots,r-1. \]
The set of $a$-admissible compositions of length greater than one is 
 \[ \mathcal{C}_a := \left\{ \left.\begin{array}{c} \mbox{compositions} \\ n=n_1+\cdots+n_r \end{array} \ \right\vert \begin{array}{c} 2\leq r \leq n \\ a_{\hn_i}>0 \mbox{ for all }i=1,\ldots,r-1 \end{array}\right\}. \]
\end{definition}

\begin{remark}
At times we may also notate a composition $n=n_1+\cdots+n_r$ as an ordered list $C=(n_1,\ldots,n_r)$.
\end{remark}

\begin{definition}[\bf $\pmb{(r-1)}$-fold residue term]\label{def:k-fold-residue-term}
Suppose that $r\geq 2$ and $C\in \mathcal{C}_a$ is given by $n=n_1+\cdots+n_r$.  Let
 \[ \delta_C:=(\delta_1,\delta_2,\ldots,\delta_{r-1}) \in \big(\Z_{\geq 0}\big)^{r-1}\] with $0\leq \delta_i \leq \lfloor a_{\hn_i}\rfloor$ for each $i=1,\ldots,r-1$.  If $C$ has length two, we write $\delta_C=\delta$.  We define the \emph{$(r-1)$-fold residue term}
\begin{align}\label{eq:kresidue-term}
 p_{T,R}^{(n)}(y;-a,\delta_C) & := 
 \int\limits_{\substack{\halpha_n=0\\\re(\alpha)=0}}
 e^{\frac{\alpha_1^2+\cdots+\alpha_n^2}{T^2/2}} 
 \cdot 
 \FR{n}(\alpha)
 \bigg(\prod_{1\leq j\neq k \leq n} \hskip -4pt
 \Gamma_R\big(\alpha_j-\alpha_k\big) \bigg) \hskip -3pt
 \\ \nonumber & \hskip 12pt 
 \cdot
 \bigg( \prod_{i=1}^{r-1} y_i^{\frac{\hn_i(n-\hn_i)}{2}+2(\halpha_{\hn_i}+\delta_i)}\bigg)
 \cdot 
 \int\limits_{\substack{\re(s_j)=-a_j\\ j\notin \{\hn_1,\ldots,\hn_{r-1}\}}}
 \bigg( \prod_{j\notin \{\hn_1,\ldots,\hn_{r-1}\}} y_j^{\frac{j(n-j)}{2}-2s_j} \bigg)
 \\ \nonumber & \hskip -24pt \cdot
 \underset{s_{\hn_1}=-\halpha_{\hn_1}-\delta_1}{\res} 
    \left(\underset{s_{\hn_2}=-\halpha_{\hn_2}-\delta_2}{\res}\left( \cdots
    \left(\underset{s_{\hn_{r-1}}=-\halpha_{\hn_{r-1}}-\delta_{r-1}}{\res} \what{n}{\alpha}{s}\right) \cdots \right) \right)\, ds\, d\alpha.
\end{align}
\end{definition}

\begin{remark}\label{rmk:zero-unless-admissible}
In the shifted integral \eqref{eq:kresidue-term}, if $-a_i>0$ for some $i$, there will be no residues coming from the integral in $s_i$ because we are not shifting past any poles.  For this reason, one only obtains residue terms $p_{T,R}^{(n)}(y;-a,\delta_C)$ in the case that $C$ is $a$-admissible.  That said, equation~\eqref{eq:kresidue-term} makes perfect sense even if $C$ is not $a$-admissible.  In this case, $p_{T,R}^{(n)}(y;-a,\delta_C)$ is identically zero.
\end{remark}

\begin{proposition}\label{prop:pTR-expansion}
Suppose that $a=(a_1,\ldots,a_{n-1})\in \R^{n-1}$.  Then there exists constants $\kappa(C)$ such that
 \[ p_{T,R}^{(n)}(y) = p_{T,R}^{(n)}(y;-a)\ +\ \sum_{C\in \mathcal{C}_a} \kappa(C) \sum_{\substack{\delta_C=(\delta_1,\ldots,\delta_{r-1})\\ 0\leq \delta_i \leq \lfloor a_{\hn_i}\rfloor}} p_{T,R}^{(n)}(y;-a,\delta_C). \]
\end{proposition}

Before giving the proof, we make some preliminary remarks and observations.

\begin{remark}\label{rmk:only-one-kresidueterm}
Notice that an element $\sigma$ of the symmetric group $S_n$ (i.e., the group of permutations of a set of $n$ elements) acts on $\alpha=(\alpha_1,\ldots,\alpha_n)$ and, by extension, on $\halpha_k$ via
 \[ \sigma\cdot \halpha_k := \alpha_{\sigma(1)}+\alpha_{\sigma(2)}+\cdots+\alpha_{\sigma(k)}. \]
We can consider the analog to \eqref{eq:kresidue-term} obtained by replacing each instance of $\halpha_m$ with $\sigma\cdot \halpha_m$:
\begin{multline*}
 \int\limits_{\substack{\halpha_n=0\\\re(\alpha)=0}}
 e^{\frac{\alpha_1^2+\cdots+\alpha_n^2}{T^2/2}} 
 \cdot 
 \FR{n}(\alpha)
 \cdot
 \bigg(\prod_{1\leq j\neq k \leq n} \hskip -4pt
 \Gamma_R\big(\alpha_j-\alpha_k\big) \bigg) \hskip -3pt
 \bigg( \prod_{i\in \{\hn_1,\ldots,\hn_{r-1}\}} y_i^{\frac{i(n-i)}{2}+2(\sigma\cdot\halpha_i+\delta_i)}\bigg) 
 \\ \nonumber 
 \cdot 
 \int\limits_{\substack{\re(s_j)=-a_j\\ j\in \{\hn_1,\ldots,\hn_{r-1}\}}} \hskip -5pt
 \bigg( \prod_{j\notin \{\hn_1,\ldots,\hn_{r-1}\}} y_j^{\frac{j(n-j)}{2}-2s_j} \bigg)
 \underset{s_{i_1}=-\sigma\cdot\halpha_{i_1}-\delta_{i_1}}{\res} 
\left(\underset{s_{i_2}=-\sigma\cdot\halpha_{i_2}-\delta_{i_2}}{\res}
 \cdots 
    \underset{s_{i_k}=-\sigma\cdot\halpha_{i_k}-\delta_{i_k}}{\res}  \what{n}{\alpha}{s}    \right)
 \\ \nonumber 
 \cdot ds\, d\alpha
\end{multline*}
We make two observations:
\begin{itemize}
 \item As $C$ varies over all compositions of length $r$ and $\sigma$ varies over all possible permutations and $\delta_C$ varies over all $\big(\Z_{\geq 0}\big)^{r-1}$, one obtains all possible $(r-1)$-fold residues coming from shifting the lines of integration in $p_{T,R}^{(n)}(y)$.  This is a consequence of Theorem~\ref{th:residue-description} below.
 \item The action of $S_n$ on ordered subsets of $\{\alpha_1,\alpha_2,\ldots,\alpha_n\}$ given by permuting the indices is trivial on $\what{n}{\al}{s}$, i.e., $\what{n}{\sigma(\al)}{s}=\what{n}{\al}{s}$, and on the function
  \[ e^{\frac{\alpha_1^2+\cdots+\alpha_n^2}{T^2/2}} 
 \cdot 
 \FR{n}(\alpha)
 \cdot
 \bigg(\prod_{1\leq j\neq k \leq n} \hskip -4pt
 \Gamma_R\big(\alpha_j-\alpha_k\big) \bigg). \]
 This implies that relabeling the variables $\al_1,\al_2,\ldots \al_n$ by $\al_{\sigma^{-1}(1)},\al_{\sigma^{-1}(2)},\ldots,\al_{\sigma^{-1}(n)}$ everywhere (1) doesn't change the value of the integral, and (2) recovers the original integral given in \eqref{eq:kresidue-term}.
\end{itemize}
\end{remark}

\begin{remark}\label{rmk:cPsize}
The constant $\kappa(C)$ is the size of the (generic) orbit of the action of $S_n$ on the set
 \[ A=\{\halpha_{\hn_1},\ldots,\halpha_{\hn_{r-1}}\}. \]
Hence, defining the stabilizer of $A$ to be 
 \[ \stab(A):= \{\sigma \in S_n\mid \sigma\cdot \halpha_m=\halpha_m\mbox{ for each }m=\hn_1,\ldots,\hn_{r-1}\}, \] 
we see that
 \[ \kappa(C) = \frac{\# S_n}{\#\stab(A)} = \frac{n!}{\prod\limits_{i=1}^{r-1}(n_i!)}. \]
Since the exact value of $\kappa(C)$ is irrelevant to our application, we omit its proof below and leave it instead to the interested reader.
\end{remark}

\begin{proof}[Proof of Proposition~\ref{prop:pTR-expansion}]
Beginning with \eqref{eq:pTRnstart}, we see that $p_{T,R}^{(n)}(y)=p_{T,R}^{(n)}(y;b)$ for any $b=(b_1,\ldots,b_{n-1})$ with $b_i>0$ for each $i=1,\ldots,n-1$.  In order to compare this with $p_{T,R}^{(n)}(y;-a)$, we successively shift the lines of integration in the variables $s_k$ for each $k$ such that $-a_k<0$ (in descending order).  If $-a_k>0$ then shifting the line of integration from $\re(s)=b_k$ to $\re(s_k)=-a_k$ doesn't change the value of the integral in $s_k$.  In other words, there is a residue term if and only if the composition $C$ is admissible.

Beginning with the fact that
 \[ p_{T,R}^{(n)}(y) = p_{T,R}^{(n)}(y;b)\quad \mbox{for any $b=(b_1,\ldots,b_{n-1})$ for which $b_j>0$ for all $j$}, \]
we may shift the line of integration in $s_{n-1}$ to $\re(s_{n-1})=-a_{n-1}$.  In doing so, provided that $a_{n-1}>0$, we pass poles at $s_{n-1}=-\sigma\cdot \halpha_1-\delta_1$ for each $0\leq \delta \leq \lfloor a_1\rfloor$.  Hence, taking into account Remark~\ref{rmk:only-one-kresidueterm}, and considering $n=(n-1)+1$ (denoted $(n-1,1)$), it follows that
\begin{align}\label{eq:pTR-shift-in-s1}
 p_{T,R}^{(n)}(y) & = p_{T,R}^{(n)}(y;(b_1,b_2,b_3,\ldots,-a_{n-1})) 
 \\ \nonumber
 & \hskip 36pt + \kappa((n-1,1)) \cdot \sum_{\delta_{(n-1,1)}}p_{T,R}^{(n)}(y;(b_1,b_2,\ldots,-a_{n-1}),\delta_{(n-1,1)}),
\end{align}
where $\kappa((n-1,1))$ is a constant (which can be verified to agree with the description given in Remark~\ref{rmk:cPsize}.)

We now shift the line of integration in $s_{n-2}$ to $\re(s_{n-2})=-a_{n-2}$.  As before, provided that $a_{n-2}>0$, the Cauchy residue theorem and Remark~\ref{rmk:only-one-kresidueterm} give
\begin{align}\label{eq:pTR-shift-in-s1s2}
 p_{T,R}^{(n)}(y) & = p_{T,R}^{(n)}(y;(b_1,\ldots,b_{n-3},-a_{n-2},-a_{n-1})) 
 \\ \nonumber & \qquad + \kappa((n-2,2)) \hskip -6pt \sum_{\delta_{(n-2,2)}} p_{T,R}^{(n)}(y;(b_1,\ldots,b_{n-3},-a_{n-2},-a_{n-1}),\delta_{(n-2,2)})  \hskip -6pt
 \\ \nonumber & \qquad +
 \kappa((n-1,1)) \sum_{\delta_{(n-1)}} p_{T,R}^{(n)}(y;(b_1,\ldots,b_{n-3},-a_{n-2},-a_{n-1}),\delta_{(n-1,1)}) 
 \\ \nonumber & \qquad +
 \kappa((n-2,1,1)) \sum_{\delta_{(n-2,1,1)}} p_{T,R}^{(n)}(y;(b_1,\ldots,b_{n-3},-a_{n-2},-a_{n-1}),\delta_{(n-2,1,1)}),
\end{align}
for constants $\kappa(C)$ for each of $C=(n-1,1),(n-2,2),(n-2,1,1)$ as claimed.

We next repeat this process shifting the integrals in $s_{n-3}$ for each of the terms on the right of \eqref{eq:pTR-shift-in-s1s2}, and then again for $s_{n-4}$ and so forth (skipping those $s_m$ for which $a_m<0$) until all of the lines of integration have been moved to $\re(s_m)=-a_m$ for every possible integral.  The claimed formula is now evident.
\end{proof}

\subsection{Example: \boldmath$\GL(4)$}\label{sec:GL4-residue-terms}\label{sec:GL4-example}

We now consider the special case of $\what{4}{\al}{s}$ where
 \[ \alpha=(\alpha_1,\alpha_2,\alpha_3,\alpha_4)\in \big(i \R\big)^4, \quad \halpha_4=0.\]
Fix $\eps>0$.  Recall that $p_{T,R}^{(4)}(y) = p_{T,R}^{(4)}(y;(\eps,\eps,\eps))$.  If we now shift the lines of integration to $\re(s)=(-a)$ where $a=(a_1,a_2,a_3)\in \R^3$, then we get additional residue terms corresponding to each composition $4=n_1+\cdots+n_r$ and each $\delta_C\in \big(\Z_{\geq 0}\big)^r$ as follows.  

In general the composition $n=n_1+\cdots+n_r$ (by abuse of notation, we also think of this as a vector $(n_1,\ldots,n_r)$ so that $\widehat{n}_k = n_1+\cdots +n_k$) corresponds to taking an $(r-1)$-fold residue in the variables $s_{\widehat{n}_1},s_{\widehat{n}_2},\ldots,s_{\widehat{n}_{r-1}}$.   Here is a table of the residues corresponding to the different compositions. 

\[
\begin{array}{|c|c|c|}
\hline
\mbox{composition $C$} & \mbox{residues in $s$-variables} & \delta_C
\\ \hline
1+3 & s_1=-\alpha_1-\delta_1 & (\delta_1) \\
2+2 & s_2=-\alpha_1-\alpha_2-\delta_2 & (\delta_2) \\
3+1 & s_3=-\alpha_1-\alpha_2-\alpha_3-\delta_3 & (\delta_3)
\\ \hline
1+1+2 & s_1=-\alpha_1-\delta_1,\  s_2=-\alpha_1-\alpha_2-\delta_2 & (\delta_1,\delta_2)
\\ 
1+2+1 & s_1=-\alpha_1-\delta_1,\  s_3=-\alpha_1-\alpha_2-\alpha_3-\delta_3 & (\delta_1,\delta_3)
\\ 
2+1+1 & s_2=-\alpha_1-\alpha_2-\delta_2, \ s_3=-\alpha_1-\alpha_2-\alpha_3-\delta_3 & (\delta_2,\delta_3)
\\ \hline
\end{array}
\]
In each case $0\leq \delta_i\leq \lfloor a_i \rfloor$.  Not included in the table are the triple residues in $s_i=-\halpha_i-\delta_i$ for each $i=1,2,3$.  These correspond to the composition $4=1+1+1+1$ and $\delta_C=(\delta_1,\delta_2,\delta_3)$.

\subsection{The integral \boldmath $ {\mathcal{I}_{T,R}^{(m)}(-a)}$ in terms of an explicit recursive formula for $\what{m}{\alpha}{s}$}\label{sec:ITRMa}

At first glance, the following definition appears to be relevant only for the shifted $p_{T,R}^{(n)}$-term, as it is essentially equal to $p_{T,R}^{(n)}((1,\ldots,1);-a)$, and not for the shifted residue terms.  However, it will turn out to be pivotal to bounding the residue terms as well.

\begin{definition}[\bf{The integral \boldmath $\mathcal{I}_{T,R}^{(m)}$}]
Let $m\geq 2$ be an integer and $a=(a_1,\ldots,a_{m-1})\in \R^{m-1}$.  Then we define
\begin{equation}\label{eq:ITRm}
 \mathcal{I}_{T,R}^{(m)}(-a) := \int\limits_{\substack{\halpha_m=0\\ \re(\alpha)=0}}
 e^{\frac{\alpha_1^2+\cdots+\alpha_m^2}{T^2/2}} 
 \cdot 
 \FR{n}(\alpha)
 \bigg(\prod_{1\leq j\neq k \leq m} \hskip -4pt
 \Gamma_R\big(\alpha_j-\alpha_k\big) \bigg)
 \int\limits_{\substack{s=(s_1,\ldots,s_{m-1})\\\re(s)=-a}} \hskip -8pt
\left| \widetilde{W}_{m,\alpha}(s)\right|\, ds\, d\alpha.
\end{equation}
\end{definition}

As alluded to above, inserting the result of Theorem~\ref{th:IS} into \eqref{eq:shifted-pTR}, we find that
\begin{align*}
 \big| p_{T,R}^{(n)}(y,-a) \big| & \ll
 \bigg(\prod_{j=1}^{n-1}
 y_j^{\frac{j(n-j)}{2}-2a_j} \bigg) \mathcal{I}_{T,R}^{(n)}(-a).
\end{align*}
Hence, giving a bound for $p_{T,R}^{(n)}(y)$ requires only that we bound $\mathcal{I}_{T,R}^{(m)}(-a)$ in the case of $m=n$.  However, much more is true: we will show that if $C$ is the composition $n=n_1+\cdots +n_r$, then $p_{T,R}^{(n)}(y;-a,\delta_C)$ can be bounded by the same product of $y_i$'s as above times a certain power of $T$ and a product of the form
 \[ \prod_{\ell=1}^{r-1} \mathcal{I}_{T,R}^{(n_\ell)}(-a^{(\ell)}), \]
for certain values $a^{(\ell)}=(a_1^{(\ell)},\ldots,a_{n_\ell-1}^{(\ell)})$ which depend on the value of $a=(a_1,\ldots,a_{n-1})\in \R^{n-1}$.

The significance of this fact should not be understated.  Without it, we would be required to treat nearly every possible composition $C$ (hence each possible residue term) individually.  Indeed, returning to the case of $n=4$, as noted in Section~\ref{sec:GL4-example} above, there were seven residue terms.  The only symmetries that we were able to exploit in \cite{GSW21} to help were that the $(1,3)$ and $(3,1)$ residues were equivalent, and the $(1,1,2)$ and $(2,1,1)$ residues were equivalent as well.  This left five individual distinct cases, each of which required several pages of work to bound.  So, although the method of this paper does require dealing with some tricky notation and combinatorics, it eliminates the need to treat each residue on its own terms.

\section{\large\bf Bounding \boldmath $\mathcal{I}_{T,R}^{(m)}$}\label{sec:ITRmBound}

Recall that for $\alpha=(\alpha_1,\ldots,\alpha_m) \in \C^m$ satisfying $\halpha_m=0$ and $a=(a_1,a_2,\ldots,a_{n-1})\in \R^{n-1}$,
\begin{equation}\label{eq:ITRn-minusa}
 \mathcal{I}_{T,R}^{(m)}(-a) := \int\limits_{\substack{\halpha_m=0\\ \re(\alpha)=0}} 
 e^{\frac{\alpha_1^2+\cdots+\alpha_m^2}{T^2/2}} 
 \cdot 
 \FR{m}(\alpha) \hskip -3pt
 \prod_{1\leq j\neq k \leq m} \hskip -4pt
 \big| \Gamma_R(\alpha_j-\alpha_k)\big|
 \int\limits_{\re(s)=-a} \hskip -5pt
 \left| \widetilde{W}_{m,\alpha}(s)\right|\, ds\, d\alpha.
\end{equation}

\begin{theorem}\label{th:ITRnaBound}
Let $\mathcal{I}_{T,R}^{(m)}(-a)$ be as above and set $D(m)=\deg(\mathcal{F}_1^{(m)}(\al))$.   Then for any $0<\varepsilon<\frac12$, 
 \[ \mathcal{I}_{T,R}^{(m)}(-a) \ll T^{\varepsilon+\frac{(m+4)(m-1)}{4} + R\cdot\big(D(m)+\frac{m(m-1)}{2}\big) - \sum\limits_{j=1}^{m-1}B(a_j)}, \]
where
 \[ B(c) = \begin{cases} 0 & \mbox{ if }c<0 \\ \lfloor c \rfloor+2(c - \lfloor c \rfloor) & \mbox{ if } 0<
\lfloor c \rfloor+\varepsilon < c \leq \lfloor c \rfloor+\frac12, \\ \lceil c \rceil & \mbox{ if } \frac12 < \lceil c \rceil - \frac12\leq c <\lceil c \rceil-\varepsilon. \end{cases} \]
The implicit constant depends on $\varepsilon$, $R$ and $m$.
\end{theorem}
Theorem~\ref{th:IS} allows us to write $\what{m}{\alpha}{s}$ in terms of an integral of the product of several  {gamma} functions and the lower rank Mellin transform $\what{m-1}{\beta}{z}$ where
 \[ \beta = \big(\beta_1,\ldots,\beta_{m-1}\big) := \left( \alpha_1+\frac{\alpha_m}{m-1},\ldots,\alpha_{m-1} + \frac{\alpha_m}{m-1}\right). \]
Using this, we are able siphon off the contribution to the integrand of \ref{eq:ITRn-minusa} that is independent of the variable $\beta$.  This in turn allows us to relate $\mathcal{I}_{T,R}^{(m)}$ to $\mathcal{I}_{T,R}^{(m-1)}$ and prove the result inductively.

\subsection{Symmetry of integration in $\alpha$}

Since the integrand of \eqref{eq:ITRn-minusa} is invariant under the action of $\sigma \in S_m$ acting on $\alpha=(\alpha_1,\ldots,\alpha_m)$, we may restrict the integration to a fundamental domain.  A choice of such a fundamental domain is
\begin{equation}\label{eq:alphaorder}
 \im(\al_1) \geq \im(\al_2) \geq \cdots \geq \im(\al_m).
\end{equation}
Hence, \eqref{eq:ITRn-minusa} is equal, up to a constant, to the same integral but restricted to $\alpha$ satisfying \eqref{eq:alphaorder}.  In the sequel we will always assume that \eqref{eq:alphaorder} holds.

\subsection{Extended exponential zero set}

Recall that Stirling's asymptotic formula (for $\sigma\in \R$ fixed and $t\in \R$ with $\lvert t \rvert\to \infty$) is given by
\begin{equation} \label{eq:Stirling}
 \Gamma(\sigma+it)\sim \sqrt{2\pi}\cdot \lvert t \rvert^{\sigma-\frac12}\, e^{-\frac{\pi}{2}\lvert t\rvert }.
\end{equation}

\begin{definition}[{\bf Exponential and Polynomial Factors of a Ratio of Gamma Functions}] \label{def:poly-exp-part}
We call $\lvert t\rvert^{\sigma-\frac12}$ the \emph{polynomial factor} of $\Gamma(\sigma+it)$, and $e^{-\frac{\pi}{2}\lvert t\rvert}$ is called the \emph{exponential factor}.  For a ratio of  {gamma} functions, the \emph{polynomial (respectively, exponential) factor} is  {composed of} the polynomial (respectively, exponential) factors of each individual Gamma function.  
\end{definition}

\begin{lemma}[\bf{Extended Exponential Zero Set}]
\label{lem:ExpZeroSet}
Assume that $\alpha\in \C^m$ is a Langlands parameter satisfying
 \[ \im(\alpha_1)\geq \im(\al_2)\geq \ldots \geq \im(\alpha_m). \]
Then the integrand of $\mathcal{I}_{T,R}^{(m)}$ (as a function of $s$) has exponential decay outside of the set $I = I_1\times I_2 \times\cdots \times I_{m-1}$, where 
 \[ I_j := \left\{s_j \left \lvert -\sum_{k=1}^{j} \im(\alpha_k) \leq \im(s_j) \leq -\sum_{k=1}^{j} \im(\alpha_{m-k+1}) \right.\right\}. \]
\end{lemma}
\begin{remark}
See \cite{GSW21} for the definition of the \emph{exponential zero set} of an integral.  The extended exponential zero set given in Lemma~\ref{lem:ExpZeroSet} contains the exponential zero set for $\mathcal{I}_{T,R}^{(m)}$.
\end{remark}
\begin{proof}
We first prove Lemma~\ref{lem:ExpZeroSet} in the case that $m=2$.  In the formula \eqref{eq:ITRn-minusa} for $\mathcal{I}_{T,R}^{(n)}$, replace $\what{2}{\alpha}{s_1}$ with $\Gamma(s_1+\alpha_1)\Gamma(s_1+\alpha_2)$.  Then assuming \eqref{eq:alphaorder}, the exponential factor is $e^{\frac{\pi}{2}\mathcal{E}(s,\alpha)}$ where
 \[ \mathcal{E}(s,\alpha) = \big\lvert \im(s_1)+\im(\al_1)\big\rvert+\big\lvert\im(s_1)+\im(\al_2)\big\rvert -2\im(\alpha_1).\]
We see, therefore, that the exponential factor $\mathcal{E}(s,\alpha)$ is negative unless
 \[ \im(s_1)+\im(\al_1) \geq 0 \quad \mbox{and} \quad \im(s_1)+\im(\al_2) \leq 0 \quad \Longleftrightarrow \quad -\im(\alpha_1) \leq \im(s_1) \leq -\im(\al_2), \]
as claimed.

Let us suppose that $m\geq 3$ and  $c=(c_1,c_2,\ldots,c_{m-1})$ with $c_j>0$ ($j=1,2,\ldots,m-1$).  In order to prove Lemma~\ref{lem:ExpZeroSet} using induction on $m$, we make use of the change of variables
 \[ \beta_j = \alpha_j+\frac{\alpha_m}{m-1}, \quad j=1,\ldots,m-1. \]
Observe that $\beta_1\cdots+\beta_{m-1}=0$.  By Lemma~\ref{lem:alphasumk} in the case that $k=m-1$, 
\[ \al_1^2+\cdots+\alpha_m^2 = \beta_1^2+\cdots \beta_{m-1}^2+ \frac{m}{m-1}\alpha_m^2. \]
Then in the integrand for $\mathcal{I}_{T,R}^{(m)}(c)$ we may substitute the formula for $\what{m}{\alpha}{s}$ given in Theorem~\ref{th:IS}.  We also use the fact (see Lemma~\ref{lem:prodGammaRkdecomp}) that
 \[ \prod_{1\leq j\neq k\leq m} \Gamma(\alpha_j-\alpha_k) = \left( \prod_{1\leq j\neq k\leq m-1} \Gamma(\beta_j-\beta_k) \right) \cdot \left( \prod_{i=1}^{m-1} \Gamma(\alpha_m-\alpha_i) \Gamma(\alpha_i-\alpha_m) \right), \]
and, via Stirling,
 \[ \prod_{i=1}^{m-1} \Gamma(\alpha_m-\alpha_i) \Gamma(\alpha_i-\alpha_m) \ll e^{\pi \im(\alpha_m)}. \]
Note that \eqref{eq:alphaorder} implies that $\im(\alpha_m)\leq 0$, hence,
\begin{multline*}
 \mathcal{I}_{T,R}^{(m)}(c) \ll \int\limits_{\re(\alpha_m)=0} e^{\frac{m}{m-1}\frac{\alpha_m^2}{T^2/2}}
 \int\limits_{\substack{\hbeta_{m-1}=0 \\ \re(\beta)=0}} 
 e^{\frac{\beta_1^2+\cdots+\beta_{m-1}^2}{T^2/2}} 
 \cdot 
 \left|\mathcal{P}_{(D(m)-D(m-1))R}(\alpha_m,\beta)\right| 
 \\ 
 \cdot
 \FR{m-1}(\beta) \hskip -3pt
 \prod_{1\leq j\neq k \leq m-1} \hskip -4pt
 \big| \Gamma_R(\beta_j-\beta_k)\big|  \int\limits_{\substack{\re(z_j)=b_j \\ 1\leq j \leq m-2}} 
 \left| \widetilde{W}_{m-1,\beta}\left(z \right)\right|
 \\ 
 \cdot \prod_{j=1}^{m-1} 
 \int\limits_{\re(s_j)=c_j} \hskip -6pt
 \left|\Gamma\Big(s_j-z_{j-1}+\tfrac{(m-j)\alpha_m}{m-1}\Big)
 \Gamma\Big(s_j-z_j-\tfrac{j\alpha_m}{m-1}\Big)
 \Gamma_R\Big(\tfrac{-m}{m-1}\alpha_m-\beta_j\Big)
 \Gamma_R\Big(\beta_j+\tfrac{m}{m-1}\alpha_m\Big) \right|
 \\ 
 \; ds_j
 \,
 dz\, d\alpha.
\end{multline*}
By the induction hypothesis, the second row of this expression has exponential decay outside of the set
\begin{equation}\label{eq:zj-bounds}
 \left\{ z=(z_1,\ldots,z_{m-2})\, \left|\, -\sum_{k=1}^j \beta_k \leq \im(z_j) \leq  -\sum_{j=1}^k \beta_{m-j} \right.\right\},
\end{equation}
for each $k=1,2,\ldots,m-2$.  (Recall that $z_0=z_{m-1}=0$.)

The assumption $\im(\alpha_j)\geq \im(\alpha_m)$ and the definition of $\beta_j$ above imply that
 \[ \im(\alpha_j-\alpha_m) = \im\Big(\beta_j+\tfrac{m}{m-1}\alpha_n\Big) \geq 0 \quad \mbox{($j=1,2,\ldots,m-1$).} \]
Thus, the exponential factor coming from the final line in the expression above is $e^{\frac{\pi}{2}\mathcal{E}(s,z,\beta,\alpha_m)}$ where
\begin{align*}
 \mathcal{E}(s,z,\beta,\alpha_n) & =
 \sum_{j=1}^{n-1} \Big( \big| \im(s_j-z_{j-1}+\tfrac{n-j}{n-1}\alpha_n)\big| + \big| \im(s_j-z_j-\tfrac{j}{n-1}\alpha_n)\big| - \im\big(\tfrac{n}{n-1}\alpha_n+\beta_j\big) \Big)
 \\ & =
 \sum_{j=1}^{n-1} \Big( \big| \im(s_j-z_{j-1}+\tfrac{n-j}{n-1}\alpha_n)\big| + \big| \im(s_j-z_j-\tfrac{j}{n-1}\alpha_n)\big|\Big) - n\im(\alpha_n).
\end{align*}
 {We know that the integral defining $ \mathcal{I}_{T,R}^{(m)}(-a)$ is convergent.  Therefore,} it must be the case that $\mathcal{E}(s,z,\beta,\alpha_m)\leq 0$.  In order to find where $\mathcal{E}=0$, i.e., where there is \emph{not} exponential decay, we seek for values $\epsilon_{1,1},\epsilon_{2,1},\ldots,\epsilon_{1,m-1},\epsilon_{2,m-1}\in \{\pm 1\}$ for which 
\begin{equation}\label{eq:ExpZeroSetFormula}
 \sum_{j=1}^{m-1} \left( \epsilon_{1,j} \im(s_j-z_{j-1}+\tfrac{m-j}{m-1}\alpha_m) + \epsilon_{2,j}\im(s_j-z_j-\tfrac{j}{m-1}\alpha_m)\right) = m\im(\alpha_m).
\end{equation}
In order for the $s$-variables to cancel it is clear that for each $j=1,2,\ldots,m-1$ it need be true that $\epsilon_j:=\epsilon_{1,j}=-\epsilon_{2,j}$.  With this assumption, equation~\ref{eq:ExpZeroSetFormula} simplifies:
\begin{align*}
 \sum_{j=1}^{m-1} \Big(\epsilon_j \im\big(z_j-z_{j-1}+\tfrac{m}{m-1}\alpha_m\big) \Big) = m\im(\alpha_m).
\end{align*}
In order for this to hold true, it is necessary that $\epsilon_j=1$ for all $j$, since otherwise, the coefficients of $\alpha_m$ on each side of the inequality wouldn't match.  On the other hand, $\epsilon_j=1$ for all $j$ is sufficient as well since
 \[ \sum_{j=1}^{m-1} \im(z_{j-1}-z_j) = \im(z_0-z_{m-1}) =0. \]
This unique solution to \eqref{eq:ExpZeroSetFormula} implies, therefore, that there is exponential decay in the integrand of $\mathcal{I}_{T,R}^{(m)}$ above unless
 $\im(z_{j-1} - \tfrac{m-j}{m-1} \alpha_m) \leq \im(s_j) \leq \im(z_j + \tfrac{j}{m-1} \alpha_m).$
The inductive assumption \eqref{eq:zj-bounds} implies that
 \[ \im(z_{j-1}-\tfrac{m-j}{m-1}\alpha_m) 
  \geq 
  -\sum_{k=1}^{j-1}\Big( \beta_k-\tfrac{\alpha_m}{m-1}\Big) - \alpha_m = -\sum_{k=1}^j \alpha_j,
 \]
and
 \[ \im(z_j + \tfrac{j}{m-1} \alpha_m) 
  \leq
  -\sum_{k=1}^j \big(\beta_k-\tfrac{\alpha_m}{m-1}\big)
  = -\sum_{k=1}^m\alpha_k,
 \]
thus yielding the desired bounds on $\im(s_j)$.

To complete the proof, we remark that if $-a<0$, in order to use the result of Theorem~\ref{th:IS}, we need to first apply the shift equations given in Corollary~\ref{cor:shiftequations} below.  This will allow us to rewrite $\mathcal{I}_{T,R}^{(m)}(-a)$ as a sum over terms all of which have the same basic form as that for $\mathcal{I}_{T,R}^{(m)}(c)$ with $c>0$.  Each of these terms has precisely the same exponential factor since this depends only on the imaginary parts of the arguments of the Gamma functions, hence the same exponential zero set is determined in general.
\end{proof}

For each $j=1,\ldots,n$, we define
\begin{equation} \label{eq:Bj}
 \mathcal{B}_j(s_j,\alpha):= \prod_{\substack{K\subseteq \{1,\ldots,n\}\\ \#K=j}} \Big( s_j + \sum_{k\in K}\alpha_k \Big).
\end{equation}
Using this, the following corollary is easily deduced.  (See \cite{GSW21} for the case of $n=4$.)
\begin{corollary}\label{cor:shiftequations}
Let $r=(r_1,\ldots,r_{n-1})\in \Z_{\geq 0}^{n-1}$.  There exists a sequence of shifts $\sigma=(\sigma_1,\ldots,\sigma_{n-1})\in \Z_{\geq 0}^{n-1}$ and polynomials $Q_{\sigma,r}(s,\alpha)$ such that
 \[ \left| \widetilde{W}_{n,\alpha}(s) \right| \ll \sum_\sigma \frac{\lvert Q_{\sigma,r}(s,\alpha)\rvert }{\prod\limits_{j=1}^{n-1} \lvert \mathcal{B}_j(s_j,\alpha)\rvert^{r_j}} \left| \widetilde{W}_{n,\alpha}(s+r+\sigma) \right|, \]
where
 \[ Q_{\sigma,r}(s,\alpha) = \prod_{j=1}^{n-1} P_{\sigma_j,r_j}(s,\alpha),  \qquad
 \deg(P_{\sigma_j,r_j}(s,\alpha)) =  r_j\left( \binom{n}{j}-2\right)-2\sigma_j. \]
\end{corollary}

\subsection{Proof of Theorem~\ref{th:ITRnaBound} in the case $m=2$}
\begin{proof}
As in the proof of Lemma~\ref{lem:ExpZeroSet}, we can replace $\what{2}{(\alpha,-\alpha)}{s}$ with $\Gamma(s+\alpha)\Gamma(s-\alpha)$ and estimate using Stirling's bound.  We may, moreover, restrict $s$ to the exponential zero set $-\im(\alpha)\leq \im(s) \leq \im(\al)$ to see that
\begin{align*}
\mathcal{I}_{T,R}^{(2)}(-a) & = \int\limits_{\substack{\halpha_n=0\\\re(\alpha)=0}} 
 e^{\frac{\alpha^2}{T^2}} 
 \cdot 
 \hskip -3pt
 \big| \Gamma_R(2\alpha) \Gamma_R(-2\alpha)  \big|
 \int\limits_{\substack{s=(s_1,\ldots,s_{n-1})\\\re(s)=-a}} \hskip -5pt
 \left| \widetilde{W}_{2,(\alpha,-\alpha)}(s)\right|\, ds\, d\alpha
 \\ & \ll
 \int\limits_{\substack{\halpha_n=0\\\re(\alpha)=0}} 
 e^{\frac{\alpha^2}{T^2}} 
 \cdot 
 \big(1+\left| 2\im(\al) \right|\big)^{R+\frac12}
  \hskip -3pt
  {\int\limits_{\substack{\re(s)=-a\\ -\im(\al)\leq \im(s)\leq \im(\alpha)}}} \hskip -5pt
  \big(1+\lvert \im(s)-\im(\alpha)\rvert\big)^{-a-\frac12}
  \\ & \hskip 120pt \cdot
  \big(1+\lvert \im(s)+\im(\alpha)\rvert\big)^{-a-\frac12}
 \, ds\, d\alpha.
\end{align*}
Due to the presence of the term $e^{\frac{\alpha^2}{T^2}}$, we may assume moreover that $\im(\alpha) \leq T^{1+\varepsilon}$.  Thus, we have the bound
\begin{align*}
 \mathcal{I}_{T,R}^{(2)}(-a) & \ll \ \int\limits_{\substack{\re(\alpha)=0\\0\leq \im(\alpha) \leq T^{\varepsilon+1}}} \big(1+2\lvert \alpha\rvert\big)^{R+\frac12}
 \hskip -12pt
 \int\limits_{\substack{\re(s)=-a \\ -\im(\alpha)\leq \im(s) \leq \im(\alpha)}}
 \hskip -12pt
 \big(1+\alpha-s\big)^{-a-\frac12}\big(1+\alpha-s\big)^{-a-\frac12}\, ds\; d\al
 \\
 & \ll
 \int\limits_{\substack{\re(\alpha)=0\\0\leq \im(\alpha) \leq T^{\varepsilon+1}}}  \big(1+2\lvert \alpha\rvert\big)^{R+\frac12-\min\big\{a+\frac12,2a\big\}}\; d\al \ll T^{\varepsilon+R+\frac32-\min\big\{a+\frac12,2a\big\}},
\end{align*}
In the statement of Theorem~\ref{th:ITRnaBound}, the claimed bound is $\mathcal{I}_{T,R}^{(2)}(-a) \ll T^{\eps+R+\frac32- B(a)}$, where  {$B(a)$ is as defined in Theorem \ref{th:ITRnaBound}}.
We have in fact proved that $\mathcal{I}_{T,R}^{(2)}(-a) \ll T^{\eps+R+\frac32- B'(a)}$, where
 \[  B'(a) = \max\big\{a+\tfrac12,2a\big\}=\begin{cases} 2a & \mbox{ if } \eps < a \leq \frac12, \\ a +\frac12 & \mbox{ if } a\geq \frac12. \end{cases} \]
If, $a<0$, then we may shift the integral over $\re(s)=-a$ to be as close to $\re(s)=0$ as desired; indeed, we may make the shift to the point that the error can be absorbed into the $\eps$ term in the power of $T$.  Therefore, since $B(a)\leq B'(a)$ for all $a>0$, the Theorem follows.
\end{proof}

\pagebreak

\subsection{Proof of Theorem~\ref{th:ITRnaBound} for general $m$}
\begin{proof}
Let $m\geq3$ and assume that Theorem~\ref{th:ITRnaBound} has been shown to be true for all integers $2\leq k <m$.  It follows from Corollary~\ref{cor:shiftequations} with $r_j= \lceil a_j \rceil$ that
\begin{align*}
 \mathcal{I}_{T,R}^{(m)}(-a) & \ll \sum_\sigma \int\limits_{\substack{\halpha_m=0\\ \re(\alpha)=0}}
 e^{\frac{\alpha_1^2+\cdots+\alpha_m^2}{T^2/2}} 
 \cdot 
 \FR{m}(\alpha) \hskip -3pt
 \prod_{1\leq j\neq k \leq m} \hskip -4pt
 \big| \Gamma_R(\alpha_j-\alpha_k)\big|
 \\ & \hskip 48pt
 \int\limits_{\substack{s=(s_1,\ldots,s_{m-1})\\\re(s)=-a}} \hskip -5pt
 \frac{\lvert \mathcal{P}_{d(m)-2\lvert\sigma\rvert}(s,\alpha)\rvert }{\prod\limits_{j=1}^{m-1} \lvert \mathcal{B}_j(s_j,\alpha)\rvert^{\lceil a_j \rceil}} \left| \widetilde{W}_{m,\alpha}(s+r+\sigma)\right|
 \, ds\, d\alpha 
\end{align*}
By Theorem~\ref{th:IS},
\begin{align*}
 \mathcal{I}_{T,R}^{(m)}(-a) & \ll
 \sum_\sigma \int\limits_{\substack{\halpha_m=0\\ \re(\alpha)=0}} \hskip -4pt
 e^{\frac{\alpha_1^2+\cdots+\alpha_m^2}{T^2/2}} 
 \cdot 
 \FR{m}(\alpha) \hskip -3pt
 \prod_{1\leq j\neq k \leq m} \hskip -4pt
 \big| \Gamma_R(\alpha_j-\alpha_k)\big|
 \int\limits_{\re(s)=-a} \hskip -5pt
 \frac{\lvert \mathcal{P}_{d(m)-2\lvert\sigma\rvert}(s,\alpha)\rvert }{\prod\limits_{j=1}^{m-1} \lvert \mathcal{B}_j(s_j,\alpha)\rvert^{\lceil a_j \rceil}}
 \\ & \hskip 48pt \int\limits_{\substack{z=(z_1,\ldots,z_{m-2})\\ \re(z)=b}} \Bigg( \prod_{j=1}^{m-1} \left|\Gamma\Big(s_j+\lceil a_j \rceil+\sigma_j-z_{j-1}+\frac{(m-j)\alpha_m}{m-1}\Big)\right|
 \\ & \hskip 96pt
 \left|\Gamma\Big(s_j+\lceil a_j \rceil+\sigma_j-z_j-\frac{j\alpha_m}{m-1}\Big)\Bigg)\right|
 \cdot
 \left| \widetilde{W}_{m-1,\beta}\left(z \right)\right|\, dz\, ds\, d\alpha. 
\end{align*}
Next, we use the functional equation for the  {gamma} function to rewrite
\begin{multline*}
\Gamma\Big(s_j+\lceil a_j \rceil+\sigma_j-z_{j-1}+\frac{(m-j)\alpha_m}{m-1}\Big)\Gamma\Big(s_j+\lceil a_j \rceil+\sigma_j-z_j-\frac{j\alpha_m}{m-1}\Big) 
\\
= \mathcal{P}_{2\sigma_j}(s,z,\alpha)\Gamma\Big(s_j+\lceil a_j \rceil-z_{j-1}+\frac{(m-j)\alpha_m}{m-1}\Big)\Gamma\Big(s_j+\lceil a_j \rceil-z_j-\frac{j\alpha_m}{m-1}\Big).
\end{multline*}
Additionally, we use the fact that the integrand has exponential decay unless $\lvert \alpha_1\rvert, \ldots,\lvert \alpha_m\rvert \leq T^{1+\varepsilon}$, and by Lemma~\ref{lem:ExpZeroSet}, each of the variables $s_j$ are bounded in terms of $\alpha$.  This means that we may replace the polynomials $\mathcal{P}_{2\sigma_j}$ with the bound $T^{\varepsilon+2\sigma_j}$.  Note that in doing so, the dependence on $\sigma$ is removed:
\begin{align*}
 \mathcal{I}_{T,R}^{(m)}(-a) & \ll T^{\varepsilon+\sum\limits_{j=1}^{m-1} \lceil a_j\rceil \left( \binom{m}{j}-2 \right)} 
 \int\limits_{\substack{\halpha_m=0\\ \re(\alpha)=0}}
 e^{\frac{\alpha_1^2+\cdots+\alpha_m^2}{T^2/2}} 
 \cdot 
 \FR{m}(\alpha) \hskip -3pt
 \prod_{1\leq j\neq k \leq m} \hskip -4pt
 \big| \Gamma_R(\alpha_j-\alpha_k)\big|
 \int\limits_{\substack{s=(s_1,\ldots,s_{m-1})\\\re(s)=-a}}
 \\ & \hskip 24pt 
 \cdot\int\limits_{\substack{z=(z_1,\ldots,z_{m-2})\\\re(z)=b}} \left( \prod_{j=1}^{m-1} \frac{\left|\Gamma\Big(s_j+\lceil a_j\rceil-z_{j-1}+\frac{(m-j)\alpha_m}{m-1}\Big)
 \Gamma\Big(s_j+\lceil a_j\rceil-z_j-\frac{j\alpha_m}{m-1}\Big)\right|}{\lvert \mathcal{B}_j(s_j,\alpha)\rvert^{\lceil a_j\rceil}}\right)
 \\ & \hskip 192pt
 \cdot\left| \widetilde{W}_{m-1,\beta}\left(z \right)\right|\, dz\, ds\, d\alpha.
\end{align*}
Notice that the conclusion of Proposition~\ref{prop:ITR-minusa-plusL} follows from the last several steps by simply replacing $s$ by $s+L$ in the integrand (or, equivalently, replace $\re(s)=-a$ by $\re(s)=-a+L$ in the domain of integration), and then at the step where the functional equation of Gamma is used to remove $\sigma$ from the  {gamma} functions, we remove $L$ in the exact same fashion.

We deduce that
\begin{align*}
 \mathcal{I}_{T,R}^{(m)}(-a) & \ll T^{\varepsilon+\sum\limits_{j=1}^{m-1} \lceil a_j \rceil \left( \binom{m}{j}-2 \right)} \cdot 
 \int\limits_{\substack{\halpha_m=0\\ \re(\alpha)=0}}
 e^{\frac{\al_1^2+\cdots+\al_m^2}{T^2/2}} 
 \cdot
 \FR{m}(\al) \hskip -6pt
 \prod_{1\leq j\neq k \leq m-1} \hskip -4pt
 \big| \Gamma_R(\al_j-\al_k)\big|  \int\limits_{\substack{z=(z_1,\ldots,z_{m-2})\\\re(z)=b}}
 \\ & \hskip -30pt
 \prod_{j=1}^{m-1} 
 \int\limits_{\re(s_j)=\lceil a_j\rceil-a_j} \hskip -10pt
 \frac{\left|\Gamma\Big(s_j-z_{j-1}-\frac{(m-j)\halpha}{m-1}\Big)
 \Gamma\Big(s_j-z_j+\frac{j\halpha}{m-1}\Big)\right|}{\lvert \mathcal{B}_j(s_j,\alpha)\rvert^{\lceil a_j\rceil}} \Big| \Gamma_R\big(\tfrac{n}{n-1}\halpha-\beta_j\big)
 \Gamma_R\big(\beta_j-\tfrac{m}{m-1}\halpha\big) \Big|
 \\ & \hskip 192pt 
 \cdot
 \left| \widetilde{W}_{m-1,\beta}\left(z \right)\right|
 ds_j\, dz\, d\alpha.
\end{align*}
Note that we have also made the change of variable $s\mapsto s_j-\lceil a_j\rceil$ for each $j=1,2,\ldots,m-1$, and we are using the notation $\halpha:=-\alpha_m$.  (Using the terminology of Lemma~\ref{lem:alphasumk} in the case of $k=m-1$, we have $\halpha=\halpha_{m-1}$.)  As in the case of $n=2$, due to the presence of the exponential terms, we see that the integral has exponential decay unless $\lvert \al_j\rvert \ll T^{1+\varepsilon}$.

\begin{lemma}\label{lem:IjBound}
Let $\alpha=(\alpha_1,\ldots,\alpha_m)$ and $\beta_j=\alpha_j-\frac{\halpha}{m-1}$ be as above.  In particular, they are purely imaginary with $\lvert \beta_k\rvert, \lvert \halpha_j\rvert < T^{1+\eps}$.  Suppose, moreover, that $\alpha$ is in $j$-general position.  Then
\begin{multline*}
 \int\limits_{\re(s_j)=\lceil a_j \rceil-a_j} \hskip -8pt
 \frac{\left|\Gamma\Big(s_j-z_{j-1}-\frac{(m-j)\halpha}{m-1}\Big)
 \Gamma\Big(s_j-z_j+\frac{j\halpha}{m-1}\Big)\right|}{\lvert \mathcal{B}_j(s_j,\alpha)\rvert^{\lceil a_j \rceil}} \Big| \Gamma_R\big(\tfrac{m}{m-1}\halpha-\beta_j\big)
 \Gamma_R\big(\beta_j-\frac{m}{m-1}\halpha\big) \Big| \; ds_j
\\
  \ll T^{\eps+R+\frac12+\max\{0,2(\lceil a_j\rceil-a_j)-1\}}
  \sum_{\substack{L\subseteq \{1,\ldots,m\} \\ \#L=j}} \prod_{\substack{K\subseteq \{1,\ldots,m\} \\ \#K=j \\ K\neq L}}\left(1+ \Big| \sum_{\ell\in L}\alpha_\ell - \sum_{k\in K}\alpha_k \Big|\right)^{-[a_j]}.
\end{multline*}
\end{lemma}
\begin{proof}
Let $\mathcal{I}_j$ denote the integral we are seeking to bound.  

The polynomial part (see Definition~\ref{def:poly-exp-part}) 
of the Gamma functions in $\mathcal{I}_j$ is
\begin{multline*}
 \lvert \mathcal{Q}_j(s,z,\alpha) \rvert  
 \ll \big(1+\im(\beta_j-\tfrac{n}{n-1}\halpha)\big)^{\varepsilon+R+\frac12}\big(1+\lvert \im(s_j-z_j)\rvert\big)^{\lceil a_j \rceil-a_j-\re(z_j)-\frac12} \\ \cdot
 \big(1+\lvert\im(s_j-z_{j-1})\rvert\big)^{\lceil a_j \rceil-a_j-\re(z_{j-1})-\frac12},
\end{multline*}
and the exponential factor (when taking all $\mathcal{I}_j$ in unison) is negative for any $s_j$ outside of the interval $I_j$ defined in Lemma~\ref{lem:ExpZeroSet}.  That lemma together with the presence of the other exponential terms in our integral allow us to take trivial bounds for the polynomial part, namely that $\mathcal{Q}_j(s,z,\alpha) \ll T^{\varepsilon+R+\frac12+\max\{0,2(\lceil a_j\rceil-a_j)-1\}}$. (Recall that $0\leq \re(z_j)$.)  Thus we see that
\begin{align*}
\mathcal{I}_j & \ll
 T^{\eps+R+\frac12+\max\{0,2(\lceil a_j\rceil-a_j)-1\}} \hskip -6pt
 \int\limits_{\substack{\re(s_j)=\lceil a_j \rceil-a_j\\\im(s_j)\in I_j}} \prod\limits_{\substack{ J\subseteq \{1,\ldots,n\} \\ \#J = j}} \Big| s_j + \sum\limits_{k\in J} \alpha_k\Big|^{-\lceil a_j \rceil}\, ds_j,
\end{align*}
The desired result now follows easily from this and the statement of Lemma~\ref{lem:simpleint}.
\end{proof}

Combining Lemma~\ref{lem:IjBound} with the bound for $\mathcal{I}_{T,R}^{(n)}(-a)$ given immediately before the statement of the lemma, and applying Lemma~\ref{lem:alphasumk}, Lemma~\ref{lem:prodGammaRkdecomp} and Lemma~\ref{lem:FRnkdecomp} (in the case that $k=n-1$ and $\gamma_1=0$), we now have the bound
\begin{align*}
 \mathcal{I}_{T,R}^{(m)}(-a) & \ll
 \sum_{\substack{L\subseteq \{1,\ldots,m\} \\ \#L=j}}  T^{\eps+(R+\frac12)(m-1)+\sum\limits_{j=1}^{m-1}\Big( \max\{0,2(\lceil a_j\rceil-a_j)-1\}+ \lceil a_j \rceil\left( \binom{m}{j}-2 \right)\Big)}
 \cdot \hskip -6pt
 \int\limits_{\re(\halpha)=0} e^{\frac{m}{m-1}\frac{\halpha^2}{2T^2}}
 \\ & \hskip 24pt 
 \cdot 
 \int\limits_{\substack{\hbeta_{m-1}=0\\ \re(\beta)=0}}
 e^{\frac{\beta_1^2+\cdots+\beta_{m-1}^2}{T^2/2}}  
 \hskip -3pt
 \cdot
 \mathcal{P}_{D(m)-D(m-1)}^R(\halpha,\beta)
 \cdot
 \FR{m-1}(\beta) 
 \hskip -3pt
 \prod_{1\leq j\neq k \leq m-1} \hskip -4pt
 \big| \Gamma_R(\beta_j-\beta_k)\big|  
 \\ & \hskip 24pt \cdot
 \prod_{j=1}^{m-1} 
 \prod_{\substack{K\subseteq \{1,\ldots,m\} \\ \#K=j \\ K\neq L}}\left(1+ \Big| \sum_{\ell\in L}\alpha_\ell - \sum_{k\in K}\alpha_k \Big|\right)^{-[a_j]}
 \int\limits_{\substack{z=(z_1,\ldots,z_{m-2})\\\re(z)=b}} 
 \left| \widetilde{W}_{m-1,\beta}\left(z \right)\right|
 \; dz\, d\beta\, d\halpha.
\end{align*}
To be more explicit, the polynomial $\mathcal{P}_{D(m)-D(m-1)}^R(\halpha,\beta)$ is the portion of $\FR{m}(\alpha)$ which involves the terms $\alpha_m$.

At this point, we combine each of the terms in the final row with the corresponding term in $\FR{m}(\alpha)$.  Strictly speaking, what is actually happening here is that this has the effect of reducing the power of each factor of $\FR{m}(\alpha)$ by at most
 \[ \max\{\lceil a_1 \rceil,\ldots,\lceil a_{m-1} \rceil \}. \]
Since each of the corresponding exponents remains positive, the net result is to reduce the overall power of $T$ by 
 \[ \eps + \sum_{j=1}^{m-1} \lceil a_j \rceil \left( \binom{m}{j}-1 \right).\]
Using this, and accounting for the integration in $\halpha$ (which may be assumed to take place only for $\lvert \im(\halpha)\rvert \leq T^{1+\eps}$), we now may write
\begin{align*}
 \mathcal{I}_{T,R}^{(m)}(-a) \ll & \
 T^{\eps+(R+\frac12)(n-1)+R(D(m)-D(m-1))+1+\sum\limits_{j=1}^{m-1}\big( \max\{0,2(\lceil a_j\rceil-a_j)-1\}- \lceil a_j \rceil\big)}
 \int\limits_{\substack{\hbeta_{m-1}=0\\ \re(\beta)=0}} 
 \\ & \hskip 6pt 
 \cdot 
  e^{\frac{\beta_1^2+\cdots+\beta_{m-1}^2}{2T^2}}  
 \hskip -3pt
 \cdot
 \FR{m-1}(\beta) 
 \hskip -3pt
 \prod_{1\leq j\neq k \leq m-1} \hskip -4pt
 \big| \Gamma_R(\beta_j-\beta_k)\big|  
 \int\limits_{\substack{z=(z_1,\ldots,z_{m-2})\\\re(z)=b}}
 \left| \widetilde{W}_{m-1,\beta}\left(z \right)\right|
 \; dz\, d\beta.
\end{align*}
Obviously, at this point we want to apply the inductive hypothesis.  Since at this point we only need to do so in the case that $b_j>0$ (i.e., $-a_j<0$) for all $j=1,\ldots,m-2$, the reduction in the powers of the exponents of any one of the factors of $\FR{m}(\alpha)$, as occurred above, leaves the overall power positive.  Therefore, there is no issue, and we can assert (additionally applying Lemma~\ref{lem:max-simplify}) the bound 
\begin{align*}
 \mathcal{I}_{T,R}^{(m)}(-a) \ll & \
 T^{\eps+(R+\frac12)(m-1)+R(D(m)-D(m-1))+1+ 
 A(m-1)}
  \cdot
 T^{R\big(D(m-1)+\frac{(m-1)(m-2)}{2}\big)-\sum\limits_{j=1}^{m-1} B(a_j) }
 \\ = & \
 T^{\eps+R\big( D(m)+\frac{m(m-1)}{2} \big) + \frac{n+1}{2} + A(m-1) - \sum\limits_{j=1}^{m-1} B(a_j) }.
\end{align*}
Taking   
 $ A(m) = \tfrac{m+1}{2}+A(m-1) $
gives the claimed bound.  Since $A(2)=\frac32$, it follows that 
 \[ A(3)=\tfrac42 +A(2)=\tfrac12(4+3),\ldots,A(m)=\tfrac12\big((m+1)+m+\cdots+3\big)=\tfrac{(m+4)(m-1)}{4}, \]
as claimed.
\end{proof}

In the course of proving Theorem~\ref{th:ITRnaBound} we also established the following result that we record here since it will be useful in its own right.

\begin{proposition}\label{prop:ITR-minusa-plusL}
Suppose that $L=(\ell_1,\ell_2,\ldots,\ell_{m-1})\in \big( \Z_{\geq 0}\big)^{m-1}$.  Then 
\begin{multline*}
 \int\limits_{\substack{\halpha_m=0\\\re(\alpha)=0}}
 e^{\frac{\alpha_1^2+\cdots+\alpha_m^2}{2T^2}} 
 \cdot 
 \FR{m}(\alpha) \hskip -3pt
 \prod_{1\leq j\neq k \leq m} \hskip -4pt
 \big| \Gamma_R(\alpha_j-\alpha_k)\big|
 \int\limits_{\substack{s=(s_1,\ldots,s_{m-1})\\\re(s)=-a}} \hskip -5pt
 \left| \widetilde{W}_{m,\alpha}(s+L)\right|\, ds\, d\alpha
 \\
 \ll T^{\varepsilon+2\lvert L \rvert} \cdot 
 \int\limits_{\substack{\halpha_m=0\\\re(\alpha)=0}} 
 e^{\frac{\alpha_1^2+\cdots+\alpha_m^2}{2T^2}}
 \cdot 
 \FR{m}(\alpha) \hskip -3pt
 \prod_{1\leq j\neq k \leq m} \hskip -4pt
 \big| \Gamma_R(\alpha_j-\alpha_k)\big|
 \int\limits_{\substack{s=(s_1,\ldots,s_{m-1})\\\re(s)=-a}} \hskip -5pt
 \left| \widetilde{W}_{m,\alpha}(s)\right|\, ds\, d\alpha.
\end{multline*}
As a shorthand for this result, we write $\mathcal{I}_{T,R}^{(m)}(-a+L) \ll T^{\varepsilon+2\lvert L\rvert}\cdot \mathcal{I}_{T,R}^{(m)}(-a)$.
\end{proposition}

\section{\bf\large Bounding \boldmath $\pTR{n}(y)$}\label{sec:pTRnBound}

In this section we prove the following.

\begin{theorem}\label{th:pTRn-bound}
Let $n\geq 2$ and $\varepsilon\in(0,\frac14)$.  Suppose that $a=(a_1,a_2,\ldots,a_{n-1})$ satisfies $\lfloor a_j \rfloor+\eps < a_j < \lceil a_j \rceil -\eps$ for each $j=1,\ldots,n-1$.  Let $\mathcal{C}$ be the set of compositions  $n=n_1+\cdots+n_r$ with $r\geq 2$.  Then, for 
 \[ \Delta_a(C):=\Big\{ \delta_C=(\delta_1,\ldots,\delta_{r-1})\in \Z^{r-1} \ \Big|\ 0\leq \delta_j<a_{\hn_j}\ (j=1,\ldots,r-1) \Big\}, \]
 {and $B(c)$ as defined in Theorem~\ref{th:ITRnaBound},} we have
\begin{equation}\label{eq:final-pTRnbound}
 \left\lvert p_{T,R}^{(n)}(y) \right\rvert \ll \left\lvert p_{T,R}^{(n)}(y;-a) \right\rvert + \sum_{C\in \mathcal{C}} \sum_{\delta_C\in \Delta_a(C)} \left\lvert p_{T,R}^{(n)}(y;-a,\delta_C) \right\rvert ,
\end{equation}
where
\begin{equation}\label{eq:final-pTRna-bound}
 \left| p_{T,R}^{(n)}(y;-a) \right| \ll \prod_{j=1}^{n-1} y_j^{\frac{n(n-j)}{2}+2a_j} \cdot T^{\eps+\frac{(n+4)(n-1)}{4}+ \frac{R}{2}\cdot\left(\binom{2n}{n}-2^n \right) - \sum\limits_{j=1}^{n-1}B( a_j)},
\end{equation}
and
\begin{multline}\label{eq:final-pTRna-deltaC-bound}
 \big\lvert p_{T,R}^{(n)}(y;-a,\delta_C) \big\rvert \ll \prod_{j=1}^{n-1} y_j^{\frac{n(n-j)}{2}+2a_j} \cdot T^{\eps+\frac{(n+4)(n-1)}{4}+ \frac{R}{2}\cdot\left(\binom{2n}{n}-2^n \right) - \sum\limits_{j=1}^{n-1}B( a_j)-\frac12\sum\limits_{k=1}^{r-1} (n_k+n_{k+1})(a_{\hn_k}-\delta_k)}.
\end{multline}
The implicit constant depends on both $\varepsilon$ and $n$.
\end{theorem}

 {
\begin{remark}\label{rmk:pTRnBound}
Note that \eqref{eq:final-pTRna-deltaC-bound} is bounded by \eqref{eq:final-pTRna-bound}.  Therefore, letting $D(n) =  \frac12 \binom{2n}{n} - \frac{n(n-1)}{2} - 2^{n-1}$ as in \eqref{eq:Dn-explicit}, Theorem~\ref{th:pTRn-bound} implies that
 \[ \left\lvert p_{T,R}^{(n)}(y) \right\rvert \ll \prod_{j=1}^{n-1} y_j^{\frac{n(n-j)}{2}+2a_j} \cdot T^{\eps+\frac{(n+4)(n-1)}{4}+ R\cdot\big(D(n) +\frac{n(n+1)}{2} \big) - \sum\limits_{j=1}^{n-1}B( a_j)} \]
\end{remark}
}

\subsection{Explicit single residue formula}

In order to bound the terms $p_{T,R}^{(n)}(y;-a,\delta_C)$ we need an explicit formula for the residues of the Mellin transform of the $\GL(n)$ Whittaker function.  The following result establishes this for the case of single residues (i.e., when the composition $C$ has length $2$) as a corollary of Conjecture~\ref{conj} combined with a theorem of Stade \cite{Stade2001} for the ``first'' residues, i.e., for those residues corresponding, in the notation of the theorem, to $\delta=0$.
\begin{theorem}\label{th:residue-description}
Let $\what{m}{\al}{s}$ be the Mellin transform of the Whittaker function on $\GL(n,\R)$ with purely imaginary parameters $\alpha=(\alpha_1,\ldots,\alpha_n)$ in general position.  Let $\sigma\in S_n$ act on $\alpha$ via
 \[ \sigma\cdot \alpha := ( \alpha_{\sigma(1)},\alpha_{\sigma(2)},\ldots,\alpha_{\sigma(n)}). \]
The poles of $\what{n}{\al}{s}$ occur, for each $1\leq m \leq n-1$, at
 \[ s_m \in \big\{ -\sigma\cdot\halpha_m - \delta \, \big| \, \sigma\in S_n, \delta\in \Z_{\geq 0}, \big\}. \]
The residue at $s_m=-\halpha_m-\delta$ is equal to a sum over shifts $L=(\ell_1,\ell_2,\ldots,\ell_{n-1})$ of terms of the form
\begin{multline*}
 \prod_{\substack{K\subseteq \{1,2,\ldots,n\} \\ \#(K\cap \{1,2,\ldots,m\})\neq m-1 \\ \#K=m}} \left( \Big(\sum_{i\in K} \alpha_i \Big) - \halpha_m-\delta\right)_\delta^{-1} \left(\prod_{i=1}^m\prod_{j=m+1}^n \Gamma(\alpha_j-\alpha_i-\delta)\right)
 \\
 \cdot \mathcal{P}_{\left(\binom{n}{m}-2\right)\delta-2\lvert L\rvert}(s,\alpha)\what{m}{\beta}{s'+L'}\what{n-m}{\gamma}{s''+L''},
\end{multline*}
where
\begin{equation}\label{eq:sprimevariables}
 s'=\left.\Big( s_j+\tfrac{j}{m}\halpha_m\Big)\right|_{1\leq j \leq m}, \quad \quad
 s''=\left.\Big( s_{m+j}+\tfrac{n-m-j}{n-m}\halpha_m\Big)\right|_{1\leq j \leq n-m},
\end{equation}
with $L'=(\ell_1,\ldots,\ell_{m-1})$ and $L''=(\ell_{m+1},\ldots,\ell_{n-1})$ being the portion of $L$ corresponding to $s'$ and $s''$ respectively.  It is the case that $\ell_{m-1}=\ell_{m+1}=0$.  Note that we take as definition that $\widetilde{W}_{1}:=1$.  The same formula holds for the residue at $s_m=-\sigma\cdot \halpha_m-\delta$ by replacing each instance of $\alpha_j$ with $\alpha_{\sigma(j)}$.
\end{theorem}

\begin{remark}\label{rmk:alt-residue-theorem}
Another way of writing the above expression for the residue would be to take the product over all $K\subseteq \{1,\ldots,n\}$ with $\# K=m$ and replace $\Gamma(\alpha_{j}-\alpha_{i}-\delta)$ with $\Gamma(\alpha_{j}-\alpha_{i})$.  The two versions are equivalent because if $K\setminus \{1,\ldots,m\} = \{j\}$, then $\{1,\ldots,m\}\setminus K =\{k\}$ and
 \[ \left( \Big(\sum_{i\in K} \alpha_i \Big) - \halpha_m-\delta\right)_\delta^{-1}\Gamma(\alpha_j-\alpha_k) = \Gamma(\alpha_j-\alpha_k-\delta). \] 
\end{remark}

\begin{proof}[Sketch of proof]
In the case that $\delta=0$, this result (for $L=(0,\ldots,0)\in \C^{n-1}$) agrees with  {\cite[Theorem~3.1]{Stade2001}}.  If $\delta>0$, we need to first apply Conjecture~\ref{conj} to rewrite the expression for $\what{n}{\alpha}{s}$ around $s_m=-\alpha_m-\delta$ as a sum over shifts $L=(\ell_1,\ldots,\ell_{n-1})\in (\Z_{\geq 0})^{n-1}$ (with $\ell_m\geq \delta$ for each $L$) of terms $\what{n}{\alpha}{s+L}$.  Of all of these terms, the only ones for which there is a pole at $s_m=-\halpha_m-\delta$ are those for which $\ell_m=\delta$, in which case we can use the above referenced theorem of Stade to write down the residue.  Doing so, we obtain the alternate expression referenced to in Remark~\ref{rmk:alt-residue-theorem}.
\end{proof}

\subsection{Explicit higher residue formulae}\label{sec:higherOrderResidues}

In order to generalize Theorem~\ref{th:residue-description}, we first establish notation related to the $(r-1)$-fold residue of $\what{n}{\alpha}{s}$ at
 \[ s_{\hn_\ell} = -\halpha_{\hn_\ell} -\delta_{\hn_\ell}, \qquad \ell=1,\ldots,r-1. \]
To this end, let $s^{(j)}:=(s_1^{(j)},\ldots,s_{n_j-1}^{(j)})$ where $s_k^{(j)}=s_{\hn_{j-1}+k}$.  By abuse of notation, we write
 \[ s:= \big(
 \underbrace{s_1^{(1)},s_2^{(1)},\ldots, s_{n_1-1}^{(1)}}_{=:s^{(1)}}\ ,\
 \underbrace{s_1^{(2)},s_2^{(2)},\ldots, s_{n_2-1}^{(2)}}_{=:s^{(2)}}\ ,\ \ldots\ ,\
 \underbrace{s_1^{(k)},s_2^{(k)},\ldots, s_{n_k-1}^{(k)}}_{=:s^{(k)}}\big) \in \C^{n-r}, \]
which agrees with the original $s=(s_1,\ldots,s_{n-1})$ but removes $s_{\hn_1},\ldots,s_{\hn_{r-1}}$.  

Similarly, if $\alpha=(\alpha_1,\ldots,\alpha_n)$, we define
 \[ \alpha^{(\ell)} := \big( \alpha_1^{(\ell)},\ldots,\alpha_{\ell}^{(\ell)}\big) \in \C^{n_\ell}, \qquad \alpha_j^{(\ell)} := \alpha_{\hn_{\ell-1}+j}-\tfrac{1}{n_\ell}\big(\halpha_{\hn_\ell}-\halpha_{\hn_{\ell-1}}\big), \]
and
 \[ \big| \al^{(j)} \big|^2 := \big(\al_1^{(j)}\big)^2+\big(\al_2^{(j)}\big)^2+\cdots +\big(\al_{n_j}^{(j)}\big)^2. \]
If $a\in \R^{n-1}$ then by $\re(s)=-a$ we mean that $\re(s_j)=-a_j$ for each $j\neq \hn_1,\ldots,\hn_{r-1}$.

\hskip 0pt

With this notation in place, we can now state a generalization of Theorem~\ref{th:residue-description}.
\begin{corollary}\label{cor:residue-description}
Let $n=n_1+\cdots +n_r$ ($r\geq 2$), and set $\hn_\ell :=\sum\limits_{j=1}^\ell n_j$ as above.  For each $\ell=1,\ldots,r-1$, let $b^{(\ell)}=(b_1^{(\ell)},b_2^{(\ell)}, \ldots, b_{n_\ell-1}^{(\ell)})$ with
 \[ b_j^{(\ell)} = \halpha_{i_{\ell-1}} + \tfrac{j}{n_{\ell}}\big( \halpha_{\hn_\ell}-\halpha_{\hn_{\ell-1}} \big) \qquad \mbox{for each $1\leq j \leq n_\ell-1$}. \]
Let $\delta_j\in \Z_{\geq 0}$ for $j=1,\ldots,r-1$.  There exist positive shifts $L=(L^{(1)},\ldots,L^{(r)})$ with $L^{(\ell)}=(L_1^{(\ell)},\ldots,L_{n_\ell-1}^{(\ell)})\in \big(\Z_{\geq 0} \big)^r$ such that the iterated residue of $\what{n}{\alpha}{s}$ at 
 \[ s_{\hn_{r-1}}=-\halpha_{\hn_{r-1}}-\delta_{r-1}\ ,\ \ldots\ ,\ s_{\hn_1}=-\halpha_{\hn_1}-\delta_1, \]
is equal to a sum over all such shifts of
\begin{multline*}
\mathcal{P}_d(s,\alpha) \left(\prod_{\ell=1}^r
\what{n_\ell}{\alpha^{(\ell)}}{s^{(\ell)}+b^{(\ell)}+L^{(\ell)}} \right)
\prod_{j=1}^{r-1} \prod_{\substack{K\subseteq \{1,2,\ldots,\hn_{j+1}\} \\ \#(K\cap \{1,\ldots,\hn_j\})\neq \hn_j-1 \\ \#K=\hn_j}} \left( \Big(\sum_{i\in K} \alpha_i \Big) - \halpha_{\hn_j}-\delta_j\right)_{\delta_j}^{-1}
\\
\prod_{1\leq k < m \leq r}\prod_{i=1}^{n_k} \prod_{j=1}^{n_m} \Gamma\bigg(\alpha_j^{(m)}-\alpha_i^{(k)}+\tfrac{1}{n_m}\big(\halpha_{\hn_m}-\halpha_{\hn_{m-1}}\big)-\tfrac{1}{n_k}\big(\halpha_{\hn_k}-\halpha_{\hn_{k-1}}\big)-\delta_m\bigg),
\end{multline*}
where
 \[ d = \left[ \sum_{\ell=1}^{r-1} \delta_\ell\left( \binom{\hn_{\ell+1}}{\hn_\ell} - 2 \right) \right] - 2\lvert L \rvert. \]
\end{corollary}
\begin{proof}
This follows easily by induction with the base case being Theorem~\ref{th:residue-description}.
\end{proof}

\begin{remark}
Although it is possible to rewrite each of the terms $\bigg(\sum\limits_{i\in K} \alpha_i\bigg)-\halpha_{\hn_\ell}-\delta_\ell$ appearing in the statement of Corollary~\ref{cor:residue-description} in terms of the variables $\alpha^{(j)}$ and $\halpha_m^{(j)}$ for various $j$ and $m$, the exact description is unnecessary for our purposes.
\end{remark}

\subsection{Proof of Theorem~\ref{th:pTRn-bound}}\label{sec:pTRn-bound}

\begin{proof}
As a first step, note that Proposition~\ref{prop:pTR-expansion} implies that \eqref{eq:final-pTRnbound} follows from \eqref{eq:final-pTRna-bound} and \eqref{eq:final-pTRna-deltaC-bound}.

As shown in Section~\ref{sec:ITRMa}, the shifted $p_{T,R}^{(n)}$ term satisfies
\begin{align*}
 \big| p_{T,R}^{(n)}(y,-a) \big| & \ll
 \bigg(\prod_{j=1}^{n-1}
 y_j^{\frac{j(n-j)}{2}-2a_j} \bigg) \mathcal{I}_{T,R}^{(n)}(-a).
\end{align*}
Combined with the bound from Theorem~\ref{th:ITRnaBound}, this gives \eqref{eq:final-pTRna-bound}.

To complete the proof, we need to show that \eqref{eq:final-pTRna-deltaC-bound} holds.  We do this in Section~\ref{sec:multi-residue-bound}.  Although this proof is valid for any $r\geq 2$, as a warmup, we first prove the special case $r=2$ (i.e., the case of single residues) in Section~\ref{sec:single-residues}.
\end{proof}

\subsection{Bounds for single residue terms}\label{sec:single-residues}

In this section\footnote{Note that this section will be superseded by Section~\ref{sec:multi-residue-bound} which will prove the bound for any admissible $C$ with $\len(C)\geq 2$.  This section treats the case $\len(C)=2$.} we bound $p_{T,R}^{(n)}(y;-a,\delta_C)$ in the case that of the composition $C=(m,n-m)$.  Since $C$ is a composition of length two, we may take (see Definition~\ref{def:k-fold-residue-term}) $\delta_C=\delta \in \Z_{\geq 0}$.

\begin{proof}[Proof of \eqref{eq:final-pTRna-deltaC-bound} when $r=2$]
Using Lemmas~\ref{lem:alphasumk}, \ref{lem:prodGammaRkdecomp} and \ref{lem:resm_extrapoly}, we can rewrite
 \[ e^{\frac{\al_1^2+\cdots\al_n^2}{T^2/2}}\mathcal{F}_R^{(n)}(\alpha) \prod_{1\leq j\neq k \leq n} \Gamma_R(\alpha_j-\alpha_k) \]
in terms of $\beta$, $\gamma$ and $\al_n$.  Thus, together with Theorem~\ref{th:residue-description}, we see that  Definition~\ref{def:k-fold-residue-term} in the case of a single residue term (i.e., $r=2$) satisfies the bound
\begin{align*}
 p_{T,R}^{(n)}(y;-a,\delta_C) & \ll \hskip -4pt 
 \int\limits_{\re(\halpha_m)=0} y_m^{\frac{m(n-m)}{2}+\halpha_m+\delta} \cdot e^{\frac{n}{m(n-m)}\frac{\halpha_m^2}{T^2/2}}
 \int\limits_{\substack{\hbeta_m=0\\\re(\beta)=0}}
 e^{\frac{|\beta|^2}{T^2/2}} 
 \int\limits_{\substack{\widehat{\gamma}_{n-m}=0\\\re(\gamma)=0}}
 e^{\frac{|\gamma|^2}{T^2/2}}
 \\ & \hskip 18pt \cdot
 \left( \FR{m}(\beta) \cdot \prod_{1\leq i\neq j\leq m} \Gamma_R(\beta_i-\beta_j)\right) \left( \FR{n-m}(\gamma) \cdot \prod_{1\leq i\neq j\leq n-m} \Gamma_R(\gamma_i-\gamma_j)\right)
 \\ & \hskip -24pt \cdot
 \prod_{i=1}^m \prod_{j=1}^{n-m} \bigg( 
 \Gamma_R\big(\beta_i-\gamma_j+\tfrac{n\halpha_m}{m(n-m)}\big)
 \Gamma_R\big(\gamma_j-\beta_i-\tfrac{n\halpha_m}{m(n-m)}\big)
 \Gamma\big( \gamma_j-\beta_i-\tfrac{n\halpha_m}{m(n-m)}-\delta \big)
 \bigg)
 \\ & \hskip 12pt \cdot
 \Bigg(\prod_{j\neq m}\int\limits_{\substack{\re(s_j)=-a_j\\j\neq m}}
 y_j^{\frac{j(n-j)}{2}-s_j} \Bigg)
 \mathcal{P}_{R\big(D(n)-D(m)-D(n-m)\big)-\delta\left( \binom{n}{m}-m(n-m)-1\right)}(s,\alpha)
 \\ & \hskip 38pt \cdot
 \mathcal{P}_{\left(\binom{n}{m}-2\right)\delta-2\lvert L\rvert}(s,\alpha) \cdot
 \what{m}{\beta}{s'+L'}\what{n-m}{\gamma}{s''+L''}
 \, ds\; d\gamma\; d\beta\; d\halpha_m.
\end{align*}
In order to have the correct power of $y_m$ we need to shift the line of integration in $\halpha_m$ to $\re(\halpha_m)=a_m-\delta$.  Note that by Lemma~\ref{lem:nopoles}, no poles are crossed in doing so, and by Lemma~\ref{lem:alphak-delta-residue-bound}, taking $\beta=\beta_i-\gamma_j$ and $z=\frac{n\halpha_m}{m(n-m)}$, we may replace the third to last line by
 \[ \mathcal{P}_{m(n-m)R-n(a_m-\delta)-m(n-m)\delta}(s,\halpha_m,\beta,\gamma). \]
Let $\lvert \beta \rvert^2 := \beta_1^2+\cdots+\beta_{m}^2$, and define $\lvert \gamma \rvert^2$ similarly.  Replacing the integral over $\halpha_m$ by $T^{\eps+1}$ and factoring out the powers of $y_j$, we see that
\begin{align*}
 \Big| p_{T,R}^{(n)}(y;-a,\delta_C) \Big| & \ll
 \bigg( \prod_{j=1}^{n-1} y_j^{\frac{j(n-j)}{2}+a_j} \bigg)
 \cdot T^{\eps+\left(\binom{n}{m}-2\right)\delta+R\big(D(n)-D(m)-D(n-m)\big)-\delta\left( \binom{n}{m}-m(n-m)-1\right)}
 \\ & \hskip 12pt
 \cdot
 T^{-2\lvert L\rvert+m(n-m)R-n(a_m-\delta)-m(n-m)\delta+1}
 \cdot 
 \int\limits_{\substack{\hbeta_m=0\\\re(\beta)=0}}
 e^{\frac{\lvert \beta \rvert^2}{T^2/2}} 
 \int\limits_{\substack{\widehat{\gamma}_{n-m}=0\\\re(\gamma)=0}}
 e^{\frac{\lvert \gamma \rvert^2}{T^2/2}}
 \\ & \hskip 30pt \cdot
 \left( \FR{m}(\beta) \cdot \prod_{1\leq i,j\leq k} \Gamma_R(\beta_i-\beta_j)\right) \left( \FR{n-m}(\gamma) \cdot \prod_{1\leq i,j\leq n-k} \Gamma_R(\gamma_i-\gamma_j)\right)
 \\ & \hskip 72pt \cdot
 \hskip -6pt
 \int\limits_{\substack{\re(s_j)=-a_j\\ 1\leq j \leq n-1\\ j\neq m}} \Big\lvert \what{m}{\beta}{s'+L'} \Big \rvert \cdot \Big\lvert \what{n-m}{\gamma}{s''+L''} \Big\rvert
 \, ds\; d\gamma\; d\beta.
\end{align*}
Note that by Proposition~\ref{prop:ITR-minusa-plusL} we may remove the dependence on the shift $L$.  Hence
\begin{align*}
 \Big| p_{T,R}^{(n)}(y;-a,\delta_C) \Big| & \ll
 \bigg( \prod_{j=1}^{n-1} y_j^{\frac{j(n-j)}{2}+a_j} \bigg)
 \cdot T^{\eps+R\big(D(n)-D(m)-D(n-m)+m(n-m)\big)+\delta(n-1)}
 \\ & \hskip 36pt
 \cdot T^{-na_m+1} 
 \int\limits_{\substack{\hbeta_m=0\\\re(\beta)=0}}
 e^{\frac{\beta_1^2+\cdots+\beta_m^2}{T^2/2}} 
 \int\limits_{\substack{\widehat{\gamma}_{n-m}=0\\\re(\gamma)=0}}
 e^{\frac{\gamma_1^2+\cdots+\gamma_n^2}{T^2/2}}
 \\ & \hskip 0pt \cdot
 \left( \FR{m}(\beta) \cdot \prod_{1\leq i\neq j\leq k} \Gamma_R(\beta_i-\beta_j)\right) \left( \FR{n-m}(\gamma) \cdot \prod_{1\leq i\neq j\leq n-k} \Gamma_R(\gamma_i-\gamma_j)\right)
 \\ & \hskip 72pt \cdot
 \hskip -6pt
 \int\limits_{\substack{\re(s_j)=-a_j\\ 1\leq j \leq n-1\\ j\neq m}} \Big\lvert \what{m}{\beta}{s'} \Big \rvert \cdot \Big\lvert \what{n-m}{\gamma}{s''} \Big\rvert
 \, ds\; d\gamma\; d\beta.
\end{align*}

By \eqref{eq:sprimevariables}, 
 \[ s_j' = s_j - \tfrac{j}{m}(\widehat{\alpha}_m-\delta) , \quad \mbox{and}\quad
    s_j''= s_{m+j} - \tfrac{n-m-j}{n-m}(\widehat{\alpha}_m-\delta). \]
Thus the integrals in $\beta$ and $\gamma$ above are essentially the product of $\mathcal{I}_{T,R}^{(m)}(-a')$ and $\mathcal{I}_{T,R}^{(n-m)}(-a'')$.  The only issue is that because, as seen in 
the fact that the variables $s'$ and $s''$ are shifted, we have
 \[ a_j' = a_j - \tfrac{j}{m}(a_m-\delta) , \quad \mbox{and}\quad
    a_j''= a_{m+j} - \tfrac{n-m-j}{n-m}(a_m-\delta). \]
Therefore, we can rewrite the previous formula as
\begin{align*}
 \Big| p_{T,R}^{(n)}(y;-a,\delta_C) \Big| & \ll
 \bigg( \prod_{j=1}^{n-1} y_j^{\frac{j(n-j)}{2}+a_j} \bigg)
 \cdot T^{\eps+R\big(D(n)-D(m)-D(n-m)+m(n-m)\big)}
 \\ & \hskip 36pt
 \cdot T^{\delta(n-1)-na_m+1} \cdot
 \mathcal{I}_{T,R}^{(m)}(-a') \cdot  \mathcal{I}_{T,R}^{(n-m)}(-a'').
\end{align*}
By Theorem~\ref{th:ITRnaBound}, we have
\begin{align*}
 \Big| p_{T,R}^{(n)}(y;-a,\delta_C) \Big| & \ll
 \bigg( \prod_{j=1}^{n-1} y_j^{\frac{j(n-j)}{2}+a_j} \bigg)
 \cdot T^{\eps+R\big(D(n)-D(m)-D(n-m)+m(n-m)\big)}
 \\ & \hskip 12pt
 \cdot T^{\delta(n-1)-na_m+1} \cdot
 T^{\varepsilon+C(m) + R\cdot\big(D(m)+\frac{m(m-1)}{2}\big) - \sum\limits_{j=1}^{m-1}B(a_j')} 
 \\ & \hskip 12pt
 \cdot  T^{\varepsilon+C(n-m) + R\cdot\big(D(n-m)+\frac{(n-m)(n-m-1)}{2}\big) - \sum\limits_{j=1}^{n-m-1}B(a_j'')},
\end{align*}
Recall that $C(k) = \frac{(k+4)(k-1)}{4}$.  Hence, using the elementary identity
 \[ C(m)+C(n-m) = C(n) - \frac{m(n-m)}{2} - 1 \] 
together with Lemma~\ref{lem:residue-powerofT-new-aj},
\begin{align*}
 \Big| p_{T,R}^{(n)}(y;-a,\delta_C) \Big| & \ll
 \bigg( \prod_{j=1}^{n-1} y_j^{\frac{j(n-j)}{2}+a_j} \bigg)
 \cdot T^{\eps+R\big(D(n)+m(n-m) + \frac{m(m-1)}{2} + \frac{(n-m)(n-m-1)}{2} \big)}
 \\ & \hskip 12pt
 \cdot T^{\delta(n-1)+C(n)-\frac{m(n-m)}{2}-na_m - \sum\limits_{j=1}^{m-1}B(a_j)+\frac{n-2}{2}(a_m-\delta+1)+B(a_m)}
 \\ & \ll
 \bigg( \prod_{j=1}^{n-1} y_j^{\frac{j(n-j)}{2}+a_j} \bigg)
 \cdot T^{\eps+C(n)+R\big(D(n) + \frac{n(n-1)}{2}\big) - \sum\limits_{j=1}^{n-1} B(a_j)}
 \\ & \hskip 48pt \cdot
 T^{\frac{n-2}{2}(\delta-a_m+1)-\frac{m(n-m)}{2} - n(\delta-a_m)+ B(a_m)-\delta}
\end{align*}
This gives the desired bound provided that the exponent of the final $T$ is negative.  Using the facts that $-\frac{m(n-m)}{2}$ is maximized when $m=1$ or $m=n-1$ and $B(a_m)\leq a_m+\frac12$, we see that the final exponent is
\begin{align}\label{eq:requals2negative}
 -\frac{n}{2}\big(a_m-\delta\big)+\frac{n-1}{2}-\frac{m(n-m)}{2} & \leq\ -\frac{n}{2}\big(a_m-\delta\big), 
\end{align}
as claimed.
\end{proof}

\subsection{Bounds for $(r-1)$-fold residues}\label{sec:multi-residue-bound}

We consider a composition $C$ of $n$ of length $r\geq 2$ given by $n=n_1+\cdots+n_r$.  We may also write $C=(n_1,\ldots,n_r)$.  Let $\hn_\ell = \sum\limits_{j=1}^\ell n_j$ as usual.

As a final piece of notation, let $\beta=(\beta_1,\ldots,\beta_r)$ be defined via
 \[ \beta_i := \halpha_{\hn_i}-\halpha_{\hn_{i-1}}\]
Note that $\sum\limits_{i=1}^r \beta_i =0$ and more generally, defining $\hbeta_m=\sum\limits_{i=1}^m\beta_i$, $\halpha_{\hn_i}=\hbeta_i$.  Since (assuming that $\halpha_n=0$) the Jacobians of the change of variables 
 \[ \alpha\mapsto (\alpha^{(1)},\halpha_{\hn_1},\alpha^{(2)},\halpha_{\hn_2},\ldots,\halpha_{\hn_{r-1}},\alpha^{(r)}) \]
and  
 \[ (\halpha_{\hn_1},\ldots, \halpha_{\hn_{r-1}}) \mapsto (\beta_1,\ldots,\beta_{r-1}) \]
are trivial, we see that (for $\beta_1+\cdots+\beta_r=0$)
\begin{equation}\label{eq:change-of-variables}
 d\alpha = d\beta \; d\alpha^{(1)} d\alpha^{(2)} \cdots d\alpha^{(r)}.
\end{equation}

\begin{proof}[Proof of \eqref{eq:final-pTRna-deltaC-bound} when $r\geq 2$]
Note that
 \[ b_j^{(\ell)} = \halpha_{\hn_{\ell-1}} + \tfrac{j}{n_{\ell}}\big( \halpha_{\hn_\ell}-\halpha_{\hn_{\ell-1}} \big) = \hbeta_{\ell-1}+\tfrac{j}{n_\ell}\beta_\ell \qquad \mbox{for each $1\leq j \leq n_\ell-1$}. \]
Recall that by Definition~\ref{def:k-fold-residue-term},
\begin{align*}
p_{T,R}^{(n)}(y;-a,\delta_C) & := \hskip -6pt
 \int\limits_{\substack{\halpha_n=0\\\re(\alpha)=0}}
 e^{\frac{\alpha_1^2+\cdots+\alpha_n^2}{T^2/2}} 
 \cdot 
 \FR{n}(\alpha)
 \bigg(\prod_{1\leq j\neq k \leq n} \hskip -4pt
 \Gamma_R\big(\alpha_j-\alpha_k\big) \bigg) \hskip -3pt
 \\ \nonumber & \hskip 12pt 
 \cdot
 \bigg( \prod_{i=1}^{r-1} y_{\hn_i}^{\frac{\hn_i(n-\hn_i)}{2}+\halpha_{\hn_i}+\delta_i}\bigg)
 \cdot 
 \int\limits_{\substack{\re(s_j)=-a_j\\ j\notin \{\hn_1,\ldots,\hn_{r-1}\}}}
 \bigg( \prod_{j\notin \{\hn_1,\ldots,\hn_{r-1}\}} y_j^{\frac{j(n-j)}{2}-s_j} \bigg)
 \\ \nonumber & \hskip -24pt \cdot
 \underset{s_{\hn_1}=-\halpha_{\hn_1}-\delta_1}{\res} 
    \left(\underset{s_{\hn_2}=-\halpha_{\hn_2}-\delta_2}{\res}\left( \cdots
    \left(\underset{s_{\hn_{r-1}}=-\halpha_{\hn_{r-1}}-\delta_{r-1}}{\res} \what{n}{\alpha}{s}\right) \cdots \right) \right)\, ds\, d\alpha.
\end{align*}

Using Remark~\ref{rmk:gen_res_extrapoly} and Corollary~\ref{cor:residue-description}, we can bound $\big\lvert p_{T,R}^{(n)}(y;-a,\delta_C) \big\rvert$ by a sum over certain shifts $L$ each of the form
\begin{align*}
 & 
 \int\limits_{\substack{\hbeta_r=0\\\re(\beta)=0}}
 e^{\big( \frac{\beta_1^2}{n_1}+\cdots+\frac{\beta_r^2}{n_r} \big)\frac{2}{T^2}}
 \cdot \Bigg(
 \prod_{j=1}^{r-1} \hskip 6pt
 y_{\hn_j}^{\frac{\hn_j(n-\hn_j)}{2}+\hbeta_j+\delta_j}
 \int\limits_{\substack{\halpha^{(j)}_{n_j}=0 \\ \re(\alpha^{(j)})=0}}
 e^{\frac{\left\lvert \alpha^{(j)}\right\rvert^2}{T^2/2}} \Bigg)
 \\  & \hskip 108pt 
 \cdot 
 \mathcal{P}_{d_1-2\lvert L\rvert}(\alpha)
 \cdot 
 \int\limits_{\substack{\re(s)=-a
 }}
 \bigg( \prod_{j\notin \{\hn_1,\ldots,\hn_{r-1}\}} y_j^{\frac{j(n-j)}{2}-s_j} \bigg) \cdot \mathcal{P}_{d_2}(s,\alpha)
\\ & \hskip 12pt
 \cdot
 \prod_{1\leq k < m \leq r}\prod_{i=1}^{n_k} \prod_{j=1}^{n_m} \Gamma\bigg(\alpha_j^{(m)}-\alpha_i^{(k)}+\tfrac{\beta_m}{n_m}-\tfrac{\beta_k}{n_k}-\delta_k\bigg) \prod_{\epsilon\in\{\pm1\}} \Gamma_R\bigg(\epsilon\Big(\alpha_j^{(m)}-\alpha_i^{(k)}+\tfrac{\beta_m}{n_m}-\tfrac{\beta_k}{n_k}\Big)\bigg)
 \\ & \hskip 90pt \cdot
 \prod_{\ell=1}^{r} \left( \mathcal{F}_R^{(n_\ell)}\big(\alpha^{(\ell)}\big)
 \bigg(\prod_{1\leq j\neq k \leq n_\ell} \hskip -4pt
 \Gamma_R\big(\alpha_j^{(\ell)}-\alpha_k^{(\ell)}\big) \bigg)
\what{n_\ell}{\alpha^{(\ell)}}{s^{(\ell)}+b^{(\ell)}+L^{(\ell)}} \right) 
\\ & \hskip 324pt
ds\ d\alpha^{(1)}d\alpha^{(2)}\cdots d\alpha^{(r)}\ d\beta,
\end{align*}
where
 \[ d_1 = \sum_{\ell=1}^{r-1} \delta_\ell\left( \binom{\hn_{\ell+1}}{\hn_\ell} - 2 \right) \]
and
 \[ d_2 = R\cdot \left( D(n) - \sum_{\ell=1}^r D(n_\ell)\right) - \sum\limits_{\ell=1}^{r-1} \left[ \delta_\ell\left( \binom{\hn_{\ell+1}}{\hn_\ell} - n_{\ell+1} \hn_\ell - 1\right)\right] \]
are the degrees coming from Remark~\ref{rmk:gen_res_extrapoly} and Corollary~\ref{cor:residue-description}, respectively, and $b^{(\ell)}$ is as in Corollary~\ref{cor:residue-description}.  Note that, in addition to using the change of variables \eqref{eq:change-of-variables}, we have used Lemma~\ref{lem:generalalphasum} and  Lemma~\ref{lem:genprodGammaRkdecomp} to break up $e^{2\lvert \alpha \rvert^2/T^2}$ and rewrite the product of $\Gamma(\alpha_j-\alpha_k)$ in terms of $\alpha^{(1)},\ldots,\alpha^{(r)}$ and $\beta$.

The next step is to shift the lines of integration in the variables $\beta_j$ for $j=1,\ldots,r-1$ (or, equivalently, $\hbeta_j$ for $j=1,\ldots,r-1$) such that the real part of the exponent of each term $y_{\hn_j}$ is $\frac{\hn_j(n-\hn_j)}{2}+a_j$.  In particular, this implies that we must shift the line of integration of $\hbeta_j$ to
\begin{equation}\label{eq:realOfBetaj}
 \re(\hbeta_j) = a_{\hn_j}-\delta_j \Longleftrightarrow \re(\beta_j) = \re(\hbeta_j-\hbeta_{j-1}) = (a_{\hn_j}-\delta_j)-(a_{\hn_{j-1}}-\delta_{j-1}).
\end{equation}
Provided that $R$ is sufficiently large, Lemma~\ref{lem:nopoles} implies that this shift can be made without passing any poles.  Moreover, Lemma~\ref{lem:alphak-delta-residue-bound} implies that
\begin{multline}\label{eq:extraGammaPoly}
 \prod_{1\leq k < m \leq r}\prod_{i=1}^{n_k} \prod_{j=1}^{n_m} \Gamma\bigg(\alpha_i^{(k)}-\alpha_j^{(m)}+\tfrac{\beta_k}{n_k}-\tfrac{\beta_m}{n_m}-\delta_m\bigg) \prod_{\epsilon\in\{\pm1\}} \Gamma_R\bigg(\epsilon\Big(\alpha_i^{(k)}-\alpha_j^{(m)}+\tfrac{\beta_k}{n_k}-\tfrac{\beta_m}{n_m}\Big)\bigg)
 \\ 
 \asymp
 \prod_{1\leq k < m \leq r}\prod_{i=1}^{n_k} \prod_{j=1}^{n_m} \Big(1+\big|\im\big(\alpha_i^{(k)}-\alpha_j^{(m)}+\tfrac{\beta_k}{n_k}-\tfrac{\beta_m}{n_m}\big)\big|\Big)^{R-\re\big( \tfrac{\beta_k}{n_k}-\tfrac{\beta_m}{n_m} \big)-\delta_m} 
\end{multline}
Note that the presence of the term $e^{\big( \frac{\beta_1^2}{n_1}+\cdots+\frac{\beta_r^2}{n_r} \big)\frac{2}{T^2}}$ implies that there is exponential decay for $\lvert \im(\beta_j)\rvert \gg T^{1+\eps}$.  As we will see momentarily, besides the polynomial terms $\mathcal{P}_{d_1}(\alpha)$, $\mathcal{P}_{d_2}(s,\alpha)$ and \eqref{eq:extraGammaPoly}, we just get a product of $\mathcal{I}^{(n_j)}_{T,R}(-c^{(j)})$ for some (to be determined) values $-c^{(\ell)}$.  The upshot is that all of these polynomials can be bounded by $T$ to the degree of the polynomial plus $\eps$.  Hence, we can bound the expression above by
\begin{multline}\label{eq:residuebound-withL}
 T^{\eps+r-1+d-2\lvert L\rvert}
 \cdot 
 \left( \prod_{j=1}^{n-1} \hskip 6pt
 y_{j}^{\frac{j(n-j)}{2}+a_j} \right)
 \cdot
 \prod_{\ell=1}^r \Bigg(
 \int\limits_{\substack{\halpha^{(\ell)}_{n_\ell}=0 \\ \re(\alpha^{(\ell)})=0}}
 e^{\frac{\left\lvert \alpha^{(\ell)}\right\rvert^2}{T^2/2}}
 \cdot
 \mathcal{F}_R^{(n_\ell)}\big(\alpha^{(\ell)}\big)
 \int\limits_{\substack{ \re(s^{(\ell)})=-a^{(\ell)} }}
\\ 
 \cdot
 \bigg(\prod_{1\leq j\neq k \leq n_\ell} \hskip -4pt
 \Gamma_R\big(\alpha_j^{(\ell)}-\alpha_k^{(\ell)}\big) \bigg)
\Big| \what{n_\ell}{\alpha^{(\ell)}}{s^{(\ell)}+b^{(\ell)}+L^{(\ell)}} \Big|\ ds^{(\ell)}\, d\alpha^{(\ell)}\Bigg),
\end{multline}
where $d=d_1+d_2+d_3$ with $d_1$ and $d_2$ as above and 
 \[  d_3 = R\cdot\sum_{\ell=1}^r n_\ell \hn_\ell - \sum_{k=1}^{r-1}\Big(  (n_k+n_{k+1})(a_{\hn_k}-\delta_k)+\delta_k n_{k+1}\hn_k\Big) \]
is the bound coming from the terms described in \eqref{eq:extraGammaPoly}, simplified using Lemma~\ref{lem:extra-Gammasum}.  Combining everything, we find that $d$ equals
 \[ R\cdot \bigg( \hskip -3pt D(n) - \sum_{\ell=1}^{r} D(n_\ell) + \sum_{1\leq k < m \leq r} n_k n_m \hskip -3pt \bigg) - \sum_{k=1}^{r-1}\Big(\delta_k + (n_k+n_{k+1})(a_{\hn_k}-\delta_k)\Big). \]

Recall that the bound on $p_{T,R}^{(n)}(y;-a,\delta_C)$ is a \emph{sum} of expressions of the form given in \eqref{eq:residuebound-withL} for various shifts $L$.  However, using Proposition~\ref{prop:ITR-minusa-plusL}, we can remove the dependence on the shifts.  Hence,
\begin{multline}\label{eq:residuebound-withoutL}
 \big\lvert p_{T,R}^{(n)}(y;-a,\delta_C) \big\rvert \ll T^{\eps+d+r-1}
 \cdot 
 \left( \prod_{j=1}^{n-1} \hskip 6pt
 y_{j}^{\frac{j(n-j)}{2}+a_j} \right)
 \cdot
 \prod_{\ell=1}^r \Bigg(
 \int\limits_{\substack{\halpha^{(\ell)}_{n_\ell}=0 \\ \re(\alpha^{(\ell)})=0}}
 e^{\frac{\left\lvert \alpha^{(\ell)}\right\rvert^2}{T^2/2}}
\\ 
 \cdot
 \mathcal{F}_R^{(n_\ell)}\big(\alpha^{(\ell)}\big)
 \hskip -6pt
 \int\limits_{\substack{ \re(s^{(\ell)})=-a^{(\ell)} }}
 \bigg(\prod_{1\leq j\neq k \leq n_\ell} \hskip -6pt
 \Gamma_R\big(\alpha_j^{(\ell)}-\alpha_k^{(\ell)}\big) \bigg)
\Big| \what{n_\ell}{\alpha^{(\ell)}}{s^{(\ell)}+b^{(\ell)}} \Big|\ ds^{(\ell)}\, d\alpha^{(\ell)}\Bigg),
\end{multline}
Hence, setting $c^{(\ell)} = a^{(\ell)}-\re\big(b^{(\ell)}\big),$ where 
 \[ b^{(\ell)} = (b_1^{(\ell)},\ldots,b_{\hn_j}^{(\ell)}), \qquad b_j^{(\ell)} = \hbeta_{\ell-1} + \tfrac{j}{n_{\ell}}\beta_{\ell}, \]
we find that
\begin{align*}
 \big\lvert p_{T,R}^{(n)}(y;-a,\delta_C) \big\rvert \ll \left( \prod_{j=1}^{n-1} \hskip 6pt
 y_{j}^{\frac{j(n-j)}{2}+a_j} \right)
 \cdot
 T^{\eps+r-1+d}
 \cdot \prod_{\ell=1}^r \mathcal{I}_{T,R}^{(\ell)}(-c^{(\ell)}).
\end{align*} 

Let $C(m):=\frac{(m+4)(m-1)}{4}$.  We now now apply Theorem~\ref{th:ITRnaBound} to each $\mathcal{I}_{T,R}^{(n_\ell)}$ to obtain
\begin{align*}
 \big\lvert p_{T,R}^{(n)}(y;-a,\delta_C) \big\rvert \ll 
 T^{\eps+r-1+d+\sum\limits_{\ell=1}^r\Big( R\big(D(n_\ell)+\frac{n_\ell(n_\ell-1)}{2}\big) + C(n_\ell) - \sum\limits_{k=1}^{n_\ell-1}B(c_k^{(\ell)})\Big)}
 \cdot
 \prod_{j=1}^{n-1} \hskip 6pt
 y_{j}^{\frac{j(n-j)}{2}+a_j}.
\end{align*} 
Now we generalize the proof of Lemma~\ref{lem:residue-powerofT-new-aj}, keeping in mind that $a<B(a)<a+\frac12$, to simplify the expression
\begin{align*}
 \sum_{\ell=1}^r \sum_{j=1}^{n_\ell-1}B(c_j^{(\ell)})
 & \geq
 \sum_{\ell=1}^r \sum_{j=1}^{n_\ell-1}\big( a_{\hn_\ell+j}-\re(\hbeta_{\ell-1}) - \tfrac{j}{n_{\ell}}\re(\beta_{\ell}) \big) \\
 & =
 \left(\sum_{j=1}^{n-1} a_j\right)-\left(\sum_{k=1}^{r-1} a_{\hn_k}\right) - \sum_{\ell=1}^r \Big[ (n_\ell-1)\re(\hbeta_{\ell-1}) + \tfrac{n_\ell-1}{2}\re(\beta_\ell) \Big] \\
 & \geq
 \sum_{j=1}^{n-1}\Big( B(a_j)-\tfrac{1}{2}\Big)-\sum_{k=1}^{r-1} a_{\hn_k} - \sum_{\ell=1}^r \Big[ (n_\ell-1)\re\left(\hbeta_{\ell} - \tfrac{1}{2}\beta_\ell\right) \Big] \\
 & =
 -\tfrac{n-1}{2} +
 \sum_{j=1}^{n-1} B(a_j)-\sum_{k=1}^{r-1} a_{\hn_k} - \sum_{\ell=1}^r \Big[ (n_\ell-1)\left(A_\ell - \tfrac{1}{2}(A_\ell-A_{\ell-1}\right) \Big] \\
 & =
 -\tfrac{n-1}{2} +
 \sum_{j=1}^{n-1} B(a_j)-\sum_{k=1}^{r-1} a_{\hn_k} - \frac12\sum_{\ell=1}^r \Big[ (n_\ell-1)\left(A_\ell + A_{\ell-1}\right) \Big].
\end{align*}
Next, we write the sum over $\ell$ as
\begin{align*}
 \sum_{\ell=1}^r \Big[ (n_\ell-1)\left(A_\ell + A_{\ell-1}\right) \Big] 
 & =
 \sum_{\ell=1}^r (n_\ell-1)A_\ell + \sum_{\ell=1}^r (n_\ell-1) A_{\ell-1} \\
 & =
 \sum_{\ell=1}^r (n_\ell-1)A_\ell + \sum_{\ell=0}^{r-1} (n_{\ell+1}-1) A_\ell \\
 & =
 (n_1-1)A_0 + (n_r-1)A_r + \sum_{\ell=1}^{r-1} (n_\ell+n_{\ell+1}-2) A_\ell \\
 & =
 \sum_{k=1}^{r-1} (n_k+n_{k+1}-2)(a_{\hn_k}-\delta_k)
\end{align*}
We plug this back in to get
\begin{align*}
- \sum_{\ell=1}^r \sum_{j=1}^{n_\ell-1}B(c_j^{(\ell)})
 \leq \tfrac{n-1}{2} -
 \sum_{j=1}^{n-1} B(a_\ell) + \sum_{k=1}^{r-1} a_{\hn_k} + \frac12\sum_{k=1}^{r-1} (n_k+n_{k+1}-2)(a_{\hn_k}-\delta_k),
\end{align*}
from which it follows that the exponent of $T$ in the bound for $\big\lvert p_{T,R}^{(n)}(y;-a,\delta_C) \big\rvert$ above is
\begin{align*}
& \eps+r-1+d+\sum\limits_{\ell=1}^r\Big( R\big(D(n_\ell)+\tfrac{n_\ell(n_\ell-1)}{2}\big) + C(n_\ell)\Big) - \sum\limits_{k=1}^{n_\ell-1}B(c_k^{(\ell)})\Big)
\\ & \qquad
= \eps + d' + R\Big(D(n)+\tfrac{n(n-1)}{2}\Big) + C(n) - \sum_{j=1}^{n-1} B(a_j),
\end{align*}
where
\begin{align*}
 d' & = r-1 +
 d'' + \tfrac{n-1}{2} + \frac12\sum_{k=1}^{r-1} (n_k+n_{k+1}-2)(a_{\hn_k}-\delta_k) - C(n) + \sum_{\ell=1}^r C(n_\ell) + \sum_{k=1}^{r-1} a_{\hn_k}
 \\ & =
 d'' + \tfrac{n-1}{2} + \frac12\sum_{k=1}^{r-1} (n_k+n_{k+1}-2)(a_{\hn_k}-\delta_k) 
 + \sum_{k=1}^{r-1} a_{\hn_k}
 -\tfrac{1}{2}\sum_{1\leq k <m\leq r} n_k n_m,
\end{align*}
and
\begin{align*}
 d'' & = d - R\cdot \bigg( \hskip -3pt D(n) - \sum_{\ell=1}^{r} D(n_\ell) + \sum_{1\leq k < m \leq r} n_k n_m \hskip -3pt \bigg)
 & =
  - \sum_{k=1}^{r-1}\Big(\delta_k + (n_k+n_{k+1})(a_{\hn_k}-\delta_k)\Big).
\end{align*}
Hence,
\begin{align*}
 d' & = 
 \frac{n-1}{2} - \sum_{k=1}^{r-1}\Big(\delta_k + (n_k+n_{k+1})(a_{\hn_k}-\delta_k)\Big) + 
 \\ & \qquad
 + \frac12\sum_{k=1}^{r-1} (n_k+n_{k+1}-2)(a_{\hn_k}-\delta_k) 
 + \sum_{k=1}^{r-1} a_{\hn_k}
 -\frac{1}{2}\sum_{1\leq k <m\leq r} n_k n_m
 \\ & =
 \frac{n-1}{2} - \sum_{k=1}^{r-1}(n_k+n_{k+1})(a_{\hn_k}-\delta_k) 
 + \frac12\sum_{k=1}^{r-1} (n_k+n_{k+1})(a_{\hn_k}-\delta_k) 
 -\frac{1}{2}\sum_{1\leq k <m\leq r} n_k n_m
 \\ & =
 \frac12 \left( n-1  
 -\sum_{k=1}^{r-1} (n_k+n_{k+1})(a_{\hn_k}-\delta_k) 
 -\sum_{1\leq k <m\leq r} n_k n_m \right)
\end{align*}
Note that if $r=2$ and $n_1=m$ and $\delta_1=\delta$, then this expression becomes
 \[ \frac{n-1}{2} - \frac{n}{2}(a_m-\delta) - \frac{m(n-m)}{2}, \]
which agrees with \eqref{eq:requals2negative}.

Therefore, to complete the proof, we need only show that
 $n-1 - \sum\limits_{1\leq k<m \leq r} n_k n_m \leq 0.$
Indeed,
 \[ n-1-
 \sum_{1\leq k<m \leq r} n_k n_m 
 = n-1-
 \sum_{k=1}^{r-1}\sum_{m=k+1}^r n_k n_m 
 = n-1-
 \sum_{k=1}^{r-1} n_k (n-\hn_k) \leq n-1-n_1(n-n_1) \leq 0,
 \]
(with the final inequality being equality if and only if $n_1=1$ or $n_1=n-1$), as desired.
\end{proof}

\begin{remark}
A critical step in the proof of \eqref{eq:final-pTRna-deltaC-bound} (either in the case of single residues, as is proved in Section~\ref{sec:single-residues} or higher order residues, as in Section~\ref{sec:multi-residue-bound}) is to shift the lines of integration in the variables $\halpha_{m}$ or $\hbeta_j$.  A feature of this work that is quite different from the case of $\GL(4)$ as proved in \cite{GSW21}, is that no poles are crossed when making these shifts.  This represents a major simplification.  Recall from the discussion of Section~\ref{sec:GL4-example} that in the case of $n=4$ there are two fundamentally different types of single residues, two different types of double residues and a triple residue.  As it turned out, when making the additional shift for each of the single and double residues, one ends up with five \emph{additional} residue terms.  Taken all together, it was necessary to complete the analysis of writing down explicitly what the residues are in terms of  {gamma} function {s}, finding the exponential zero set, applying Stirling's formula and then obtaining a bound for ten(!) separate residues integrals.  All of this was in addition to performing these steps for the shifted $p_{T,R}^{(4)}$ term.
 
\end{remark}

\appendix

\section{Auxiliary results}

In an effort to avoid obstructing the flow of the argument in the main body of this paper, we will include here the many technical results that are used throughout.  {We remind the reader that the notational conventions that are used throughout the paper and this appendix are given in Section~\ref{sec:notation}.}

\begin{lemma}\label{lem:yIwasawa-of-yprime}
Suppose that $w= w_{(n_1,n_2,\ldots,n_r)}$ for some composition $n=n_1+\cdots +n_2$ with $r\geq 2$.  Then, if $y=(y_1,\ldots,y_{n-1})$,  $w\, y\, w^{-1}$ is equal to
\begin{multline*}
 \Bigg( \underbrace{y_{n-\hn_1+1}\, ,\, y_{n-\hn_1+2}\, ,\, \ldots\, ,\, y_{n-1}}_{\mbox{\scriptsize{$n_1-1$ terms}}}
 \, ,\,
 \left( \prod_{k=n-\hn_2}^{n-1} y_k \right)^{-1}
 \, , \, \ldots \\
 \ldots
 \, , \,
 \left( \prod_{k=n-\hn_i}^{n-\hn_{i-2}-1} y_k \right)^{-1}
 \, , \, 
 \underbrace{y_{n-\hn_i+1}\, ,\, y_{n-\hn_i+2}\, ,\, \ldots\, ,\, y_{n-\hn_{i-1}-1}}_{\mbox{\scriptsize{$n_i-1$ terms}}}
 \, , \,
 \left( \prod_{k=n-\hn_{i+1}}^{n-\hn_{i-1}-1} y_k \right)^{-1}
 \, , \, \ldots
 \\
 \ldots\, ,\, 
 \left( \prod_{k=1}^{n-\hn_{s-2}-1} y_k \right)^{-1}
 \, , \,
 \underbrace{y_{n-\hn_1+1}\, ,\, y_{n-\hn_1+2}\, ,\, \ldots\, ,\, y_{n-1}}_{\mbox{\scriptsize{$n_r-1$ terms}}}
 \Bigg)
\end{multline*}
In particular, 
 \[ \big\lVert wyw^{-1} \big\rVert ^{a_k} = \prod_{i=1}^r \prod_{j=1}^{n_i} y_{n-\hn_i+j}^{-a_{\hn_{i-1}}+a_{\hn_{i-1}+j}-a_{\hn_i}}.
 \]
\end{lemma}
\begin{proof}
Let $w=w_{(n_1,n_2,\ldots,n_r)}$ as above.  In order to carefully analyze $y' = w\, y\, w^{-1}$, we define 
 $x_i := \prod\limits_{j=1}^{i} y_j.$
This notation implies that $y = \diag(x_{n-1},x_{n-2},\ldots,x_1,1)$.
Now, let us think of the matrix $y$ as a block diagonal of the form $y = \diag(A_1,A_2,\ldots,A_r)$
where 
 \[ A_i = \diag( x_{n-\hn_{i-1}-1},x_{n-\hn_{i-1}-2},\ldots,x_{n-\hn_{i-1}-n_i} ) \in \GL(n_i,\R). \]
Thus,
 \[ y' = w\, y \, w^{-1} =\diag(A_r,A_{r-1},\ldots,A_1)
   = x_{n-n_1} \diag(B_r,B_{r-1},\ldots,B_1) . \]
Let $1\leq i \leq r$ and $0\leq j \leq n_i-1$ and set
 \[ z_{\hn_{i-1}+j} := \frac{x_{n-\hn_i+j}}{x_{n-n_1}}. \]
Then $(y'_1,y'_2,\ldots,y'_{n-1})$, the Iwasawa $y$-variables of $y'$, satisfy $y_i' = z_i/z_{i-1}$.  For $j\neq 0$, therefore, we see that
 \[ y'_{\hn_{i-1}+j} = \frac{x_{n-\hn_i+j}}{x_{n-\hn_i+j-1}} = \frac{\prod\limits_{k=1}^{n-\hn_i+j}y_k}{\prod\limits_{\ell=1}^{n-\hn_i+j-1}y_\ell} = y_{n-\hn_i+j}, \]
and for $j=0$,
 \[ y'_{\hn_{i}} = \frac{x_{n-\hn_{i+1}}}{x_{n-\hn_{i+1}-1}} = \frac{x_{n-\hn_{i+1}}}{x_{n-\hn_{i}+n_i-1}} = \frac{\prod\limits_{k=1}^{n-\hn_{i+1}}y_k}{\prod\limits_{\ell=1}^{n-\hn_{i-1}-1}y_\ell} = \left( \prod_{k=1}^{n_i+n_{i+1}-1} y_{n-\hn_{i+1}+k}\right)^{-1}, \]
from which the statement of the lemma follows directly.
\end{proof}

\begin{definition}
We say that $\alpha=(\alpha_1,\ldots,\alpha_n)\in\C^n$ is in \emph{$j$-general position} if the set 
 \[ \left\{\left. \sum_{k\in J} \alpha_k\, \right|\, J\subseteq \{1,\ldots,n\},\, \#J=j \right\} \]
consists of $\binom{n}{j}$ distinct elements.  We say that $\alpha$ is \emph{in general position} if it is in $j$-general position for each $j=1,\ldots,n-1$.
\end{definition}

\begin{lemma}\label{lem:simpleint}
Suppose that there exists $\eps>0$ such that for each $j=1,\ldots,n-1$, the real part of $s_j$ is bounded by at least $\eps$ from any integer.  Assume that $\alpha$ is in $j$-general position, $\re(\alpha_i)=0$ for each $i=1,\ldots,n-1$ {, and $r_j\in\Z_{\ge0}$}.  Assume that
 \[ \im(\al_1)\geq \im(\al_2) \geq \cdots \geq \im(\al_n), \]
and let $I_j=[-\im(\alpha_1+\cdots+\al_j),-\im(\alpha_n+\cdots+\alpha_{n-j+1})]$.  If $r_j\geq 2$, then
\begin{align*}
 \int\limits_{\substack{\re(s_j)=\sigma_j \\ \im(s_j)\in I_j}} \hskip -3pt
 \prod_{\substack{J\subseteq \{1,\ldots,n\} \\ \#J=j}} \hskip -3pt
 \Big| s_j+ \sum_{k\in J} \alpha_k \Big|^{-r_j}\, ds_j
 \ll \hskip -2pt \sum_{\substack{L\subseteq \{1,\ldots,n\} \\ \#L=j}} \prod_{\substack{K\subseteq \{1,\ldots,n\} \\ \#K=j \\ K\neq L}}\left(1+ \Big| \sum_{\ell\in L}\alpha_\ell - \sum_{k\in K}\alpha_k \Big|\right)^{-r_j}.
\end{align*}
If $r_j=1$ there is an extra power of $\eps$ in the exponent (in which case the implicit constant will depend on $\eps$), and if $r_j=0$, the integral is bounded by 
 \[ \Big(1 + \sum_{k=1}^{j} \alpha_k - \sum_{\ell=1}^j \alpha_{n+1-\ell} \Big). \]
\end{lemma}
\begin{remark}
The implicit $\ll$--constant depends on $\sigma_j$, but in applications this will always be bounded.
\end{remark}

\begin{proof}
The bound in the case of $r_j=0$ is obvious, so we may assume henceforth that $r_j\geq 1$.  Consider the set
 \[ \mathcal{A}_j := \left\{\left. \sum_{k\in J} \alpha_k\, \right|\, J\subseteq \{1,\ldots,n\},\, \#J=j \right\}. \]
For a fixed choice $\alpha$ in $j$-general position, let $A_1$ be the element of $\mathcal{A}_j$ that has the greatest imaginary part, $A_2$ the next greatest imaginary part and so on.  Hence $-\im(A_1)<-\im(A_2)<\cdots <-\im(A_{\binom{n}{j}})$.

Write $s_j=\sigma_j+it_j$. Note that $I_j=[-\im(A_1),-\im(A_{\binom{n}{j}})]$.  Upon applying Lemma~A.3 from \cite{GSW21}, one obtains the bound
\begin{align*}
 \int\limits_{I_j} 
 \prod_{\substack{J\subseteq \{1,\ldots,n\} \\ \#J=j}} 
 \Big| s_j+ \sum_{k\in J} \alpha_k \Big|^{-r_j}\, ds_j
 \ll 
 \Big(1+\im(A_1)-\im(A_{\binom{n}{j}})\big)^{\eps}
 \hskip -2pt \prod_{k=1}^{\binom{n}{j}-1} \Big(1+ \im(A_k-A_{k+1})\Big)^{-r_j}.
\end{align*}
This is one of the possible summands on the right hand side of the statement of the lemma.  Hence, regardless of the specific ordering which may arise for the given choice of $\alpha$, the claim follows.
\end{proof}

\begin{lemma}\label{lem:max-simplify}
Let $a\in \R$.  Then
 \[ \max\{0,2(\lceil a\rceil-a)-1\} - \lceil a \rceil \leq \begin{cases} -\lceil a\rceil & \mbox{ if }a\in  {(\lceil a \rceil-\frac12,\lceil a \rceil ]} \\ -\lfloor a\rfloor-2\big(a-\lfloor a\rfloor\big) & \mbox{ if }a\in(\lfloor a \rfloor,\lfloor a \rfloor+\frac12]. \end{cases} \]
\end{lemma}
\begin{proof}
First, let us assume that $a\in  {(\lceil a \rceil-\frac12,\lceil a \rceil]}$.  Then $\lceil a\rceil-a<\frac12$, hence
 \[ \max\{0,2(\lceil a\rceil-a)-1\} - \lceil a \rceil = -\lceil a \rceil. \]
On the other hand, assuming that $a\in (\lfloor a \rfloor,\lfloor a \rfloor+\frac12]$, we see that
\begin{align*}
\max\{0,2(\lceil a\rceil-a)-1\} - \lceil a \rceil 
& = \lceil a\rceil-2a-1
 = \lfloor a \rfloor - 2a
 = -\lfloor a \rfloor - 2\big( a-\lfloor a\rfloor \big),
\end{align*}
as claimed.
\end{proof}

\begin{lemma}\label{lem:residue-powerofT-new-aj}
Suppose that $a_1,\ldots,a_n\in \R_{>0}$.  Let
 \[ B(a) :=  \begin{cases} 0 & \mbox{ if }a<0 \\ \lfloor a \rfloor+2(a - \lfloor a \rfloor) & \mbox{ if } 0<
\lfloor a \rfloor+\varepsilon < a \leq \lfloor a \rfloor+\frac12, \\ \lceil a \rceil & \mbox{ if } \frac12 < \lceil a \rceil - \frac12\leq a <\lceil a \rceil-\varepsilon. \end{cases} \]
Then for any $\delta_m\in \Z_{\geq0}$ with $0<a_m-\delta_m$,
\begin{multline*}
\phantom{xxx} \sum_{j=1}^{m-1} B\left(a_j-\tfrac{j}{m}(a_m-\delta_m)\right) + \sum_{j=1}^{n-m-1} B\left(a_{m+j}-\tfrac{n-m-j}{n-m}(a_m-\delta_m)\right) 
 \\ 
 \geq
 \left(\sum_{j=1}^{n-1} B(a_j) \right) -\tfrac{n-2}{2}(a_m-\delta_m+1)-B(a_m).\phantom{xxx}
\end{multline*}
\end{lemma}
\begin{proof}
We consider first the case of $r-\frac12\leq a_j<r$ for  {some $r\in \Z$ and} all $j=1,2,\ldots,n-1$.  For any $a\in \R$, note that
\begin{equation}\label{eq:BaBounds}
 a\leq B(a) \leq a+\tfrac12,
\end{equation}
hence
\begin{align*}
 \sum_{j=1}^{m-1} B(a_j-\tfrac{j}{m}(a_m-\delta_m)) & \geq \left(\sum_{j=1}^{m-1} a_j\right)-\tfrac{m-1}{2}(a_m-\delta_m)
  \geq 
 \left(\sum_{j=1}^{m-1} B(a_j)-\tfrac12\right)-\tfrac{m-1}{2}(a_m-\delta_m) 
 \\ & =
 \left(\sum_{j=1}^{m-1} B(a_j)\right)-\tfrac{m-1}{2}(a_m-\delta_m+1).
\end{align*}
Combining this with the other terms (which  {are} easily shown to satisfy the analogous bound), the desired result is immediate.
\end{proof}

\begin{figure}[h]%
\begin{center}
\begin{tikzpicture}[scale=.65]
\draw[<->] (0,-1.5)--(0,4.5);
\draw[<->] (-1,0)--(4.5,0);
\draw[thick,red,dotted] (0,-1)--(-1/2,-1);
\foreach \x in {0,1,2,3}
    \draw[red,fill] (\x+1/2,\x+1) ellipse (2pt and 2pt);
\draw[red,fill] (-1/2,-1) ellipse (2pt and 2pt);
\foreach \x in {0,1,2,3}
    \draw[red,fill=white] (\x+1/2,\x) ellipse (2pt and 2pt);\draw[thick,blue,dotted] (-1/2,0)--(3.5,4);
\foreach \y in {0,1,2,3}
    \draw[very thick] (\y-1/2,\y)--(\y,\y)--(\y+1/2,\y+1);
\foreach \x in {1,2,3}
    \draw[thick,red,dotted] (\x-1/2,\x)--(\x+1/2,\x);
\draw[red,fill=white] (0,-1) ellipse (2pt and 2pt);
\foreach \y in {0,1,2,3}
    \draw[fill=white] (\y,\y) ellipse (2pt and 2pt);
\foreach \x in {0,1,2,3}
    \draw[blue,fill=white] (\x,\x+1/2) ellipse (2pt and 2pt);
\foreach \y in {1,2,3,4}
    \draw (1/8,\y) -- (-1/8,\y) node[anchor=east] {$\y$};
\draw (-1/8,-1) -- (1/8,-1) node[anchor=west] {$-1$};
\foreach \x in {1,2,3,4}
    \draw (\x,1/8) -- (\x,-1/8) node[anchor=north] {$\x$};
\end{tikzpicture}
\end{center}
\caption{Comparing graph of $B(x)$ (thick black) to $B_4(x)$ (dotted red) and $B_3(x)$ (dotted blue) bounds.}\label{fig:Bfunction}
\end{figure}
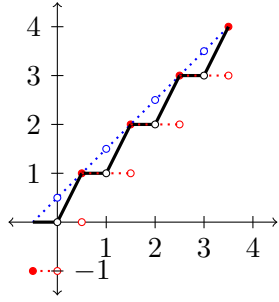

\begin{remark} 
The function $B(x)$ appears prominently in Theorem~\ref{th:pTRn-bound} and is critical in bounding the geometric side of the Kuznetsov trace formmula.  Its graph is shown in Figure~\ref{fig:Bfunction} in comparison to two other functions $B_4$ and $B_3$.

In the case of $\GL(4)$, the function $B_4$ appears \cite{GSW21} (see Theorem~4.0.1) as a bound for the $p_{T,R}$ function.  Indeed, making necessary adjustments due to a different choice of normalization factors (see Remark~\ref{rmk:GL3-GL4-compare}), the result of \cite{GSW21} is that
 \[ \big\lvert p_{T,R}^{(4)}(1;-a) \big\rvert \ll T^{\eps+27R+12-\sum\limits_{i=1}^3 B_4(a_i)}. \]
Theorem~\ref{th:pTRn-bound} establishes the same result but with $B_4$ replaced by $B$.  Although the improvement is slight, we remark that it is essential in Lemma~\ref{lem:residue-powerofT-new-aj} and evidently allows the inductive method of the present paper to lead to the same asymptotic orthogonality relation as was established directly in \cite{GSW21}.

With a bit of work, one can show that the function $B_3$, also graphed in Figure~\ref{fig:Bfunction}, appeared in \cite{GK2013} as a bound for
  \[ \big\lvert p_{T,R}^{(3)}(1;-a) \big\rvert \ll T^{\eps+6R+7-\sum\limits_{i=1}^2 B_3(a_i)}. \]
Although this looks to be an improvement on our result here, the method of \cite{GK2013} contained an error which the present method (and the method of \cite{GSW21}) corrects.

\end{remark}

\begin{lemma}\label{lem:aplusb-bound}
Let $\eps>0$.  Then for any $\rho\in \frac12 + \Z$ there exists $0<\eps'<\frac12$ sufficiently small such that, setting $\delta= \frac{2\eps'}{n^2}$, if $a=(a_1,\ldots,a_{n-1})$ where
 \[ a_j := \rho + \frac{j(n-j)}{2}\left(1+\delta\right), \]
and, for $w=w_{(n_1,\ldots,n_r)}$,  $b(a,w)=b=(b_1,\ldots,b_{n-1})$ where
 \[ b_{n-\hn_i+j} := a_{\hn_{i-1}}-a_{\hn_{i-1}+j}+a_{\hn_i} \pm \frac{\delta}{2}, \]
(meaning that $a$ and $b$ satisfy \eqref{eq:choice-of-a} and \eqref{eq:choice-of-b}, respectively), then, letting $B$ be the function defined in Theorem~\ref{th:ITRnaBound},
 \[ \sum_{j=1}^{n-1} \big( B(a_j)+B(b_j) \big) \geq \left \lfloor \frac{n-1}{2} \right\rfloor + n\rho + \Phi(n_1,\ldots,n_r) - \eps, \]
where
 \[ \Phi(n_1,\ldots,n_r) := \sum_{k=1}^{r-1} (n_k +n_{k+1}) \frac{(n-\hn_k)\hn_k}{2}. \]
\end{lemma}

\begin{proof}
We first note that although the bound $B(x)\geq x$ holds for any $x\in \R$, for any $\eps>0$, $B(x) \geq x+\frac12-\eps$ provided that $x$ is sufficiently close to a half integer.  Lemma~\ref{lem:even-odd} (as justified in Remark~\ref{rmk:AplusB}) asserts that if $n$ is  {odd} then $n-1$ elements from the set of all the possible values of $a_k$ and $b_k$ are indeed within $\eps$ of a half integer, and if $n$ is  {even} then $n-2$ of values have this property.  Hence,
\begin{equation}\label{eq:B-ak-plus-bk-bound}
 \sum_{k=1}^{n-1} \big( B(a_k) + B(b_k) \big) \geq  {\left\lfloor \frac{n-1}{2} \right\rfloor} + \sum_{k=1}^{n-1} ( a_k + b_k) - \eps.
\end{equation}
Since
 $b_{n-\hn_{i}+j} = a_{\hn_{i-1}}-a_{\hn_{i-1}+j}+a_{\hn_i} \pm \tfrac{\delta}{2},$
we see that 
 \[ \sum_{j=1}^{n_i} \big( b_{n-\hn_i+j} + a_{\hn_{i-1}+j} \big) \sim n_i\big( a_{\hn_{i-1}}+a_{\hn_i}\big). \]
Therefore, summing over $i$, we see (making use of the fact that $a_0=a_n=0$) that
\begin{align*}
 \sum_{k=1}^{n-1} \big( b_k+a_k \big) & 
 = \sum_{i=1}^r n_i\big( a_{\hn_{i-1}}+a_{\hn_i}\big) = \sum_{i=1}^{r-1}(n_{i}+n_{i+1})a_{\hn_i}
 \\ &
 = \sum_{k=1}^{r-1} (n_k+n_{k+1})\left( \rho + \frac{(n-\hn_k)\hn_k}{2} + \eps'\right)
 \\ &
 \sim \rho\big( 2n - n_1 - n_r \big) + \underbrace{\sum_{k=1}^{r-1} (n_k+n_{k+1}) \frac{(n-\hn_k)\hn_k}{2}}_{=: \Phi(n_1,\ldots,n_r)}.
\end{align*} 
Combining this with \eqref{eq:B-ak-plus-bk-bound}, the desired result is now immediate.
%
\end{proof}

\begin{lemma}\label{lem:even-odd}
Let $C=(n_1,\ldots,n_r)$ be a composition of $n$ with $r\geq 2$.  Suppose that $\rho \in \frac12 + \Z$.  Set $a_0:=0$, $a_n:=0$ and for each $1\leq k \leq n-1$ 
 we have $a_k := \rho + \frac{k(n-k)}{2}$ and
 for each $1\leq i \leq r$ and $1\leq j \leq n_i$ we let
 $b_{i,j} := a_{\hn_{i-1}} - a_{\hn_{i-1}+j} + a_{\hn_i}.$
 \vskip 4pt
Then
 \[ \# \big\{ k \mid a_k \notin \Z \big\} + \# \big\{ (i,j) \mid b_{i,j} \notin \Z \big\}
  = \begin{cases} 2n-n_1-n_r-1 & \mbox{ if $n$ is odd}, \\ \frac{n}{2} - 1 + \big\lfloor \frac{n_1}{2} \big\rfloor + \big\lfloor \frac{n_r}{2} \big\rfloor + \sum\limits_{i=2}^{r-1} \big\lceil \frac{n_i}{2} \big\rceil  & \mbox{ if $n$ is even}. \end{cases} \]
\end{lemma}
\begin{remark}\label{rmk:AplusB}
Note that the quantity given in Lemma~\ref{lem:even-odd} in the case of $n$ odd is $2n-n_1-n_r-1 \geq n-1$ for any composition $C$ (with equality precisely when $r=2$).  If $n$ is even then
 \[ \tfrac{n}{2} - 1 + \big\lfloor \tfrac{n_1}{2} \rfloor + \big\lfloor \tfrac{n_r}{2} \big\rfloor + \sum\limits_{i=2}^{r-1} \big\lceil \tfrac{n_i}{2} \big\rceil \geq \tfrac{n}{2} + \tfrac{n_1}{2} + \tfrac{n_r}{2} - 2 + \sum_{i=2}^{r-1} \tfrac{n_i}{2} = n-2. \]
Equality in this case occurs precisely when $n_1$ and $n_r$ are both odd and all other $n_i$ are even.
\end{remark}

\begin{proof}
For notational purposes, set
\begin{align*}
 A(n) & := \# \big\{ 1\leq k \leq n-1 \mid a_k \notin \Z \big\},
 \\ 
 B(C) & := \# \big\{ (i,j),\ 1\leq i \leq r,\ 1\leq j \leq n_i\ \mid\ b_{i,j}\notin \Z \big\}.
\end{align*}

We first consider the case of $n$ odd, for which $\frac{k(n-k)}{2}\in \Z$ for all integers $k$.  Therefore, $A(n) = n-1$.  As for $B(C)$, note that $b_{i,j}$ is equal to $\rho$ plus an integer as long as $i\neq 1,r$.  Otherwise, $b_{1,j},b_{r,j}\in \Z$.  Hence $B(C) = n-n_1-n_r$.

In the case of $n$ even, $\frac{k(n-k)}{2}\in \Z$  {exactly when} $k$ is even.  Hence $A(n) = \frac{n}{2}-1$.  To the end of finding $B(C)$, we introduce the notation
 \[ B_i(C) := \# \{ 1\leq j \leq n_i \mid b_{i,j} \notin \Z \}, \]
for which it is clear that $B(C) = \sum\limits_{i=1}^r B_i(C)$.
 
The cardinality of $B_i(C)$ depends, obviously, on the integrality of $b_{i,j}$.  To determine this, we first assume that $i=1$.  Then
 \[ b_{1,j} = -\frac{j(n-j)}{2} + \frac{n_1(n-n_1)}{2}. \]
Therefore (since $n$ is even), $b_{i,j}\in \Z$ if and only if $j\equiv n_1\pmod{2}$.  This implies that
 \[ B_1(C):= \begin{cases} 
   \frac{n_1-1}{2} & \mbox{ if $n_1$ is odd}, \\ 
   \frac{n_1}{2} & \mbox{ if $n_1$ is even}. \end{cases}, \]
or more concisely, $\# B_1(C) = \lfloor \frac{n_1}{2} \rfloor$.  The determination of $B_r(C)$ is similar: $\# B_r(C) = \lfloor \frac{n_r}{2} \rfloor$.

For $1<i<r$, we see that 
 {\begin{align*}
 b_{i,j} & = \rho + \frac{\hn_{i-1}(n-\hn_{i-1})}{2} - \frac{(\hn_{i-1}+j)(n-\hn_{i-1}-j)}{2} + \frac{(\hn_{i-1}+n_i)(n-\hn_{i-1}-n_i)}{2}
 \\ & = 
 \rho + \hn_{i-1}(n-\hn_{i-1}-n_i) - \frac{(\hn_{i-1}+j)(n-\hn_{i-1}-j)}{2} + \frac{n_i(n-n_i)}{2}
 \\ & \equiv 
 \frac12 + \frac{(\hn_{i-1}+j)(n-\hn_{i-1}-j)}{2} + \frac{n_i(n-n_i)}{2} \pmod{\Z}.
\end{align*}}
We see again that the integrality of $b_{i,j}$ depends on the parity of $n_i$.  If $n_i$ is odd, 
 \[ B_i(C) = \# \left\{ 1\leq j \leq n_i \left \lvert \frac{(\hn_{i-1}+j)(n-\hn_{i-1}-j)}{2} \notin \Z \right.\right\}, \]
and if $n_i$ is even, 
 \[ B_i(C) = \# \left\{ 1\leq j \leq n_i \left \lvert \frac{(\hn_{i-1}+j)(n-\hn_{i-1}-j)}{2} \in \Z \right.\right\}. \]
One can check, arguing case by case as above, that in any event, the answer is $B_i(C) = \lceil \frac{n_i}{2} \rceil$.
\end{proof}

\begin{lemma}\label{lem:Phi-properties}
Suppose that $(n_1,\ldots,n_r) \in \C^r$.  The function
 \[ \Phi(n_1,\ldots,n_r) := \sum_{k=1}^{r-1} (n_k+n_{k+1}) \frac{(n_1+\cdots + n_k)(n_{k+1}+\cdots + n_r)}{2} \]
is invariant under permutations, i.e., for any $\sigma \in S_r$, 
 we have $\Phi(n_1,\ldots,n_r) = \Phi(n_{\sigma(1)},\ldots,n_{\sigma(r)}).$
In particular, if $P=n_1+\cdots +n_r$ is a partition of $n$ then $\Phi(P):=\Phi(n_1,\ldots,n_r)$ is well defined. Moreover, among all partitions $P$ of $n$ (with $r\geq 2$), 
 \[ \Phi(P)\geq \Phi(n-1,1)=\Phi(1,n-1) = \frac{n(n-1)}{2}. \]
\end{lemma}
\begin{proof}
Suppose that $n=n_1+\cdots+n_r=m_1+\cdots + m_r$ where
 \[ m_j = \begin{cases} n_j & \mbox{ if } j\neq k, k+1, \\
 n_{k+1} & \mbox{ if }j=k, \\
 n_k & \mbox{ if } j = k+1. \end{cases} \]
Then one can show by an elementary (albeit tedious) computation that
$\Phi(n_1,\ldots,n_r) = \Phi(m_1,\ldots,m_r).$
In other words, $\Phi$ is invariant under any transposition $\tau \in S_r$, hence invariant under all  of $S_r$.

Suppose that $n=n_1+\cdots+n_r$.  If $n_k = n_k'+n_k''$ for some $1\leq k \leq r$, then one  {shows} via a straightforward computation that
 $\Phi(n_1,\ldots,n_{k-1}, n_k',n_k'',n_{k+1}, \ldots n_{r-1}) - \Phi(n_1,\ldots,n_r) = \frac{n_k n_k' n_k ''}{2}.$
If $n=n_1+\cdots + n_r$ with $r>2$, it then follows, setting $n_0 := \min\{n_1,n_2,\ldots,n_r \}$, that
 \[ \Phi(n_1,\ldots,n_r) > \Phi(n_0, n-n_0) = \frac{n n_0 (n-n_0)}{2}. \]
Among all $1\leq n_0 \leq \frac{n}{2}$, the right hand side is minimized when $n_0=1$.
\end{proof}

 {Recall that where $\Gamma_R(z):= \frac{\Gamma\left(\frac{\frac12+R+z}{2}\right)}{\Gamma(z)}$, as defined at the beginning of Section~\ref{sec:IntRepOfpTR}.}

\begin{lemma}\label{lem:nopoles}
If $\delta\in \Z$ and $\beta\in i\R$ are fixed, then the function 
 $\Gamma_R(\beta+z)\Gamma_R(-\beta-z)\Gamma(-\beta-z-\delta)$
is holomorphic for all $z$ with $\lvert \re(z) \rvert < R$.
\end{lemma}
\begin{proof}
The fact that $\lvert z\rvert <R$ implies that $\Gamma_R(\pm z)$ is holomorphic is immediate, so the only question is what happens at the (simple) poles of $\Gamma(-\beta-z-\delta)$.  But these occur at $z=-\beta+k$ for some integer $k$ which correspond to zeros of $\Gamma_R(\beta+z)$ or $\Gamma_R(-\beta-z)$.
\end{proof}

\begin{lemma}\label{lem:alphak-delta-residue-bound}
For $\delta\in \Z$ fixed and $z,\beta\in \C$ and $\lvert \re(z+\beta)+\delta\rvert < R$, we have the bound
 \[ \Gamma_R(\beta+z)\Gamma_R(-\beta-z)\Gamma(-\beta-z-\delta) \asymp \big( 1 + \lvert \im(\beta+z)\rvert\big)^{R-\re(\beta+z)-\delta}. \]
\end{lemma}
\begin{proof}
This follows immediately from the Stirling bound $\lvert\Gamma(\sigma+it)\rvert \sim \sqrt{2\pi}\lvert t \rvert^{\sigma-\frac12}e^{\pi\lvert t\rvert/2}$.
\end{proof}

\begin{definition} \label{partitioning-alpha} Let $\alpha=(\alpha_1, \ldots,\alpha_n) \in\mathbb C^n$ be Langlands parameters satisfying $\widehat{\alpha}_n = 0.$ Let $n=n_1+\cdots+n_r$ be a partition of $n$ with $n_1,\ldots,n_r\in
\Z_+.$ Then for each $\ell = 1,\ldots,r$ we define
$\alpha^{(\ell)} := \big( \alpha_1^{(\ell)},\ldots,\alpha_{n_\ell}^{(\ell)}\big) \in \C^{n_\ell}$ where
\[\alpha_j^{(\ell)} := \alpha_{\hn_{\ell-1}+j}-\tfrac{1}{n_\ell}\big(\halpha_{\hn_\ell}-\halpha_{\hn_{\ell-1}}\big), \qquad\quad \lvert \alpha^{(\ell)}\rvert^2:= \sum\limits_{j=1}^{n_\ell} \big(\alpha_j^{(\ell)}\big)^2. \]
\end{definition}
\begin{remark}
Note that $\sum\limits_{j=1}^{n_\ell} \alpha_j^{(\ell)}=0$ for each $\ell.$ In particular  $n_\ell =1$ implies $\alpha_1^{(\ell)}=0.$
\end{remark}

\begin{lemma}\label{lem:generalalphasum}
We have  $\lvert \alpha \rvert^2 = \sum\limits_{i=1}^n \alpha_i^2 = \sum\limits_{\ell=1}^r \left( \big\lvert \alpha^{(\ell)}\big\rvert^2 + \frac{1}{n_\ell} \big( \halpha_{\hn_\ell} - \halpha_{\hn_{\ell-1}} \big)^2 \right).$
\end{lemma}
\begin{proof}
Computing directly, and using the fact that $\sum\limits_{j=1}^{n_\ell} \alpha_j^{(\ell)}=0$, we find that
\begin{align*}
 \sum_{j=1}^n \alpha_j^2 &
 = \sum_{\ell=1}^r \sum_{j=1}^{n_\ell} \alpha_{\hn_{\ell-1}+j}^2
 = \sum_{\ell=1}^r \sum_{j=1}^{n_\ell} \big( \alpha_j^{(\ell)}+\tfrac{1}{n_\ell}\big(\halpha_{\hn_\ell}-\halpha_{\hn_{\ell-1}}\big) \big)^2
 \\ &
 = \sum_{\ell=1}^r \sum_{j=1}^{n_\ell} \left( \big(\alpha_j^{(\ell)}\big)^2+\tfrac{2}{n_\ell}\alpha_j^{(\ell)}\big(\halpha_{\hn_\ell}-\halpha_{\hn_{\ell-1}}\big) + \tfrac{1}{n_\ell^2}\big( \halpha_{\hn_\ell}-\halpha_{\hn_{\ell-1}} \big)^2 \right)
 \\ &
 = \sum_{\ell=1}^r \left( \big\lvert \alpha^{(\ell)} \big\rvert^2 + \tfrac{1}{n_\ell}\big( \halpha_{\hn_\ell}-\halpha_{\hn_{\ell-1}} \big)^2 \right),\end{align*}
as claimed.
\end{proof}

\begin{lemma}\label{lem:alphasumk}
Suppose that $n\geq 2$ and $\alpha_1,\alpha_2,\ldots,\alpha_n \in \C$ satisfies
$\alpha_1+\alpha_2+\cdots \alpha_n =0.$
Set $\halpha_k = \sum\limits_{j=1}^k \alpha_j$ for fixed $k\in\{1,2,\ldots,n\}$, and define 
 $\beta_j:= \al_j-\frac{1}{k}\halpha_k,
 \qquad 
 \gamma_j := \al_{j+k}+\frac{1}{n-k}\halpha_k.$
Then
 \[ \sum_{i=1}^n \al_i^2=\sum_{i=1}^k \beta_i^2+\sum_{i=1}^{n-k} \gamma_i^2+\frac{n}{k(n-k)}\halpha_k^2. \]
\end{lemma}
\begin{proof}
This is easily deduced as a special case of Lemma~\ref{lem:generalalphasum} in the case that $r=2$, $n_1=k$, $n_2=n-k$, $\beta=\alpha^{(1)}$ and $\gamma=\alpha^{(2)}$.
\end{proof}

\begin{lemma}\label{lem:genprodGammaRkdecomp}
We continue the notation of Lemma~\ref{lem:generalalphasum}.  Then
\begin{multline*}
 \prod_{1\leq i\neq j\leq n} \Gamma_R(\alpha_i-\alpha_j) = \prod_{\ell=1}^r \Bigg( \prod_{1\leq i,j\leq n_\ell} \Gamma_R\big(\alpha_i^{(\ell)}-\alpha_j^{(\ell)}\big)\Bigg)
 \\ \cdot
 \prod_{1\leq k < m \leq r}\prod_{i=1}^{n_k} \prod_{j=1}^{n_m} \prod_{\epsilon\in\{\pm1\}} \Gamma_R\left(\epsilon\bigg(\alpha_i^{(k)}-\alpha_j^{(m)}+\tfrac{1}{n_k}\big(\halpha_{\hn_k}-\halpha_{\hn_{k-1}}\big)-\tfrac{1}{n_m}\big(\halpha_{\hn_m}-\halpha_{\hn_{m-1}}\big)\bigg)\right).
\end{multline*}
\end{lemma}
\begin{proof}
Note that if $k\neq m$, then for any $1\leq i \leq n_k$ and $1\leq j\leq n_m$,
 \[  \alpha_{\hn_{k-1}+i}-\alpha_{\hn_{m-1}+j} = \alpha_i^{(k)}-\alpha_j^{(m)}+\tfrac{1}{n_k}\big(\halpha_{\hn_k}-\halpha_{\hn_{k-1}}\big)-\tfrac{1}{n_m}\big(\halpha_{\hn_m}-\halpha_{\hn_{m-1}}\big),\]
and for any $1\leq i\neq j\leq n_\ell$ we have
 $\alpha_{\hn_{\ell-1}+i}-\alpha_{\hn_{\ell-1}+j} = \alpha_i^{(\ell)}-\alpha_j^{(\ell)}.$
This immediately implies the desired formula.  
\end{proof}

\begin{lemma}\label{lem:extra-Gammasum}
Suppose that $(\beta_1,\ldots,\beta_r)$ satisfies  {$\hat{\beta}_r=0$}.  Suppose that $n=n_1+\cdots+n_r$ and set $\hn_k = \sum\limits_{j=1}^k n_j$.  Then
 $ \sum\limits_{1\leq k < m \leq r} \sum\limits_{i=1}^{n_k} \sum\limits_{j=1}^{n_m} \left( \frac{\beta_k}{n_k}-\frac{\beta_m}{n_m}\right) = \sum\limits_{j=1}^{r-1} (n_j+n_{j+1})\hbeta_j. $
\end{lemma}
\begin{proof}
We calculate
\begin{align*}
 \sum_{1\leq k < m \leq r} \sum_{i=1}^{n_k} \sum_{j=1}^{n_m} \left( \frac{\beta_k}{n_k}-\frac{\beta_m}{n_m}\right) 
 & =
 \sum_{m=2}^r \sum_{k=1}^{m-1} \sum_{i=1}^{n_k} \left( n_m \frac{\beta_k}{n_k}-\beta_m\right) 
  =
 \sum_{m=2}^r \sum_{k=1}^{m-1} \big( n_m\beta_k-n_k\beta_m\big) 
 \\ & = 
 \sum_{m=2}^r \big( n_m\hbeta_{m-1}-\hn_{m-1}\beta_m\big)
  = 
 \sum_{m=2}^r \big( n_m\hbeta_{m-1}-\hn_{m-1}(\hbeta_m-\hbeta_{m-1})\big)
 \\ & = 
 \sum_{m=2}^r \big( (n_m+\hn_{m-1})\hbeta_{m-1}-\hn_{m-1}\hbeta_m\big)
  = 
 \sum_{m=2}^r \big( \hn_m\hbeta_{m-1}-\hn_{m-1}\hbeta_m\big).
\end{align*}
This final sum telescopes to give $\sum\limits_{j=1}^{r-1} \big(\hn_{j+1}-\hn_{j-1}\big)\hbeta_{j}$.  Since $\hn_{j+1}-\hn_{j-1}=n_j+n_{j+1}$, this implies the claimed result.
\end{proof}

The following result can be interpreted as a consequence---by counting (half) the number of  {gamma} factors on each side of the equality---of Lemma~\ref{lem:genprodGammaRkdecomp}.  Alternatively, proving it independent of Lemma~\ref{lem:genprodGammaRkdecomp} gives further evidence that the product decomposition is correct.

\begin{lemma} \label{combinatorial} 
Let $n=n_1+\cdots+n_r$.  We have
  $\sum\limits_{1\le k <k'\le r} \hskip -5pt n_k\cdot n_{k'} \;+ \sum\limits_{k=1}^r \frac{n_k(n_k-1)}{2}  =  \frac{n(n-1)}{2}.$
\end{lemma}
\begin{proof} We use induction on $r$.  If $r=1$, the formula obviously holds.  Let $n=m+n_r$ where $m=n_1+\cdots+n_{r-1}$.  Then, by induction,
\begin{align*}
 \frac{n(n-1)}{2}& = \frac{(m+n_r)(m+n_r-1)}{2}
  =
 \frac{m(m-1)}{2} + \frac{mn_r + (m-1)n_r}{2} + \frac{n_r^2}{2}
 \\ & =
 \sum_{k=1}^{r-1} \frac{n_k(n_k-1)}{2} + \sum_{1\leq k<k'\leq r-1} n_k \cdot n_{k'} + mn_r + \frac{n_r(n_r-1)}{2}.
\end{align*}
Since $mn_r = n_1n_r+n_2n_r+\cdots n_{r-1}n_r$, it is evident that the desired formula holds.
 \end{proof}
 
\begin{lemma}\label{lem:combinatorial2}
Suppose $n=n_1+\cdots+n_r$.  Then
 $n^2 + \sum\limits_{\ell=1}^r\Big( \frac{n_\ell(n_\ell-1)}{2} - n_\ell \hn_\ell \Big) = \frac{n(n-1)}{2}. $
\end{lemma}
\begin{proof}
If $r=1$ the result is obviously true.  Suppose that the result holds for $r=k$.  Write $n=n_1+\cdots+n_k+n_{k+1}=\hn_k+n_{k+1}$.  Then
\begin{align*}
 n^2 + \sum_{\ell=1}^{k+1}\Big( \frac{n_\ell(n_\ell-1)}{2} - n_\ell \hn_\ell \Big) & =
 n^2 + \sum_{\ell=1}^k\Big( \frac{n_\ell(n_\ell-1)}{2} - n_\ell \hn_\ell \Big) + \frac{n_{k+1}(n_{k+1}-1)}{2} - n_{k+1} n\\
 & =
 n^2-\hn_k^2 + \left( \hn_k^2 + \sum_{\ell=1}^k\Big( \frac{n_\ell(n_\ell-1)}{2} - n_\ell \hn_\ell \Big) + \frac{n_{k+1}(n_{k+1}-1)}{2} - n_{k+1} n \right) \\
 & =
 n^2 - \hn_k^2 + \frac{\hn_k(\hn_k+1)}{2}  + \frac{n_{k+1}(n_{k+1}-1)}{2} - n_{k+1} n \\
 & =
 n^2 - \hn_k^2 + \frac{\hn_k(\hn_k+1)}{2}  + \frac{(n-\hn_k)(n-\hn_k-1)}{2} - (n-\hn_k),
\end{align*}
which can easily be shown now to simplify to $\frac{n(n-1)}{2}$, as claimed.
\end{proof}

\begin{remark}
Note that Lemma~\ref{combinatorial} and Lemma~\ref{lem:combinatorial2} are equivalent provided that
\begin{equation}\label{eq:combinatorial}
 n^2 - \sum_{\ell=1}^r n_\ell\hn_\ell = \sum_{1\leq k< k' \leq r} n_k n_{k'}.
\end{equation}
This can be verified by expanding the left hand side as follows:
\begin{align*}
 n^2 - \sum_{\ell=1}^r n_\ell\hn_\ell 
  =
 (n_1+\cdots+n_r)\hn_r - \sum_{\ell=1}^r n_\ell\hn_\ell
  = 
 \sum_{\ell=1}^r \big(n_\ell(n-\hn_\ell)-n_\ell\hn_\ell\big)
  =
 \sum_{\ell=1}^r n_\ell(n-n_\ell).
\end{align*}
That this final expression is equal to right hand side of \eqref{eq:combinatorial} is clear.
\end{remark}

\begin{lemma}\label{lem:prodGammaRkdecomp}
Let $\alpha_1,\alpha_2,\ldots,\alpha_n\in \C$ satisfy $\alpha_1+\alpha_2+\cdots\al_n=0$.  Set $\halpha_k:=\sum\limits_{j=1}^k\al_j$, and let $\beta_i$ ($1\leq i \leq k$) and $\gamma_j$ ($1\leq j\leq n-k$) be as in the previous lemma.  We have
\begin{multline*}
 \prod_{1\leq i\neq j\leq n} \Gamma_R(\alpha_i-\alpha_j) = \left(\prod_{1\leq i\neq j\leq k} \Gamma_R(\beta_i-\beta_j)\right) \left(\prod_{1\leq i\neq j\leq n-k} \Gamma_R(\gamma_i-\gamma_j)\right)
 \\ \cdot
 \prod_{i=1}^k \prod_{j=1}^{n-k} \Gamma_R\left(\beta_i-\gamma_j+\frac{n}{k(n-k)}\halpha_k\right)\Gamma_R\left(\gamma_j-\beta_i-\frac{n}{k(n-k)}\halpha_k\right)
\end{multline*}
\end{lemma}
\begin{proof}
This is easily deduced as a special case of Lemma~\ref{lem:genprodGammaRkdecomp} when $r=2$, $n_1=k$, $n_2=n-k$, $\beta=\alpha^{(1)}$ and $\gamma=\alpha^{(2)}$.
\end{proof}

 {We recall the definition of the polynomial given in Definition~\ref{def:FRn}:
 \[ \mathcal{F}_R^{(n)}(\alpha) := \prod_{j=1}^{n-2} \;\underset{\#K=\#L=j}{\prod_{K,L\,\subseteq\,(1,2,\ldots,n)}}\left(1+\sum_{k\in K}\alpha_k - \sum_{\ell\in L}\alpha_\ell     \right)^{\frac{R}{2}}. \]
Also, we remind the reader of the polynomial notation $\mathcal{P}$ given in Definition~\ref{def:polynomialnotation}.}

\begin{lemma}\label{lem:FRnkdecomp}
Let $n=n_1+\cdots+n_r$, $\alpha$, and $\alpha^{(\ell)}$   be as in Definition \ref{partitioning-alpha}.  Set $D(n)=\deg(\mathcal{F}_1^{(n)}(\alpha))$.  Then
\[
 \mathcal{F}_R^{(n)}(\alpha) = \mathcal{P}_{dR}(\alpha) \underset{n_\ell\ne 1}{\prod_{\ell=1}^r }\FR{n_\ell}(\alpha^{(\ell)})\quad \mbox{where} \quad
 d = d(n_1,\ldots,n_r) = D(n)- \underset{n_\ell\ne 1}{\sum_{\ell=1}^r} D(n_\ell). 
\]
\end{lemma}
\begin{proof}
This follows from the fact that if $I,J\subseteq\{ 1, 2,\ldots, n_\ell\}$ with $\#I=\#J$ then 
 \[ \left(\sum_{i\in I} \alpha_i^{(\ell)} \right) - \left(\sum_{j\in J} \alpha_j^{(\ell)} \right) = \left(\sum_{i\in I} \alpha_{\hn_{\ell-1}+i} \right) - \left(\sum_{j\in J} \alpha_{\hn_{\ell-1}+j} \right). \]
Therefore, each $\FR{n_\ell}(\alpha^{(\ell)})$ constitutes a unique factor of $\FR{n}(\alpha)$ for each $\ell=1,\ldots,r$. 
\end{proof}

\begin{lemma}\label{lem:resm_extrapoly}
Suppose that $\delta\in \Z_{\geq 0}$ and $R>\delta$.  Then
 \[ \FR{n}(\alpha) \cdot \hskip-17pt\prod_{\substack{K\subseteq \{1,2,\ldots,n\} \\ \#(K\cap \{1,2,\ldots,m\})\neq m-1 \\ \#K=m}} \left( \Big(\sum_{i\in K} \alpha_i \Big) - \halpha_m-\delta\right)_\delta^{-1}
 \ll
 \FR{m}(\beta) \cdot \FR{n-m}(\gamma) \cdot \mathcal{P}_d(\alpha), \]
where $d = R\big(D(n)-D(m)-D(n-m)\big)-\delta\big( \binom{n}{m}-m(n-m)-1\big)$.
\end{lemma}
\begin{proof}
Let $M:=\{1,2,\ldots,m\}$.  Then
\[ \# \big\{ K\subseteq \{1,2,\ldots,n\} \, \big\vert \, \#K = m,\, \#(K\cap M) \neq 0,1\big\} = \binom{n}{m} - m(n-m) - 1. \]
From the definition of $\FR{n}(\alpha)$ given in Definition~\ref{def:FRn}, we see that for each such $K$ there are factors
 \[ \left(1 + \sum_{i\in K} \alpha_i - \sum_{j\in M} \alpha_j \right)^{R/2}\left(1 - \sum_{i\in K} \alpha_i + \sum_{j\in M} \alpha_j \right)^{R/2} \]
of $\FR{n}(\alpha)/[\FR{m}(\beta)\FR{n-m}(\gamma)]$ for which
 \[ \frac{ \left( 1 + \sum\limits_{i\in K} \alpha_i - \sum\limits_{j\in M} \alpha_j \right)^{R/2} \left( 1 - \sum\limits_{i\in K} \alpha_i + \sum\limits_{j\in M} \alpha_j \right)^{R/2} }{ \left( \sum\limits_{i\in K} \alpha_i - \sum\limits_{j\in M} \alpha_j - \delta \right)_\delta } \ll \left( 1 - \Big(\sum\limits_{i\in K} \alpha_i - \sum\limits_{j\in M} \alpha_j \Big)^2 \right)^{\frac{R-\delta}{2}}. \]
This bound holds because the degree of the Pochhammer symbol in the denominator is $\delta$, and by assumption, the degree of the numerator is $R>\delta$.  Combining all such terms with the remaining factors of $\FR{n}(\alpha)/[\FR{m}(\beta)\FR{n-m}(\gamma)]$, gives a polynomial of degree $d$.
\end{proof}

\begin{remark}\label{rmk:gen_res_extrapoly}
Let $n=n_1+n_2+\cdots+n_r$ and $\hn_\ell = \sum\limits_{i=1}^\ell n_i$. The result of Lemma~\ref{lem:resm_extrapoly} clearly generalizes to the case of taking multiple residues at $s_{\hn_\ell}=-\halpha_{\hn_\ell}-\delta_\ell$ for each $\ell=1,\ldots,r-1$ (in reverse order).  In this case, taking the product on the left hand side over all of the terms we obtain
 \[ \mathcal{F}_R^{(n)}(\alpha) \cdot \prod_{\ell=1}^{r-1} \prod_{\substack{K\subseteq \{1,2,\ldots,\hn_{\ell+1}\} \\ \#(K\cap \{1,\ldots,\hn_\ell\})\neq \hn_\ell-1 \\ \#K=\hn_\ell}} \left( \Big(\sum_{i\in K} \alpha_i \Big) - \halpha_{\hn_\ell}-\delta_\ell\right)_{\delta_\ell}^{-1} \ll \mathcal{P}_d(\alpha) \cdot \prod_{\ell=1}^{r} \mathcal{F}_R^{(n_\ell)}\big(\alpha^{(\ell)}\big), \]
where 
 \[ d=R\cdot \left( D(n) - \sum_{\ell=1}^r D(n_\ell)\right) - \sum\limits_{\ell=1}^{r-1} \left[ \delta_\ell\left( \binom{\hn_{\ell+1}}{\hn_\ell} - n_{\ell+1} \hn_\ell - 1\right)\right]. \]
\end{remark}

\section*{Acknowledgments}

Eric Stade  would like to thank Taku Ishii for many helpful converations, and for the ideas constituting the proof of Conjecture \ref{conj} in the case $n=5$. Michael Woodbury would like to thank the University of Colorado for wonderful accommodations while hosting him during the Spring 2022 semester.  {We would also like to thank the referees for many helpful comments.}

\vskip 24pt
\bibliographystyle{amsalpha}
\bibliography{biblio}

\end{document}